\documentclass[11pt]{article}
\usepackage{amsmath, amsfonts}

\numberwithin{equation}{section}

\def\B1{B_{1/2}}
\def\Bli{{\lambda_i}}
\def\Bl1{{\lambda_1}}
\def\Box{\hfill\rule{2.5mm}{2.5mm}}

\def\D{{\cal D}}
\def\DD{{G}}
\def\F{{\cal F}}
\def\H{{\cal H}}

\def\Int{\int_{-1}^1}
\def\Jac{\mathop{\rm Jac}}
\def\L{{\cal L}}
\def\LL{{\hat L}}
\def\N{{\mathbb {N}}}

\def\R{{\mathbb {R}}}
\def\RR{{\cal R}}
\def\SS{{\cal S}}
\def\T{{\cal T}}

\def\argth{\mathop {\rm argth}}

\def\bk{{\bar \kappa_0}}
\def\bl1{{\bar \lambda_1}}
\def\bli{{\bar \lambda_i}}

\def\build#1_#2^#3{\mathrel{
\mathop{\kern 0pt#1}\limits_{#2}^{#3}}}

\def\cd{{\bar d}}
\def\chr{{R}}
\def\cs{{S}}
\def\ct{{t}}

\def\cy{{Y}}
\def\d{\displaystyle}
\def\dd{{d_0}}

\def\diag{\mathop{\rm diag}}

\def\ep{\epsilon}
\def\fe{\varphi_{d,\nu,\epsilon}}
\def\fee{\varphi_{d,\nu,\epsilon_0}}
\def\ff{F^{d,\nu}_{1,1}}

\def\h1{\mathop{\rm H^1_{\rm loc,\rm u}}}
\def\hk{{\hat k}}
\def\htsa{{S_{\hat k}}}
\def\czeta{{\bar \zeta}}

\def\iint{\displaystyle\int_{-1}^1}

\def\k{{m}}
\def\l2{\mathop{\rm L^2_{\rm loc,\rm u}}}

\def\ll{{\hat L}}
\def\lzero{{L_{\k+1}}}
\def\m{{A_-}}

\def\p{{\bar p}}
\def\pd{\partial_d}
\def\pnu{\partial_\nu}
\def\ps{\partial_s}
\def\pS{\partial_S}
\def\py{\partial_y}
\def\pY{\partial_Y}

\def\r{{\bf r}}
\def\sgn{\mathop{\rm sgn}}
\def\span{\mathop{\rm span}}

\def\tB{{b}}

\def\tmu{{\bar \nu}}

\def\tr{{r}}

\def\ts{{s}}
\def\tsa{{s_{\hat k}}}

\def\tT{{{\cal T}}}

\def\tx{{x}}

\def\ty{{y}}

\newcommand{\aref}[1]{(\ref{#1})}

\newcommand{\vc}[2]{
\left(
\begin{array}{l}
#1\\
#2
\end{array}
\right)
}

\title{\bf Isolatedness of characteristic points at blow-up for a semilinear wave equation in one space dimension
\footnote{Both authors are supported by a grant from the french Agence Nationale de la Recherche, project ONDENONLIN, reference ANR-06-BLAN-0185.} }
\author{Frank Merle\\
{\it \small Universit\'e de Cergy Pontoise and IHES}\\
Hatem Zaag\\
{\it \small CNRS LAGA Universit\'e Paris 13}}
\date{October 4, 2010}

\newtheorem{cor}{Corollary}[section]
\newtheorem{definition}[cor]{Definition}
\newtheorem{cl}[cor]{Claim}
\newtheorem{lem}[cor]{Lemma}
\newtheorem{prop}[cor]{Proposition}

\newtheorem{propo}{Proposition}
\newtheorem{theor}[propo]{Theorem}

\begin{document}

\maketitle

{\small {\bf Abstract}: We consider the semilinear wave equation with power nonlinearity in one space dimension. 
We consider an arbitrary blow-up solution $u(x,t)$, the graph $x\mapsto T(x)$ of its blow-up points and $\SS\subset \R $ the set of all characteristic points. We show that $\SS$ is locally finite. 
}
\medskip



\medskip


\section{Introduction}
\subsection{Main results}
We consider the one dimensional semilinear wave equation
\begin{equation}\label{equ}
\left\{
\begin{array}{l}
\partial^2_{t} u =\partial^2_{x} u+|u|^{p-1}u,\\
u(0)=u_0\mbox{ and }u_t(0)=u_1,
\end{array}
\right.
\end{equation}
where $u(t):x\in\R \rightarrow u(x,t)\in\R$, $p>1$, $u_0\in \rm H^1_{\rm loc,u}$
and $u_1\in \rm L^2_{\rm loc,u}$ with
 $\|v\|_{\rm L^2_{\rm loc,u}}^2=\d\sup_{a\in \R}\int_{|x-a|<1}|v(x)|^2dx$ and $\|v\|_{\rm H^1_{\rm loc,u}}^2 = \|v\|_{\rm L^2_{\rm loc,u}}^2+\|\nabla v\|_{\rm L^2_{\rm loc,u}}^2$.\\
We solve equation \aref{equ} locally in time in the space ${\rm H}^1_{\rm loc,u}\times {\rm L}^2_{\rm loc, u}$ (see Ginibre, Soffer and Velo \cite{GSVjfa92}, Lindblad and Sogge \cite{LSjfa95}).
For the existence of blow-up solutions, we have the following blow-up criterion from Levine \cite{Ltams74}: If $(u_0,u_1)\in H^1\times L^2(\R)$ satisfies
\[
\int_{\R}\left(\frac 12 |u_1(x)|^2+\frac 12|\partial_x u_0(x)|^2-\frac 1{p+1}|u_0(x)|^{p+1}\right)dx<0,
\]
then the solution of \aref{equ} cannot be global in time. More blow-up results can be found in
Caffarelli and Friedman \cite{CFtams86}, \cite{CFarma85},
Alinhac \cite{Apndeta95}, \cite{Afle02} and Kichenassamy and Littman \cite{KL1cpde93}, \cite{KL2cpde93}.

\bigskip

If $u$ is a blow-up solution of \aref{equ}, we define (see for example Alinhac \cite{Apndeta95}) a 1-Lipschitz curve $\Gamma=\{(x,T(x))\}$ 
such that the maximal influence domain $D$ of $u$ (or the domain of definition of $u$) is written as 
\begin{equation}\label{defdu}
D=\{(x,t)\;|\; t< T(x)\}.
\end{equation}
$\bar T=\inf_{x\in \R}T(x)$ and $\Gamma$ are called the blow-up time and the blow-up graph of $u$. 
A point $x_0$ is a non characteristic point 
if
\begin{equation}\label{nonchar}
\mbox{there are }\delta_0\in(0,1)\mbox{ and }t_0<T(x_0)\mbox{ such that }
u\;\;\mbox{is defined on }{\cal C}_{x_0, T(x_0), \delta_0}\cap \{t\ge t_0\}
\end{equation}
where ${\cal C}_{\bar x, \bar t, \bar \delta}=\{(x,t)\;|\; t< \bar t-\bar \delta|x-\bar x|\}$. We denote by $\RR$ (resp. $\SS$) the set of non characteristic (resp. characteristic) points.

\bigskip

In Proposition 1 in \cite{MZajm10}, we proved the {\it existence of solutions} of \aref{equ} such that
\[
\SS\neq \emptyset.
\]
 For general blow-up solutions, we proved the following facts about $\RR$ and $\SS$ in \cite{MZcmp08} and \cite{MZajm10} (see  Theorem 1 (and the following remark) in \cite{MZcmp08}, see Propositions 5 and 8 in \cite{MZajm10}):

\bigskip

{\label{old}\it (i) $\RR$ is a non empty open set, and $x\mapsto T(x)$ is of class $C^1$ on $\RR$;

(ii) $\SS$ is a closed set with empty interior, and given $x_0\in\SS$, if $0<|x-x_0|\le \delta_0$, then
\begin{equation}\label{chapeau0}
0< T(x)- T(x_0)+|x-x_0| \le \frac{C_0|x-x_0|}{|\log(x-x_0)|^{\frac{(k(x_0)-1)(p-1)}2}}
\end{equation}
for some $\delta_0>0$ and $C_0>0$, where $k(x_0)\ge 2$ is an integer. Moreover, $x\mapsto T(x)$ is left-differentiable and right-differentiable, with $T_l'(x_0)=1$ and $T_r'(x_0)=-1$.
}

\bigskip

Following our earlier work \cite{MZajm10}, we aim in this paper at proving that $\SS$ is made of isolated points. We have the following theorem, which is the main result of our analysis:
\begin{theor}\label{thfini} {\bf ($\SS$ is made of isolated points)} Consider $u(x,t)$ a blow-up solution of \aref{equ}. The set of characteristic points $\SS$ is made of isolated points.
\end{theor}
{\bf Remark}: The fact that the elements of $\SS$ are isolated points is not elementary. Direct arguments give no more than the fact that $\SS\neq \R$ (a point $x_0$ such that $T(x_0)$ is the blow-up time is non characteristic).
The first step of the proof has been done in \cite{MZajm10} where we proved that $\SS$ has an empty interior and that in similarity variables, the solution splits in a non trivial decoupled sum of (at least 2) solitons with alternate signs (see Proposition \ref{thsing} below for a statement). 
The second step is the heart of the present paper. It consists in using this 
decomposition and 
a good understanding of the dynamics of the equation in similarity variables (see equation \aref{eqw} below) near a decoupled sum of ``generalized'' solitons. More details on the strategy of the proof can be found in Section \ref{strategy} below. In fact, this is the first time where flows near an unstable sum of solitons are used and where such a result is obtained.\\
{\bf Remark}: In higher dimensions $N\ge 2$, our result would be that the $(N-1)$-dimensional Hausdorff measure of $S$ is bounded. To prove that, we strongly need to characterize all selfsimilar solutions of equation \aref{equ} in the energy space. This is the main obstruction to extend this result to higher dimensions.\\
Similar results are already available in the context of blow-up solutions for the semilinear heat equation with a subcritical power nonlinearity. They are due to Vel\'azquez who proved in \cite{Viumj93} that the $(N-1)$-Hausdorff dimension of the blow-up set is bounded. His proof is based on energy arguments and the ``infinite speed of propagation''. 
Note that in that parabolic case, there is only one value for the local energy limit at blow-up points, and this value is higher than the local energy limit at non blow-up points. Therefore, in order to conclude in \cite{Viumj93}, it was enough to compute directly the local energy at nearby points in well chosen directions.
As we said earlier, the local energy limit at blow-up has only one value, whether the blow-up behavior is stable or unstable. Only the speed of convergence to the profile is different between the stable and the unstable behaviors.
%
%
%
Here, with equation \aref{equ} which is hyperbolic, the argument cannot be as direct. This is due to two reasons:\\
- the unstable behavior (near characteristic points) and the stable behavior (near non characteristic points) have different energy levels (see Propositions \ref{threg} and \ref{thsing} below);\\
- we need information outside the light cone, which is beyond the reach of the finite speed of propagation. See Section \ref{strategy} for more details on the proof.\\ 
%
%
{\bf Remark}: The fact that $\SS$ is made of isolated points certainly does not hold in general for quasilinear wave equations. Indeed, in \cite{Anum06}, Alinhac gives an explicit solution $u(x,t)$ for the following nonlinear wave equation
\[
\partial^2_{t} u =\partial^2_{x} u+\partial_x u\partial_t u,
\]
whose domain of definition is
\[
D=\R\times [0,\infty) \backslash \{(x,t)\;|\; t\ge 1,\;|x|\le t-1\}
\]
(when $0\le t<1$, $u(x,t) = 4 \arctan\left(\frac x{1-t}\right)$).
In this example, we clearly see that $\RR =\{0\}$, $\SS = \R^*$, and the boundary of $D$ is characteristic (i.e. has slope $\pm 1$) on $\SS$. 
%

\bigskip

We also have the following estimate of the blow-up set, which improves \aref{chapeau0} proved in \cite{MZajm10}: 
\begin{theor}\label{pdes}{\bf (Description of $T(x)$ for $x$ near $x_0$)} If $x_0\in \SS$ and $0<|x-x_0|\le \delta_0$, then 
\begin{equation*}
\d\frac{1}{C_0|\log(x-x_0)|^{\frac{(k(x_0)-1)(p-1)}2}}\le T'(x)+\frac{x-x_0}{|x-x_0|} \le \frac{C_0}{|\log(x-x_0)|^{\frac{(k(x_0)-1)(p-1)}2}}
\end{equation*}
for some $\delta_0>0$ and $C_0>0$, where $k(x_0)\ge 2$ is an integer.
\end{theor}
{\bf Remark}: Integrating the estimate of Theorem \ref{pdes}, we see that
\begin{equation}\label{tx}
\d\frac{|x-x_0|}{C_0|\log(x-x_0)|^{\frac{(k(x_0)-1)(p-1)}2}}\le T(x)- T(x_0)+|x-x_0| \le \frac{C_0|x-x_0|}{|\log(x-x_0)|^{\frac{(k(x_0)-1)(p-1)}2}}.
\end{equation}
 The upper bound in this estimate was already known in \cite{MZajm10}. On the contrary, the lower bound (which is of the same size) is a new contribution in this work. Moreover, both the lower and the upper bounds on $T'(x)$ are new. There is also a dynamical version for this result, where we find the convergence of the similarity variables version $w_x(y,s)$ to $\kappa(T'(x),\cdot)$ (respectively defined in \aref{defw} and \aref{defkd} below), with a speed controlled in terms of $x$ near $x_0$ (see the proof of Theorem 2 in Section \ref{secfini} for more details). This also provides the behavior (or the ``picture'') of $u(x,t)$ for $(x,t)$ close to $(x_0, T(x_0))$, 
even outside
the backward light cone of vertex  $(x_0, T(x_0))$. For a similar result where the picture of the solution is derived near the singularity for the semilinear heat equation, see Herrero and Vel\'azquez \cite{HVcpde92}, \cite{HVihp93}  and Fermanian, Merle and Zaag \cite{FKMZma00} for the case of an isolated blow-up point, see also Vel\'azquez \cite{Vcpde92}, \cite{Vtams93}, Zaag \cite{Zihp02}, \cite{Zcmp02}, \cite{Zdmj06}, \cite{Zmme02} and Khenissy, Rebai and Zaag \cite{KRZcont10} for the case of non isolated blow-up points.\\
{\bf Remark}: From the shape of the solution near $(x_0, T(x_0))$, we can recover the topology of the solution. Indeed, the shape gives the integer $k(x_0)\ge 2$, and $k(x_0)-1$ is the number of sign changes of the solution near $(x_0, T(x_0))$ as shown in Proposition \ref{thsing} below.

\bigskip

As a consequence of our analysis, particularly the lower bound on $T(x)$ in \aref{tx}, we have the following corollary on the blow-up speed in the backward light cone with vertex $(x_0, T(x_0))$ where $x_0$ is a characteristic point.
\begin{cor}[Blow-up speed at a characteristic point]\label{corspeed}For all $x_0\in\SS$ and $t\in[0, T(x_0))$, we have
\[
\frac{|\log(T(x_0)-t)|^{\frac{k(x_0)-1}2}}{C(T(x_0)-t)^{\frac 2{p-1}}}\le \sup_{|x-x_0|<T(x_0)-t}|u(x,t)|\le \frac{C |\log(T(x_0)-t)|^{\frac{k(x_0)-1}2}}{(T(x_0)-t)^{\frac 2{p-1}}}.
\]
\end{cor}
{\bf Remark}: Note that by definition of the maximal influence domain (or the domain of definition) $D$ of $u$ in \aref{defdu}, we don't even know whether the solution blows up on the boundary $\Gamma$ of $D$. In other words, given $x_0\in\R$, we don't know whether $\sup_{|x-x_0|<T(x_0)-t}|u(x,t)|\to \infty$ as $t\to T(x_0)$. When $x_0\in\RR$, this fact follows from our convergence result to the profile proved in \cite{MZjfa07} (see Corollary 4 page 49 there). When $x_0\in\SS$, the same fact holds, as one may derive from the corollary above. However, this latter case needs much more work (involving 3 papers \cite{MZjfa07}, \cite{MZajm10} and the present one), based on the decomposition into a decoupled sum of solitons (see Proposition \ref{thsing} below for a statement).
Now, if we refine this result, we will see that there is a difference in the size of the $L^\infty$ norm, whether $x_0\in\RR$ or $x_0\in\SS$. Indeed, 
 Since $k(x_0)\ge 2$ when $x_0\in\SS$, it is clear from this corollary that the blow-up speed in $L^\infty$ is higher than the case when $x_0\in\RR$ where the rate is a given by the associated ODE:
\[
\frac{(T(x_0)-t)^{-\frac 2{p-1}}}C \le \sup_{|x-x_0|<T(x_0)-t}|u(x,t)|\le C(T(x_0)-t)^{-\frac 2{p-1}}
\]
(see Corollary 4 page 49 in \cite{MZjfa07} for an argument). Note that the same analysis can be carried out with the $H^1$ norm (averaged on the section of the light cone) instead of the $L^\infty$, leading to the same results.

\subsection{Formulation in similarity variables}
Let us recall in the following the similarity variables and some energy estimates from our earlier work.

\bigskip

 Given some $x_0\in\R$, a natural tool is to introduce the following change of variables:
\begin{equation}\label{defw}
w_{x_0}(y,s)=(T(x_0)-t)^{\frac 2{p-1}}u(x,t),\;\;y=\frac{x-x_0}{T(x_0)-t},\;\;
s=-\log(T(x_0)-t).
\end{equation}
The function $w=w_{x_0}$ satisfies the 
following equation for all $y\in (-1,1)$ and $s\ge -\log T(x_0)$:
\begin{equation}\label{eqw}
\partial^2_{s}w= \L w-\frac{2(p+1)}{(p-1)^2}w+|w|^{p-1}w
-\frac{p+3}{p-1}\partial_sw-2y\partial^2_{y,s} w
\end{equation} 
\begin{equation}\label{defro}
\mbox{where }\L w = \frac 1\rho \py \left(\rho(1-y^2) \py w\right)\mbox{ and }
\rho(y)=(1-y^2)^{\frac 2{p-1}}.
\end{equation}
Introducing $V=(V_1, V_2)=(w, \partial_s w)$, we rewrite \aref{eqw} as a first order equation
\begin{equation}\label{eqw1}
\frac \partial {\partial s}\vc{V_1}{V_2}=\vc{V_2}{\L V_1-\frac{2(p+1)}{(p-1)^2}V_1+|V_1|^{p-1}V_1
-\frac{p+3}{p-1}V_2-2y\partial_{y} V_2}.
\end{equation}
Introducing 
\begin{eqnarray}
\H& =& \left\{(q_1,q_2) 
\;\;|\;\;\|(q_1,q_2)\|_{\H}^2\equiv \int_{-1}^1 \left(q_1^2+\left(q_1'\right)^2  (1-y^2)+q_2^2\right)\rho dy<+\infty\right\},\label{defnh0}\\
\H_0& =& \{ q_1\;|\; \|q_1\|_{\H_0}^2\equiv\int_{-1}^1 \left(q_1^2+\left(q_1'\right)^2  (1-y^2)\right) \rho dy < + \infty\},\label{h0}
\end{eqnarray}
we see that the Lyapunov functional for equation \aref{eqw1}
\begin{equation}\label{defenergy}
E(V)= \iint \left(\frac 12 V_2^2 + \frac 12  \left(\partial_y V_1\right)^2 (1-y^2)+\frac{(p+1)}{(p-1)^2}V_1^2 - \frac 1{p+1} |V_1|^{p+1}\right)\rho dy
\end{equation}
is defined for $V=(V_1,V_2) \in \H$.
When there is no ambiguity, we may write $E(w)$ instead of $E(w,\ps w)$.

\bigskip

Let us now introduce for all $|d|<1$ the following stationary solutions of \aref{eqw} (or solitons) defined by 
\begin{equation}\label{defkd}
\kappa(d,y)=\kappa_0 \frac{(1-d^2)^{\frac 1{p-1}}}{(1+dy)^{\frac 2{p-1}}}\mbox{ where }\kappa_0 = \left(\frac{2(p+1)}{(p-1)^2}\right)^{\frac 1{p-1}} \mbox{ and }|y|<1.
\end{equation}
Furthermore, we introduce the following family of explicit solutions of equation \aref{eqw} (related to $\kappa(d,y)$ \aref{defkd} by the selfsimilar transformation \aref{defw}) defined by 
\[
\kappa_1^*(d,\mu e^s,y)= \d\kappa_0\frac{(1-d^2)^{\frac 1{p-1}}}{(1+dy+\mu e^s)^{\frac 2{p-1}}}.
\]
Note that $\kappa_1^*(d,\mu e^s,y)\to \kappa(d)$ in $\H$ as $s\to -\infty$. Moreover ,\\
- when $\mu=0$, we recover the stationary solutions $\kappa(d)$ defined in \aref{defkd};\\
- when $\mu>0$, the solution exists for all $(y,s) \in (-1,1)\times \R$ and converges to $0$ in $\H$ as $s\to \infty$ (it is a heteroclinic connection between $\kappa(d)$ and $0$);\\
- when $\mu<0$, the solution exists for all 
$(y,s) \in (-1,1)\times \left(-\infty, \log\left(\frac {|d|-1}\mu\right)\right)$ and blows up at time $s=\log\left(\frac {|d|-1}\mu\right)$.\\
 We also introduce for all $d\in (-1,1)$ and $\nu >-1+|d|$, $\kappa^*(d,\nu,y) = (\kappa_1^*, \kappa_2^*)(d,\nu,y)$ where 
\begin{equation}\label{defk*}
\kappa_1^*(d,\nu, y) =
\d\kappa_0\frac{(1-d^2)^{\frac 1{p-1}}}{(1+dy+\nu)^{\frac 2{p-1}}},\;\;
\kappa_2^*(d,\nu, y) = \nu \pnu \kappa_1^*(d,\nu, y) =
\d-\frac{2\kappa_0\nu}{p-1}\frac{(1-d^2)^{\frac 1{p-1}}}{(1+dy+\nu)^{\frac {p+1}{p-1}}}.
\end{equation}
In our paper, we refer to these functions as ``generalized solitons'' or solitons for short.
Note that  
for any $\mu\in\R$, 
$\kappa^*(d,\mu e^s,y)$ is a solution to equation \aref{eqw1}.
In Lemma \ref{lemboundk*} in Appendix \ref{appprop} below, we give some properties of $\kappa^*(d,\nu,y)$.

\bigskip


Now, we recall the following results on non characteristic points:
\begin{prop}[Case of non characteristic points]\label{threg}$ $\\
(i) {\bf (Energy level criterion for $\RR$)} If for some $x_0$ and $s_0\ge -\log T(x_0)$, we have 
\[
E(w_{x_0}(s_0))<2 E(\kappa_0),
\]
then $x_0\in \RR$.\\
(ii) {\bf (A trapping criterion for $\RR$)} There exist $\epsilon_0>0$ and $C_0>0$ such that if for some $x_0\in\R$, $s_0\ge -\log T(x)$, $\theta=\pm 1$, $d\in(-1,1)$ and $\epsilon\in(0, \epsilon_0]$, we have 
\begin{equation}\label{trap}
\left\|\vc{w_{x_0}(s_0)}{\partial_s w_{x_0}(s_0)}-\theta\vc{\kappa(d)}{0}\right\|_{\H}\le \epsilon,
\end{equation}
then $x_0\in\RR$ and $\left|\argth T'(x_0)- \argth d\right|\le C_0 \epsilon$.
Moreover, for all $s\ge s_0$:
\begin{equation*}
\left\|\vc{w_{x_0}(s)}{\partial_s w_{x_0}(s)}-\theta\vc{\kappa(T'(x_0))}{0}\right\|_{\H}\le C_0 e^{-\mu_0(s-s_0)}
\end{equation*}
for some positive $\mu_0$ and $C_0$ independent from $x_0$.\\ 
(iii) There exist $\mu_0>0$ and $C_0>0$ such that for all $x_0\in\RR$, there exist 
$\theta(x_0)=\pm 1$ and $s_0(x_0)\ge - \log T(x_0)$ such that for all $s\ge s_0$:
\begin{equation}\label{profile}
\left\|\vc{w_{x_0}(s)}{\partial_s w_{x_0}(s)}-\theta(x_0)\vc{\kappa(T'(x_0))}{0}\right\|_{\H}\le C_0 e^{-\mu_0(s-s_0)}.
\end{equation}
Moreover, $E(w_{x_0}(s)) \to E(\kappa_0)$ as $s\to \infty$.
\end{prop}
{\bf Remark}: The notation $\argth$ stands for the inverse of the hyperbolic tangent function.\\
{\it Proof}: All these results have been proved in \cite{MZjfa07}, \cite{MZcmp08} and \cite{MZajm10}. For the reader's convenience, we indicate where the proofs can be found in those papers:\\
(i) See Corollary 7 in \cite{MZajm10}.\\
(ii)
 Since we proved in Corollary 7 in \cite{MZajm10} that
\begin{equation}\label{corcriterion}
\forall x_0\in\R,\;\;\forall s\ge -\log T(x_0),\;\;E(w_{x_0}(s)) \ge E(\kappa_0)
\end{equation}
(this comes actually from the monotonicity of $E(w)$ and the convergence in (iii) of this proposition and Proposition \ref{thsing} below),
our statement follows from Theorem 3 page 48 in \cite{MZjfa07}, (i) of this proposition and Theorem 1 page 58 in \cite{MZcmp08}.\\
(iii) See Corollary 4 page 49, Theorem 3 page 48 in \cite{MZjfa07} and Theorem 1 page 58 in \cite{MZcmp08}.\Box

\medskip

Now, we recall our results for characteristic points:
\begin{prop}\label{thsing}{\bf (Case of characteristic points)} If $x_0\in\SS$, then it holds that
\begin{equation}\label{cprofile00}
\left\|\vc{w_{x_0}(s)}{\ps w_{x_0}(s)} - \vc{\d\sum_{i=1}^{k(x_0)} \theta_i\kappa(d_i(s))}0\right\|_{\H} \to 0\mbox{ and }E(w_{x_0}(s))\to k(x_0)E(\kappa_0)
\end{equation}
as $s\to \infty$, for some 
\begin{equation*}
k(x_0)\ge 2,\;\;\theta_i=\theta_1(-1)^{i+1}
\end{equation*}
 and continuous $d_i(s)=-\tanh \zeta_i(s)\in (-1,1)$ for $i=1,...,k(x_0)$. Moreover, for some $C_0>0$, for all $i=1,...,k(x_0)$ and $s$ large enough, we have
\begin{equation}\label{equid}
\left|\zeta_i(s)-\left(i-\frac{(k(x_0)+1)}2\right)\frac{(p-1)}2\log s\right|\le C_0. 
\end{equation}
\end{prop}
{\bf Remark}: The integer $k(x_0)\ge 2$ has a geometrical interpretation: it appears also in the upper bound estimate on $T(x)$ for $x$ near $x_0$ given in \aref{chapeau0}.
\subsection{The strategy of the proof}\label{strategy}
In the following, we will explain the strategy of the proof of Theorems \ref{thfini} and \ref{pdes} and Corollary \ref{corspeed}. Consider $u(x,t)$ a blow-up solution of equation \aref{equ} and $x_0\in \SS$. 
The decomposition of $w_{x_0}(y,s)$ 
is given in Proposition \ref{thsing} (up to replacing $u(x,t)$ by $-u(x,t)$, we may assume that $\theta_1=-1$).

\medskip
 
%
%
To prove that $x_0$ is an isolated characteristic point, the only tools we have are the energy criterion and the trapping result of (i) and (ii) in Proposition \ref{threg}. In order to use these tools, we have to find the behavior of $w_x$ for $x$ near $x_0$. 
A simple idea for that is to start from the decomposition \aref{cprofile00} for $w_{x_0}$ and the fact that the blow-up set is locally different from a straight line (it is in fact corner shaped; see the result of \cite{MZajm10} stated in \aref{chapeau0}), and use the transformation \aref{defw} first to recover the behavior of $u(x,t)$, then the behavior of $w_x(y,s)$ for $x$ near $x_0$. Two problems arise in this simple idea:

\medskip

- we can't have information on $w_x(y,s)$ for all $y\in(-1,1)$, since information on the whole interval $(-1,1)$ would involve information on $w_{x_0}(y,s)$ for $|y|\ge 1$, and this is unavailable (at least at time $s$) because of the finite speed of propagation;

\medskip 

- the relation between $w_{x_0}$ and $w_x$ we get from \aref{defw} depends explicitly on the value of $T(x)$ which is an unknown. The value of $T(x)$ specified by the range in \aref{chapeau0} changes the range of $s$ for which we have information. 

\medskip

To overcome these problems, we first use \aref{cprofile00} and continuity arguments to show that for $x$ close enough to $x_0$ and 
$s=\lzero$ large enough, $w_x$ is close to a sum of $k$ solitons 
\[
\sum_{i=1}^k(-1)^i\kappa^*_1(\cd_i(\lzero),\tmu_i(\lzero))
\]
where $\kappa^*_1$ is defined in \aref{defk*}. Remember that $\kappa^*_1(d,\pm e^s)$ are heteroclinic orbits connecting $\kappa(d)$ 
to $0$ or to $\infty$.
Here, we are going to use essential facts of the theory of "solitons", namely
that this decomposition is stable in time as time increases and that there are no collisions between them. 
%
Then, the idea is to use the equation \aref{defw} satisfied by $w_x$ to propagate this decomposition from $s=\lzero$ to $|\log|x-x_0||+L$ where $L$ is large, and prove (roughly speaking) the following (see Proposition \ref{propall} below for a precise statement): 
%
%
\begin{equation}\label{xdecomp}
\sup_{\lzero\le \ts \le |\log|x-x_0||+L}\left\|\vc{w_{\tx}(\ts)}{\ps w_{\tx}(\ts)}- \sum_{i=1}^{k}(-1)^{i}\kappa^*\left(\cd_i(\ts),\tmu_i(\ts)\right)\right\|_{\H}\to 0 \end{equation}
as $L_0\to \infty$, $L\to \infty$ and $x\to x_0$
for some parameters $(\cd_i(\ts),\tmu_i(\ts))$.\\
Then, it happens that at time $s=|\log|x-x_0||+L$, all the solitons $\kappa^*(\cd_i(\ts),\tmu_i(\ts))$ for $i=2,...,k$ become small (or vanish) for $L$ large and $|x-x_0|$ small (see Claim \ref{clvanish} for a precise statement), so that only the first soliton is left in \aref{xdecomp}. This crucial property is derived from the decomposition of $w_{x_0}(y,s)$ given in Proposition \ref{thsing} and the relation between $w_{x_0}$ and $w_x$ derived from \aref{defw}. 
 Since 
\[
\forall s\ge -\log T(x), \;\;E(w_x(s))\ge E(\kappa_0)
\]
(see \aref{corcriterion}),
the first soliton has to be a pure soliton of the form $-\kappa(\cd^*_1,0)$ given in \aref{defkd}, for some explicit $\cd^*_1=\cd_1^*(x)$, leading to the following estimate:
\[
\left\|\vc{w_{\tx}(s^*)}{\ps w_{\tx}(s^*)}+\vc{\kappa\left(\cd^*_1\right)}{0}\right\|_{\H}\le \epsilon_0\mbox{ where }s^*=|\log|x-x_0||+L
\]
and $\epsilon_0>0$ is defined in the trapping result stated in (ii) of Proposition \ref{threg}. Applying that trapping result, we derive two facts:\\
- the point $x$ is non characteristic (this is the conclusion of Theorem \ref{thfini});\\
- the slope $T'(x)$ satisfies $|\argth T'(x) - \argth \cd_1^*(x)|\le C \epsilon_0$, which gives the desired estimate in Theorem \ref{pdes}. By integration, we get the approximation for $T(x)$ stated in \aref{tx}. Thus, we reach the conclusions of Theorems \ref{thfini} and \ref{pdes}. As for corollary 4, it simply follows from the estimates used for Theorems \ref{thfini} and \ref{pdes}, as we will see in the short section \ref{secspeed} below.

\medskip

In conclusion, two ideas are crucial in the propagation technique which yields \aref{xdecomp}:

\medskip

- the sum of $k$ solitons $\kappa^*(\cd_i, \tmu_i)$ is stable in time if they do not collide (that is when $|\argth\frac{\tmu_{i+1}}{1+\cd_{i+1}}-\argth\frac{\tmu_i}{1+\cd_i}|$ is large);

\medskip

- the partial information we get from the simple idea based on the algebraic transformation from $w_{x_0}$ to $w_x$ and derived from \aref{defw} indicates that the $(k-1)$ solitons for $i=2,...,k$  become close to zero.

\bigskip

The paper is organized as follows. We first give in section \ref{secmod} a modulation technique in similarity variables to control the solution near a decoupled sum of generalized solitons $\pm \kappa^*(\cd_i, \tmu_i)$. 
We then prove Theorems \ref{thfini} and \ref{pdes} in Section \ref{secfini} and Corollary \ref{corspeed}.

\section{Modulation of solutions of equation \aref{equ} near a decoupled sum of generalized solitons}\label{secmod}

\subsection{The modulation result}
In this section, we consider a function close to a decoupled sum of generalized solitons $\kappa^*$ \aref{defk*} and decompose it in a unique way so to respect some orthogonality conditions. In fact, we will extend our former results proved near the sum of decoupled pure solitons $\kappa(d)$ in \cite{MZajm10}. In some sense, we now add one parameter to the solitons (the parameter $\nu$) which gives the heteroclinic orbit $\kappa^*(d,\nu)$ between the pure soliton $\kappa(d)$ and $0$ or $\infty$.\\
%
%
%
 We first introduce for $l=0$ or $1$, for any $d\in (-1,1)$ and $r\in \H$, 
\begin{equation}\label{defpdi}
\pi_l^d(r) =\phi\left(W_l(d), r\right)
\end{equation}
where
\begin{eqnarray}
-&&\phi(q,r)= \int_{-1}^1 \left(q_1r_1+q_1' r_1' (1-y^2)+q_2r_2\right)\rho dy
=\int_{-1}^1 \left(q_1\left(-\L r_1+r_1\right) +q_2 r_2\right)\rho dy,\nonumber\\
-&&W_l(d,y)= (W_{l,1}(d,y), W_{l,2}(d,y))\label{defPhi}
\end{eqnarray}
with
\begin{equation}\label{defWl2-0}
W_{1,2}(d,y)(y)= c_1(d) \frac{(1-d^2)^{\frac 1{p-1}}(1-y^2)}{(1+dy)^{\frac 2{p-1}+1}},\;\;
W_{0,2}(d,y) = c_0\frac {(1-d^2)^{\frac 1{p-1}}(y+d)}{(1+dy)^{\frac 2{p-1}+1}},
\end{equation}
 and $W_{l,1}(d,y)\in \H_0$ is uniquely determined as the solution of 
\begin{equation}\label{eqWl1-0}
-\L r + r =\left(l - \frac{p+3}{p-1}\right)W_{l,2}(d) - 2 y\py W_{l,2}(d)+ \frac 8{p-1} \frac{W_{l,2}(d)}{1-y^2}.
\end{equation}
See \aref{defcid} below for the definition of $c_1(d)$ and $c_0$.
We now claim the following:
\begin{prop}[A modulation technique]\label{lemode0}
For all $A\ge 1$, there exist $E_0(A)>0$ and $\epsilon_0(A)>0$ such that for all  $E\ge E_0$ and $\epsilon\le \epsilon_0$, 
if $\k\in[1,k]$,
$v\in \H$ and for all $i=1,...,\k$, $(\bar d_i,\bar \nu_i)\in(-1,1)\times \R$ are such that
\begin{equation}\label{condnu}
-1+\frac 1A \le \frac{\bar \nu_i}{1-|\bar d_i|}\le A,\;\;
\bar \zeta_{i+1}^*-\bar \zeta_i^*\ge E\mbox{ and }\|\bar q\|_{\H}\le \epsilon
\end{equation}
where $\bar q = v-\d\sum_{j=1}^{\k}(-1)^j \kappa^*(\bar d_j,\bar \nu_j)$ and $\bar d_i^*=\frac{\bar d_i}{1+\bar \nu_i} = -\tanh \bar \zeta_i^*$, then, there exist $(d_i,\nu_i)$ such that for all $i=1,...,\k$ and $l=0,1$,
\begin{eqnarray}
-&&\pi_l^{d_i^*}(q)=0
\mbox{ where }q=v-\sum_{j=1}^{\k}(-1)^j \kappa^*(d_j,\nu_j),
\label{ortho}\\
-&&\left| \frac{\nu_i}{1-|d_i|}- \frac{\bar \nu_i}{1-|\bar d_i|}\right|+|\zeta_i^*-\bar \zeta_i^*|\le C(A)\|\bar q\|_{\H}\le C(A)\epsilon,\label{proche}\\
-&&-1+\frac 1{2A} \le \frac{\nu_i}{1-|d_i|}\le A+1,\;\;
\zeta_{i+1}^*-\zeta_i^*\ge \frac E2\mbox{ and }\|q\|_{\H} \le C(A)\epsilon\label{condnu1}
\end{eqnarray}
%
where $d_i^*=\frac{d_i}{1+\nu_i} = -\tanh \zeta_i^*$. 
\end{prop}
{\bf Remark}: When $\bar \nu_i=0$ for all $i=1,...,\k$, we have already proved and used this result in \cite{MZajm10} (see Lemma 3.9 there), with only one orthogonality condition for each $i$, namely $\pi_0^{d_i}(q)=0$ (note that in this case, $d_i=d^*_i$). Here, the fact that $\nu_i$ is allowed to move enables us to add an additional constraint in the problem, namely $\pi_1^{d^*_i}(q)=0$.

\medskip

Our modulation lemma follows from the following technical claim: 
\begin{cl} \label{clnond'}If $(d_i, \nu_i)$ for $i=1,...,\k$ satisfy \aref{condnu1} for some $A\ge 1$ and $E>0$ large enough, then we have for some $C^*=C^*(A)>0$:\\
(i) for $i=j$,
\[
\begin{array}{rclrcl}
 \pi_0^{d_i^*}(\pnu \kappa^*(d_i,\nu_i))&=& 0,&-\frac {C^*}{1-{d_i^*}^2}\le \pi_1^{d_i^*}(\pnu \kappa^*(d_i,\nu_i))&\le& -\frac 1{C^*(1-{d_i^*}^2)},\\
|\pi_1^{d_i^*}(\pd \kappa^*(d_i,\nu_i))|&\le &\frac{C^*}{1-{d_i^*}^2},&-\frac {C^*}{1-{d_i^*}^2}\le \pi_0^{d_i^*}(\pd \kappa^*(d_i,\nu_i))&\le& -\frac 1{C^*(1-{d_i^*}^2)}
\end{array}
\]
where $d_i^*=\frac{d_i}{1+\nu_i}$.\\
(ii) for $i\neq j$ and $l=0$ or $1$,
\[
\left|\pi_l^{d_i^*}(\pnu \kappa^*(d_j,\nu_j))\right|+\left|\pi_l^{d_i^*}(\pd \kappa^*(d_j,\nu_j))\right|\le \frac{C^*J_{\k}}{1-{d^*_j}^2}
\]
where 
\begin{equation}\label{defjbar}
J_{\k} = \sum_{n=1}^{\k-1} e^{-\frac 1{p-1}(\zeta_{n+1}^*-\zeta_n^*)}.
\end{equation}
\end{cl}

\bigskip

Let us first use Claim \ref{clnond'} to prove Proposition \ref{lemode0}. Then, we will prove Claim \ref{clnond'}. Before that, let us recall from page 86 in \cite{MZjfa07} the values of the constants $c_1(d)$ and $c_0$ appearing in the definition of $W_{l,2}(d,y)$ \aref{defWl2-0}
\begin{eqnarray}
\frac 1{c_1(d)} &=&(1-d^2)^{\frac {p+1}{p-1}}(I_2(d)+\frac{2(p+1)}{p-1}I_3(d))\in[\frac 1C, C],\label{defcid}\\
\frac 1{c_0} &=&\frac {4(1-d^2)^{\frac 2{p-1}}}{p-1} \int_{-1}^1 \frac{(y+d)^2}{(1+dy)^{\frac 4{p-1}+2}} \frac \rho{1-y^2} dy =\frac 4{p-1}\iint Y^2 (1-Y^2)^{\frac 2{p-1}-1}dY\nonumber
\end{eqnarray}
(we performed the change of variables $Y=\frac{y+d}{1+yd}$ in the second line), with 
\begin{equation}\label{defIn}
I_n(d)=\int_{-1}^1 \frac{(1-y^2)\rho }{(1+dy)^{\frac 4{p-1}+n}}.
\end{equation}
We recall also from Lemma 4.4 page 85 in \cite{MZjfa07} that for all $ d\in(-1,1)$,
\begin{equation}\label{normw}
\|W_l(d)\|_{\H} +(1-d^2)\|\pd W_l(d)\|_{\H}+\|F_l(d)\|_{\H}  \le C\mbox{ and }\phi(F_m(d), W_l(d))=\delta_{m,l}
\end{equation}
where
\begin{equation}\label{deffld}
F_1(d,y)=(1-d^2)^{\frac p{p-1}}\vc{(1+dy)^{-\frac 2{p-1}-1}}{(1+dy)^{-\frac 2{p-1}-1}},\;\; F_0(d,y)=(1-d^2)^{\frac 1{p-1}}\vc{\d\frac{y+d}{(1+dy)^{\frac 2{p-1}+1}}}{0}.
\end{equation}
Let us also introduce the following change of variables
\begin{equation}\label{chng}
(d,\nu) \mapsto (\zeta, \eta)= (\zeta(d,\nu), \eta(d,\nu))= (-\argth d, \frac \nu{1-d^2})
\end{equation}
and its inverse
\begin{equation}\label{defjj}
(\zeta, \eta)\mapsto (d,\nu) =(d(\zeta, \eta),\nu(\zeta, \eta)) =(-\tanh \zeta, \eta(1-\tanh^2 \zeta)).
\end{equation}

\medskip

{\it Proof of Proposition \ref{lemode0} assuming Claim \ref{clnond'}}: Considering an integer $\k \in[1,k]$, we define
\begin{equation}\label{defPsi}
\begin{array}{llll}
\Psi:&\H\times(\R\times (-1,\infty))^{\k}&\rightarrow&\R^{2\k}\\
&(v,(\zeta_i,\eta_i)_{i=1,...,\k})&\mapsto& (\pi_1^{d_i^*}(q),\pi_0^{d_i^*}(q))_{i=1,...,\k}
\end{array}
\end{equation}
with $q=v-\d\sum_{j=1}^{\k}(-1)^j \kappa^*(d_j,\nu_j)$, $(d_j, \nu_j)=(d(\zeta_j, \eta_j),\nu(\zeta_j, \eta_j)) $ defined in \aref{defjj}
and ${d^*_j}= \frac{d_j}{1+\nu_j}$.\\
Given some positive $A$, $E$ and $\epsilon$ together with $v\in\H$ and $(\bar d_i, \bar \nu_i)\in((-1,1)\times \R)^m$ satisfying \aref{condnu}, the conclusion follows from the application of the implicit function theorem to $\Psi$ near the point
\[
(\bar v,(\bar \zeta_i,\bar \eta_i)_{i=1,...,\k})=\left(\sum_{i=1}^{\k}(-1)^i \kappa^*(\bar d_i,\bar \nu_i), (\zeta(\bar d_i,\bar \nu_i),\eta(\bar d_i,\bar \nu_i))_{i=1,...,\k}\right).
\]
Three facts have to be checked:\\
1- Note first that
\begin{equation}\label{im1}
\Psi\left(\bar v, (\bar \zeta_i,\bar \eta_i)_{i=1,...,\k}\right)=0.
\end{equation}
2- Then, we compute from the definitions \aref{defPsi}, \aref{defpdi} and \aref{defPhi} of $\Psi$, $\pi_l^d$ and $\phi$, for $l=0,1$ and $i,j\in\{1,...,\k\}$, 
\begin{eqnarray}
D_v\Psi(v, (\zeta_j, \eta_j)_{j=1,...,\k})(u) &=&(\pi_1^{{d^*_j}}(u),\pi_0^{{d^*_j}}(u))_{j=1,...,\k},\label{Dv}\\ 
\frac{\partial}{\partial \zeta_j} \pi_l^{d_i^*}(q)&=& (-1)^{j+1} \pi_l^{d_i^*}(\partial_\zeta \kappa^*(d_j, \nu_j))+\delta_{i,j}\frac{\partial {d^*_j}}{\partial \zeta_j}\phi(\pd W_l({d^*_j}), q),\nonumber\\
 \frac{\partial}{\partial \eta_j} \pi_l^{d_i^*}(q)&=& (-1)^{j+1} \pi_l^{d_i^*}(\partial_\eta \kappa^*(d_j, \nu_j))+\delta_{i,j}\frac{\partial {d^*_j}}{\partial \eta_j}\phi(\pd W_l({d^*_j}), q),\nonumber
\end{eqnarray}
where we differentiate the function $\kappa^*(d,\nu)$ with respect to $(\zeta,\eta)$ defined by the change of variables \aref{chng}.\\
Assume now that $(v,(\zeta_i,\eta_i)_{i=1,...,\k})$ satisfies condition \aref{condnu1} where $q$ is defined in \aref{ortho}. Note from the definition \aref{deflambda} of $\lambda(d,\nu)$ and \aref{condnu1} that we have
\begin{equation}\label{boundlnu}
\forall i=1,...,{\k},\;\;(3A)^{-\frac 2{p-1}} \le \lambda(d_i,\nu_i) \le \left(2A\right)^{\frac 2{p-1}}\mbox{ and }-1+\frac 1{2A} \le \nu_i \le 2A.
\end{equation}
Using the definitions \aref{defpdi} and \aref{defPhi} of $\pi_l^d$ and $\phi$, together with \aref{normw}, we see from \aref{Dv} that
\[
\|D_v\Psi(v, (\zeta_i, \eta_i)_{i=1,...,\k})\|\le C.
\]
3- Let $M$ be the Jacobian matrix of $\Psi$ with respect to 
$(\zeta_i, \eta_i)_{i=1,...,\k}$. Since we have from the change of variables \aref{defjj}
\begin{equation}\label{jac}
\Jac_{(\zeta, \eta)}(d, \nu) = (1-d^2) 
\left(
\begin{array}{cc}
-1&0\\
\d\frac{2\nu d}{1-d^2}&1
\end{array}
\right)
\end{equation}
and from \aref{boundlnu}
\[
\left|\frac{\partial d_i^*}{\partial d_i}\right|+\left|\frac{\partial d_i^*}{\partial \nu_i}\right|\le \frac 2{1+\nu_i} \le 4A,
\]
we see from (ii) of Claim \ref{clnond'} together with  \aref{normw} that \aref{Dv} yields
\begin{equation}\label{oh}
M=M_0+M_R\mbox{ where }|M_R|\le C J_{\k}+C\|q\|_{\H},
\end{equation}
and $M_0=\diag(P(d_1,\nu_1),...,(-1)^{\k+1}P(d_\k,\nu_\k))$ is a $(2\k)\times (2\k)$ matrix, with 
\[
P(d,\nu)= 
\left(
\begin{array}{cc}
\pi_1^{d^*}(\partial_\zeta \kappa^*)&\pi_1^{d^*}(\partial_\eta \kappa^*)\\
\pi_0^{d^*}(\partial_\zeta \kappa^*)&\pi_0^{d^*}(\partial_\eta \kappa^*)
\end{array}
\right)=
\left(
\begin{array}{cc}
\pi_1^{d^*}(\partial_d \kappa^*)&\pi_1^{d^*}(\partial_\nu \kappa^*)\\
\pi_0^{d^*}(\partial_d \kappa^*)&\pi_0^{d^*}(\partial_\nu \kappa^*)
\end{array}
\right)
\Jac_{(\zeta, \eta)}(d, \nu)
\]
 and $d^*=\frac d{1+\nu}$. 
Since we easily check from \aref{condnu1}, \aref{boundlnu} and \aref{jac} that
\begin{equation}\label{id9A}
\frac 1{12A^2}\le \frac{1-{d^*_i}^2}{1-d_i^2}\le 8A^2
\end{equation}
and 
\begin{equation}\label{jac-0}
\mbox{the matrix }\Jac_{(\zeta,\eta)}(d,\nu)\mbox{ is invertible with }\|\Jac_{(\zeta,\eta)}(d,\nu)^{-1}\|\le \frac{C(A)}{1-d^2}
\end{equation}
 whenever $\frac{|\nu|}{1-|d|}\le A+1$, using 
(i) of Claim \ref{clnond'}, we see that $\pi_0^{d_i^*}(\partial_\nu \kappa^*(d_i, \nu_i))=0$ and the matrix $P(d_i, \nu_i)$ is invertible with $\|P(d_i, \nu_i)^{-1}\|\le C(A)$. Therefore, the matrix $M_0$ is invertible too and $\|M_0^{-1}\|\le C(A)$. From \aref{oh}, the definition \aref{defjbar} of $J_{\k}$ and condition \aref{condnu1}, we see that for $E$ large enough and $\epsilon$ small enough, 
\[
\mbox{the matrix }M\mbox{ is also invertible with }\|M^{-1}\|\le C(A).
\]
{\it Conclusion}:
From 1-, 2- and 3-, we see that the implicit function theorem applies and given $v\in\H$, $(\bar d_i, \bar \nu_i)$ for $i=1,...,\k$ satisfying \aref{condnu}, we get the existence of other parameters $(d_i,\nu_i)$ such that for all $i=1,...,m$ and $l=0$ or $1$, 
\[
|\zeta_i-\bar \zeta_i|+|\eta_i-\bar \eta_i|\le \|\bar q\|_{\H}\le C(A)\epsilon
\]
and \aref{ortho} holds, where $\bar q$ is defined after \aref{condnu}
and the parameters are defined by the change of variables \aref{chng}.\\
%
Since 
\[
(\zeta,\eta)\mapsto \kappa^*(d,\nu)\mbox{ and }(\zeta, \eta) \mapsto \left(-\argth\left(\frac d{1+\nu}\right), \frac \nu{1-|d|}\right)
\]
are Lipschitz (see (ii) of Lemma \ref{lemboundk*} for the first and use composition for the second), we see that \aref{proche} holds, hence \aref{condnu1} holds too.
This gives the conclusion of Proposition \ref{lemode0}, assuming Claim \ref{clnond'}. It remains to prove Claim \ref{clnond'} in order to finish the proof of Proposition \ref{lemode0}.

\subsection{Proof of Claim \ref{clnond'}}
This section is dedicated to the proof of Claim \ref{clnond'}. 

\medskip

{\it Proof of Claim \ref{clnond'}}:
Since $(d_i, \nu_i)$ for $i=1,...,\k$ satisfy \aref{condnu1}, note that \aref{boundlnu} holds here. 
In particular, $\nu_i>-1$ and $d_i^*=\frac{d_i}{1+\nu_i}$ as well as $\zeta_i^*=-\argth d_i^*$ are well defined.
%
%
For simplicity, we drop down the subscript $i$ in $d_i$, $\nu_i$ and $d_i^*$. Using the definition \aref{defk*} of $\kappa^*$, we write
\begin{eqnarray*}
\pnu \kappa^*(d,\nu,y) &=& \vc{-\frac{2\kappa_0}{p-1}\frac{(1-d^2)^{\frac 1{p-1}}}{(1+\nu+dy)^{\frac 2{p-1}+1}}}{-\frac{2\kappa_0}{p-1}\frac{(1-d^2)^{\frac 1{p-1}}}{(1+\nu+dy)^{\frac 2{p-1}+1}}+\frac{2(p+1)}{(p-1)^2}\kappa_0 \nu \frac{(1-d^2)^{\frac 1{p-1}}}{(1+\nu+dy)^{\frac 2{p-1}+2}}},\\
\pd \kappa^*(d,\nu,y)&=&\vc{-\frac{2d\kappa_0}{p-1}\frac{(1-d^2)^{\frac 1{p-1}-1}}{(1+\nu+dy)^{\frac 2{p-1}}}- \frac{2y\kappa_0}{p-1}\frac{(1-d^2)^{\frac 1{p-1}}}{(1+\nu+dy)^{\frac 2{p-1}+1}}}
{\frac{4d\kappa_0\nu}{(p-1)^2}\frac{(1-d^2)^{\frac 1{p-1}-1}}{(1+\nu+dy)^{\frac 2{p-1}+1}}+ \frac{2(p+1)\kappa_0}{(p-1)^2}y\nu\frac{(1-d^2)^{\frac 1{p-1}}}{(1+\nu+dy)^{\frac 2{p-1}+2}}}.
\end{eqnarray*}
 Using the definitions \aref{deflambda} and \aref{deffld} of $\lambda$ and $F_l({d^*})$ together with \aref{boundlnu}, we see after straightforward calculations that
\begin{eqnarray}
\frac 1{L_1} \pnu \kappa^*(d,\nu,y) &=& -F_1({d^*},y) + \frac{p+1}{(p-1)}\frac \nu{(1+\nu)} \vc{0}{\frac{(1-{d^*}^2)^{\frac p{p-1}}}{(1+{d^*} y)^{\frac 2{p-1}+2}}}\label{dnu}\\
\frac 1{L_0} \pd \kappa^*(d,\nu,y) &=& -F_0({d^*,y})-\frac{d^*\nu(\nu+2)}{1-{d^*}^2} F_1({d^*,y})+\vc{0}{\DD}\label{expdd}
\end{eqnarray}
where (from the definition \aref{deflambda} of $\lambda$ and \aref{boundlnu})
\begin{eqnarray}
L_1& =& \frac{2\kappa_0(1-d^2)^{\frac 1{p-1}}(1+\nu)^{-\frac {p+1}{p-1}}}{(p-1)(1-{d^*}^2)^{\frac p{p-1}}} =\frac{2\kappa_0\lambda(1+\nu)^{-1}}{(p-1)(1-{d^*}^2)}\in I^*,
\label{defL1}\\
L_0 &=&\frac{2\kappa_0(1-d^2)^{\frac 1{p-1}-1}(1+\nu)^{-\frac {p+1}{p-1}}}{(p-1)(1-{d^*}^2)^{\frac 1{p-1}}}=\frac{2\kappa_0\lambda^{2-p}(1+\nu)^{-3}}{(p-1)(1-{d^*}^2)}\in I^*,
\label{defL0}\\
\DD&=& d\nu \frac{F_{1,2}({d^*},y)}{1-{d^*}^2}\left(\frac{\nu+2}{\nu+1}+\frac 2{p-1}\right)+\left(\frac{p+1}{p-1}\right) \frac \nu {\nu+1}y(1-d^2)\frac{(1-{d^*}^2)^{\frac 1{p-1}}}{(1+{d^*}y)^{\frac 2{p-1}+2}}\label{defD}
\end{eqnarray}
and $I^*=[\frac 1{C^*(1-{d^*}^2)}, \frac {C^*}{(1-{d^*}^2)}]$).

\medskip

(i) We first project $\pnu \kappa^*$ then $\pd \kappa^*$.

\medskip

- Projecting identity \aref{dnu} with $\pi_l^{{d^*}}$ defined in \aref{defpdi} and $l=1$ or $0$, we write from the duality relation in \aref{normw}
\begin{equation}\label{linda*}
\pi_l^{{d^*}}\left(\pnu \kappa^*(d,\nu)\right)= L_1\left(-\delta_{l,1}+
\frac{p+1}{(p-1)}\frac \nu{(1+\nu)}(1-{d^*}^2)^{\frac p{p-1}}\iint\frac{W_{l,2}({d^*},y)\rho}{(1+{d^*}y)^{\frac 2{p-1}+2}}dy\right).
\end{equation}
When $l=0$, we write from the definition \aref{defWl2-0} of $W_{0,2}({d^*})$, 
\begin{eqnarray*}
(1-{d^*}^2)^{\frac p{p-1}}\iint\frac{W_{l,2}({d^*},y)\rho}{(1+{d^*}y)^{\frac 2{p-1}+2}}dy&=&c_0(1-{d^*}^2)^{\frac{p+1}{p-1}}\iint\frac{(y+{d^*})(1-y^2)^{\frac 2{p-1}}}{(1+{d^*}y)^{\frac 4{p-1}+3}} dy\\
&=& c_0 \iint Y(1-Y^2)^{\frac 2{p-1}}dY=0
\end{eqnarray*}
where $Y=\frac{y+{d^*}}{1+y{d^*}}$, hence the first estimate of (i) follows from \aref{linda*}.\\
When $l=1$, we write from \aref{linda*} and the definitions  \aref{defWl2-0}, \aref{defcid} and \aref{defIn} of $W_{1,2}({d^*})$, $c_1({d^*})$ and $I_3(d^*)$,
\begin{eqnarray*}
&&\pi_1^{{d^*}}\left(\pnu \kappa^*(d,\nu)\right)\\
&=& -c_1({d^*})(1-{d^*}^2)^{\frac {p+1}{p-1}}L_1\left(\frac{(1-{d^*}^2)^{-\frac {p+1}{p-1}}}{c_1({d^*})}-\frac{(p+1)\nu}{(p-1)(1+\nu)}I_3({d^*})\right)\\
&=& - c_1({d^*})(1-{d^*}^2)^{\frac {p+1}{p-1}}L_1\left(I_2({d^*})+\frac{p+1}{p-1} (2-\frac \nu{1+\nu})I_3({d^*})\right)
\end{eqnarray*}
Since $1\le 2-\frac \nu{1+\nu}\le 1+2A$ from \aref{boundlnu} and
\begin{equation}\label{I2I3}
\frac 1C \le I_n({d^*})(1-{d^*}^2)^{\frac 2{p-1}+n-2}\le C
\end{equation}
 from \aref{defIn}, 
using \aref{defL1}, we see that the second estimate of (i) follows.

\medskip

- Projecting identity \aref{expdd} with $\pi_l^{{d^*}}$ defined in \aref{defpdi} where $l=0$ or $1$, we write from the duality identity \aref{normw},
\begin{equation}\label{proj0}
\pi_l^{{d^*}}(\partial_d \kappa^*(d,\nu))= L_0(-\delta_{l,0}-\frac{d^*\nu(\nu+2)}{1-{d^*}^2}\delta_{l,1}+\iint W_{l,2}({d^*}) \DD \rho dy).
\end{equation}
When $l=1$, since 
\[
\frac{1-d^2}{1+{d^*}y} \le 2\frac{(1-d^2)}{1-{d^*}^2}=2\lambda^{p-1}(1+\nu)^2 \le C^*
\]
from the definition \aref{deflambda} of $\lambda(d,\nu)$ and \aref{boundlnu}, we write from the definitions \aref{defD} and \aref{deffld} of $\DD$ and $F_1({d^*})$ together with \aref{boundlnu},
\begin{equation}\label{boundD}
|\DD|\le C^*|\nu|\frac{(1-{d^*}^2)^{\frac 1{p-1}}}{(1+{d^*}y)^{\frac 2{p-1}+1}}.
\end{equation}
Using the definition \aref{defWl2-0} of $W_{1,2}({d^*},y)$, we then write from \aref{boundD} and the definition \aref{defcid} of $I_2({d^*})$,
\begin{equation}\label{ines}
\left|\iint W_{1,2}({d^*},y) \DD\rho dy\right| \le C^*c_1({d^*})(1-{d^*}^2)^{\frac 2{p-1}}I_2({d^*})\le C^*
\end{equation}
where we use in the last bound \aref{defcid} and \aref{I2I3}.
Since $\left|\frac{d^*\nu(\nu+2)}{1-{d^*}^2}\right|\le C^*$ by \aref{condnu1} and \aref{boundlnu}, gathering \aref{ines}, \aref{proj0} and \aref{defL0} gives the third estimate in (i) of Claim \ref{clnond'}.\\
 When $l=0$, using the change of variables $Y=\frac{y+{d^*}}{1+y {d^*}}$, we write from the definition \aref{defWl2-0} and \aref{deffld} of $W_{0,2}({d^*})$ and $F_1({d^*})$, 
\begin{eqnarray}
\iint W_{0,2}({d^*})F_{1,2}({d^*}) \rho dy &=& 
c_0(1-{d^*}^2) \iint \frac{Y(1-Y^2)^{\frac 2{p-1}}}{1-Yd^*} dY,\nonumber\\
\iint W_{0,2}({d^*})y\frac{(1-{d^*}^2)^{\frac 1{p-1}}}{(1+{d^*} y)^{\frac 2{p-1}+2}}\rho dy &=& 
\frac{c_0}{1-{d^*}^2} \iint Y\frac{(Y-d^*)}{1-Yd^*}(1-Y^2)^{\frac 2{p-1}} dY.\label{pierre}
\end{eqnarray}
Using the fact that $\iint Y(1-Y^2)^{\frac 2{p-1}}dY=0$, we see that
\begin{eqnarray}
\iint \frac{Y(1-Y^2)^{\frac 2{p-1}}}{1-Yd^*} dY&=&d^*\iint \frac{Y^2(1-Y^2)^{\frac 2{p-1}}}{1-Y{d^*}} dY,\nonumber\\
\frac{1}{1-{d^*}^2} \iint Y\frac{(Y-d^*)}{1-Yd^*}(1-Y^2)^{\frac 2{p-1}} dY&=&\iint \frac{Y^2(1-Y^2)^{\frac 2{p-1}}}{1-Y{d^*}} dY.\label{secret}
\end{eqnarray}
Since we have from the change of variables $y=-Y$,
\begin{eqnarray*}
\iint \frac{Y^2(1-Y^2)^{\frac 2{p-1}}}{1-Y{d^*}} dY&=& \frac 12 \left(\iint \frac{Y^2(1-Y^2)^{\frac 2{p-1}}}{1-Y{d^*}} dY+\iint \frac{y^2(1-y^2)^{\frac 2{p-1}}}{1+y{d^*}} dy\right)\\
&=&\iint \frac{Y^2(1-Y^2)^{\frac 2{p-1}}}{1-Y^2{d^*}^2} dY,
\end{eqnarray*}
it follows from \aref{defD}, \aref{secret} and \aref{pierre} that
\[
\iint W_{0,2}(d^*) \DD \rho dy= c_0\frac \nu{1+\nu} \left\{ \frac{d^2}{1+\nu}+\frac{p+1}{p-1}\right\}\iint \frac{Y^2(1-Y^2)^{\frac 2{p-1}}}{1-Y^2{d^*}^2} dY.
\]
Using \aref{proj0} and the definition \aref{defcid} of $c_0$, we write
\begin{eqnarray}
&&\frac 1{c_0L_0}\pi_0^{{d^*}}(\partial_d \kappa^*(d,\nu))\label{sonia}\\
&=&-\frac 4{p-1}\iint Y^2(1-Y^2)^{\frac 2{p-1}-1}dY+(1-x) \left(xd^2+\frac{p+1}{p-1}\right)\iint \frac{Y^2(1-Y^2)^{\frac 2{p-1}}}{1-x^2d^2Y^2} dY\nonumber
\end{eqnarray}
where 
\[
x=\frac 1{\nu+1}\in\left(\frac 1{3A}, \frac 1{1-(1-\frac 1{2A})(1-|d|)}\right)\subset
\left(0,\frac 1{|d|}\right)
\]
from \aref{boundlnu} and \aref{condnu1}.
Since we have 
\begin{equation}\label{aziza}
\iint \frac{Y^2(1-Y^2)^{\frac 2{p-1}}}{1-x^2d^2Y^2} dY\le \frac 4{(3p+1)(1-x^2d^2)}\iint Y^2(1-Y^2)^{\frac 2{p-1}-1}dY
\end{equation}
(see below for the proof),
it follows from \aref{sonia} that
\begin{equation}\label{laurence}
-\frac 4{p-1}I\le\frac 1{c_0L_0}\pi_0^{{d^*}}(\partial_d \kappa^*(d,\nu))\le 
-\frac 4{p-1}Ig_d(x)\mbox{ where }I=\iint Y^2(1-Y^2)^{\frac 2{p-1}-1}dY 
\end{equation}
and
\[
g_d(x) = 1-\frac{(1-x)(xd^2(p-1)+p+1)}{(3p+1)(1-x^2d^2)}\mbox{ with }x\in\left(0, \frac 1{|d|}\right).
\]
Differentiating $g_d$, we easily write
\begin{eqnarray*}
&&(3p+1)(1-x^2d^2)^2g_d'(x) =(p+1-d^2(p-1))(d^2x^2+1)-4d^2x\\
&=& (p+1-d^2(p-1))(d^2x^2-2d^2x+1)+2d^2 x (p-1)(1-d^2)>0,
\end{eqnarray*}
hence 
\[
g_d(x)\ge g_d(0) = 1-\frac{p+1}{(3p+1)} = \frac {2p}{(3p+1)}>0,
\]
and the forth estimate  in (i) follows from \aref{laurence} and the bound \aref{defL0} on $L_0$. It remains then to prove \aref{aziza} to finish the proof of (i) Claim \ref{clnond'}.

\medskip

{\it Proof of \aref{aziza}}: 

Putting $t=Y^2$, we write for $\alpha>-1$
\begin{equation}\label{k1}
\iint Y^2(1-Y^2)^\alpha dY =\int_0^1 t^{1/2}(1-t)^\alpha dt= B\left(\frac 32, \alpha+1\right),
\end{equation}
where $B(x,y)$ is the standard Beta function. Since we know that for $\alpha = \frac 2{p-1}$, 
\[
\frac{B\left(\frac 32, \alpha+1\right)}{B\left(\frac 32, \alpha\right) } = \frac \alpha{\alpha +\frac 32} = \frac 4{3p+1}\mbox{ and }1-x^2d^2Y^2 \ge 1-x^2 d^2
\]
for all $Y\in(-1,1)$, estimate \aref{aziza} follows from \aref{k1}. This concludes the proof of (i) in Claim \ref{clnond'}.

\medskip 

(ii) Take $i\neq j$ and $l=0$ or $1$. 
Using the expressions \aref{defWl2-0}, \aref{eqWl1-0}, \aref{defcid} and \aref{defkd} of $W_{l,2}$, $W_{l,1}$, $c_l(d)$  and $\kappa(d,y)$, we write
\begin{equation}\label{fani}
|W_{l,2}(d_i^*,y)|+|\L W_{l,1}(d_i^*,y) -W_{l,1}(d_i^*,y)|\le C\frac{\kappa(d_i^*,y)}{1-y^2}.
\end{equation}
Since $|F_{1,l}({d_j^*},y)|\le C\kappa({d_j^*},y)$ by definitions \aref{deffld} and \aref{defkd}, using \aref{dnu}-\aref{defD} together with \aref{boundD}, \aref{boundlnu} and \aref{condnu1}, we write
\begin{equation}\label{bordev}
|\pnu \kappa^*(d_j,\nu_j,y)| +|\pd \kappa^*(d_j,\nu_j,y)|\le \frac C{1-{d_j^*}^2}(F_{1,1}({d_j^*},y)+F_{0,1}({d_j^*},y))\le C \frac{\kappa({d_j^*},y)}{1-{d_j^*}^2}.
\end{equation}

Using the expressions \aref{defpdi} and \aref{defPhi} of $\pi_l^{d_i^*}$ and $\phi$ together with
\aref{fani} and \aref{bordev}, we write
\begin{eqnarray*}
&&\left|\pi_l^{d_i^*}(\pnu \kappa^*(d_j,\nu_j))\right|+\left|\pi_l^{d_i^*}(\pd \kappa^*(d_j,\nu_j))\right|\le \frac C{1-{d^*_j}^2}\iint \kappa(d_i^*,y)\kappa({d^*_j},y) \frac \rho{1-y^2} dy\\
&\le& \frac C{1-{d^*_j}^2} |\zeta^*_j-\zeta_i^*|e^{-\frac 2{p-1}|\zeta^*_j-\zeta_i^*|} \le \frac C{1-{d^*_j}^2} J_{\k}
\end{eqnarray*}
with $J_{\k}$ defined in \aref{defjbar} and where we used Lemma \ref{lemtech} for the last inequality. 
This concludes the proof of Claim \ref{clnond'} and Proposition \ref{lemode0} too.\Box
%
%
%

\section{Characteristic points are isolated (Proofs of Theorems \ref{thfini} and \ref{pdes})}\label{secfini}
This section is dedicated to the proof of Theorems \ref{thfini} and \ref{pdes} and Corollary \ref{corspeed}. Consider $u(x,t)$ a blow-up solution of equation \aref{equ}. Note that we already know from Proposition 5 in \cite{MZajm10} (stated in page \pageref{old} here) that the interior of $\SS$ is empty. We will prove in this section that all the points of $\SS$ are isolated. 

\bigskip

Consider $x_0\in \SS$. From translation invariance of equation \aref{equ}, we can assume that $x_0=T(x_0)=0$, hence,
\begin{equation}\label{ff}
0\in \SS\mbox{ and }T(0)=0.
\end{equation}
In the following, we use the notation $W(Y,S)$ instead of $w_0(y,s)$ defined in the selfsimilar transformation in \aref{defw}. 
Up to replacing $u(x,t)$ by $-u(x,t)$, we see from Proposition \ref{thsing} that for some integer $k=k(0)\ge 2$, for some continuous functions $D_i(S)$, $C_0>0$ and $S_0\in \R$, we have 
\begin{eqnarray}
&&\left\|\vc{W(S)}{\ps W(S)} - \vc{\d\sum_{i=1}^{k} (-1)^{i}\kappa(D_i(S))}0\right\|_{\H} \to 0\mbox{ as }S\to \infty,\label{cprofile**}\\
\forall S\ge S_0,&&\left|\argth D_i(S)-\frac{\gamma_i}2\log S\right|\le C_0\mbox{ where }\gamma_i=(p-1)\left(\frac{k+1}2-i\right).\label{equid**}
\end{eqnarray}
Introducing for $\tx \neq 0$, $\tB=\tB(\tx)$ by
\begin{equation}\label{deftB**}
-\frac{T(\tx)}{|\tx|}=1-\tB(x)\mbox{ and }l=l(x)=|\log|\tx||, 
\end{equation}
we see from \aref{chapeau0} that
\begin{equation}\label{chap**}
0< \tB\le \frac {C_0}{l^{\gamma_1}}.
\end{equation}

\medskip

 We proceed in three steps:

\medskip

- In Section \ref{sub1*}, we consider $\tx<0$ close to $0\in\SS$ and $\k\le k$, and use a modulation technique to study $w_{\tx}$ near the sum $\sum_{i=1}^\k(-1)^i\kappa^*(d_i,\nu_i)$..

\medskip

- In Section \ref{sub2*}, we state and prove a proposition which directly implies  Theorems \ref{thfini} and \ref{pdes}. 

\medskip

- In Section \ref{secspeed}, we prove Corollary \ref{corspeed}.
\subsection{Dynamical behavior of  $w_{\tx}$  for small $\tx$}\label{sub1*}
Consider some positive $L_\hk$,....,$L_{k+1}$ where $\hk =E\left(\frac {k+1}2\right)$.
Then, for each 
$x<0$ with $|x|$ small, we introduce $s_m(x)$ as follows:
\begin{eqnarray}
s_{k+1}=s_{k+1}(x) &= &L_{k+1},\label{defsk+1}\\
s_{\k}=s_{\k}(x) &=& l+\gamma_m\log l +L_m\mbox{ if }\hat k+1\le \k \le k,\label{defsm}
\end{eqnarray}
where 
$S_0$ and $\gamma_\k$ are defined in \aref{equid**}.\\
Then, from \aref{equid**} and the continuity of $D_1$, we define $\tsa=\tsa(x)$ 
maximal such that
\begin{equation}\label{deftsa}
\tsa\in[L_{k+1}, l+L_\hk]\mbox{ and }\forall \ts\in [L_{k+1}, \tsa],\;\;\frac{\tmu_1(x,\ts)}{1-\cd_1(x,\ts)}\ge -1+ e^{-(p-1)L_\hk},
\end{equation}
where 
\begin{equation}\label{defs2}
\begin{array}{rcl}
\tmu_i(s)&=&\tmu_i(x,\ts)=[\tB-(1-\cd_i(x,\ts))]\tx e^{\ts},\;\;\cd_i(s)=\cd_i(x,\ts)=D_i(\cs(x,s))\\
\mbox{and}&&
T(0)-e^{-\cs(x,s)}=T(\tx) - e^{-\ts}.
\end{array}
\end{equation}
Note from \aref{ff} and \aref{deftB**} that we have
\begin{equation}\label{defs2*}
-e^{-\cs(x,s)}=\tx(1-\tB) - e^{-\ts}\mbox{ and }\cs(s)=\cs(x,s)=-\log[|\tx|(1-\tB)+e^{-\ts}].
\end{equation}
Following this, we also introduce for all $\k\in[\hat k,k+1]$,
\begin{equation}\label{defSm}
\cs_\k=\cs_\k(x) = \cs(x,s_\k(x))=-\log[|\tx|(1-\tB)+e^{-s_\k(x)}].
\end{equation}

\bigskip

From \aref{deftsa}, one of the following cases occurs: 
\begin{eqnarray}
&\mbox{\bf (Case 1)}& \tsa = l+L_\hk,\label{ts2max}\\
&\mbox{\bf (Case 2)}& \tsa <l+L_\hk\mbox{ and }\label{defs2'}\frac{\tmu_1(x,\ts)}{1-\cd_1(x,\ts)}=-1+ e^{-(p-1)L_\hk}.
\end{eqnarray}
Note that a step of the proof is devoted to the proof that Case 2 does not occur and that $\tsa=l+L_\hk$, provided that 
the constants $L_\k$ for $\hk\le \k \le k+1$
are fixed large enough and $|x|$ small enough. Note also that for $|x|$ small enough, we have
\begin{equation}\label{ordre}
s_{k+1} \le s_k \le ...\le s_{\hat k+1}\le \tsa \le l+L_\hk.
\end{equation}
Indeed, all these inequalities are obvious by definition except the inequality 
$s_{\hat k+1}\le \tsa$ 
which follows from the fact that
$\frac{|\tmu_1(s_{\hat k+1})|}{1-|\cd_1(s_{\hat k+1})|}\le Ce^{L_{\hk+1}} l^{\gamma_{\hat k+1}}\to 0$ as $x\to 0$ (to prove this, see the definition \aref{defs2} of $\tmu_1$ and $\cd_1$, \aref{equid**}, the bound \aref{chap**} on $\tB$ and the definition \aref{defsm} of $s_{\hat k+1}$).

\bigskip

We now give the following decomposition result for $(w_{\tx}(\ts), \partial_{\ts} w_{\tx}(\ts))$:
\begin{prop}\label{propall}{\bf (Decomposition of $(w_{\tx}(\ts), \partial_{\ts} w_{\tx}(\ts))$ for $\ts\in[s_{k+1},\tsa]$)}
For all $\k \in [\hat k,k]$, 
if $L_\k>0$ is large enough, then it holds that
\begin{equation}\label{deftN}
\lim_{\lzero\to \infty}
\left(\lim_{x\to 0^-}
\sup_{\ts_{\k+1} \le \ts \le \ts_\k}\left\|\vc{w_{\tx}(\ts)}{\ps w_{\tx}(\ts)}- \sum_{i=1}^{\k}(-1)^{i}\kappa^*\left(\cd_i(\ts),\tmu_i(\ts)\right)\right\|_{\H}\right)= 0.
\end{equation}
\end{prop}
{\bf Remark}: What we actually prove is that if $L_m$ is large enough, $0<\epsilon\le \epsilon_0(L_m)$, $\lzero\ge L_{\k+1,0}(L_\k,\epsilon)$ and $-a_0(L_\k,\epsilon, L_{\k+1})\le x<0$ for some positive $\epsilon_0$, $L_{\k+1,0}$ and $a_0$, then
\[
\sup_{\ts_{\k+1} \le \ts \le \ts_\k}\left\|\vc{w_{\tx}(\ts)}{\ps w_{\tx}(\ts)}- \sum_{i=1}^{\k}(-1)^{i}\kappa^*\left(\cd_i(\ts),\tmu_i(\ts)\right)\right\|_{\H}\le \epsilon.
\]
{\bf Remark}: For all $i\in[1,k]$ and $s\in[s_{k+1},s_{\hat k}]$, we have
\begin{equation}\label{domi}
1+\tmu_i(s) -|\cd_i(s)|>0
\end{equation}
which means by definition \aref{defk*} of $\kappa^*$ that $\kappa^*\left(\cd_i(\ts),\tmu_i(\ts)\right)$ is well defined for all $y\in(-1,1)$. Indeed, if $\cd_i(s)\ge 0$, then we write from the definition \aref{defs2} of $\tmu_i(s)$, the expansion \aref{equid**} of $D_i$ and the definition \aref{deftsa} of $\tsa$ that $1+\tmu_i(\ts) -|\cd_i(\ts)|\ge 1+\tmu_1(\ts) -\cd_1(\ts)>0$ and \aref{domi} follows. If $\cd_i(s)<0$, then we see from \aref{defs2} and the bound \aref{chap**} on $\tB$ that $\tmu_i(s)\ge 0$ and \aref{domi} follows. 

\bigskip

{\it Proof of Proposition \ref{propall}}: In order to prove Proposition \ref{propall}, we proceed in two steps:\\
- In Step 1, we use the selfsimilar transformation \aref{defw} to change the estimate \aref{cprofile**} on $W(\cy,\cs)$ into an estimate on $w_\tx(y,s)$ where $\tx <0$ is close to $0$ and $s\in[s_{k+1},\tsa]$, valid only for $\ty\in(\ty_1(\ts),1)$ where 
\begin{equation}\label{defy1}
\ty_1(\ts)=-1-2\tB \tx e^{\ts}>-1.
\end{equation}
- In Step 2, we use the modulation technique of Section \ref{secmod} to extend this estimate to the whole interval $(-1,1)$ and prove Proposition \ref{propall}.

\bigskip

{\bf Step 1: Decomposition of $(w_{\tx}(\ts), \partial_{\ts} w_{\tx}(\ts))$ for $s\in[s_{k+1}, \tsa]$ and $y\in(\ty_1(s),1)$}

Here, we transform the estimate \aref{cprofile**} through the change of variables \aref{defw} to get the following: 
\begin{lem}\label{lemleft}{\bf (Decomposition of $(w_{\tx}(\ts), \partial_{\ts} w_{\tx}(\ts))$ on $(\ty_1(\ts),1)$)} For all $\k\in[\hk,k]$ and $L_\k>0$, it holds that
\[
\lim_{\lzero\to \infty}
\left(\lim_{x\to 0^-}\sup_{s_{\k+1}\le \ts \le s_\k}\left\|\vc{w_{\tx}(\ts)}{\ps w_{\tx}(\ts)}- \sum_{i=1}^k(-1)^{i}\kappa^*\left(\cd_i(\ts), \tmu_i(\ts)\right)\right\|_{\H(\ty>\ty_1(\ts))}\right)= 0,
\]
where $\ty_1(\ts)$ is defined in \aref{defy1}.
\end{lem}
{\bf Remark}: Note that all the solitons $\kappa^*(\cd_i(s), \tmu_i(s))$ appearing in this statement are well defined thanks to the remark following Proposition \ref{propall}. Note also that we count all the solitons in the estimate of Lemma \ref{lemleft}, unlike Proposition \ref{propall} where we count only the solitons from $i=1$ to $i=\k$. As a matter of fact, thanks to the vanishing property of Lemma \ref{clvanish}, we will use in the proof of Claim \ref{clproche} in Appendix \ref{apptrans} the fact that Lemma \ref{lemleft} holds in fact with a sum running from $i=1$ to $i=\k$ too.\\
{\it Proof}:
%
%
%
 Introducing the transformation
\begin{equation}\label{defw**}
\tT(v)(\ty, \ts)=(1-(1-\tB) \tx e^{\ts})^{-\frac 2{p-1}} v(\cy,\cs),\;\;
\cy =\frac{\ty +\tx e^{\ts}}{1-(1-\tB)\tx e^{\ts}}\;\;
\cs= \ts -\log(1-(1-\tB)\tx e^{\ts}),
\end{equation}
we see from the selfsimilar transformation \aref{defw} and the definition \aref{deftB**} of $\tB$ that
\begin{equation}\label{rach1}
\tT(W)=w_{\tx}.
\end{equation}
Using the definitions \aref{defkd} and \aref{defk*} of $\kappa(d,y)$ and $\kappa_1^*(d,\nu,\ty)$, we see that for all $d\in (-1,1)$
and $\ts \in \R$, if $1+[\tB-(1-d)]\tx e^{\ts}>|d|$, then
\begin{equation}\label{rach2}
\forall \ty \in (-1,1),\;\;
\tT(\kappa(d,\cdot))(\ty, \ts)=\kappa^*_1\left(d,[\tB-(1-d)]\tx e^{\ts},\ty\right).
\end{equation}
Therefore, defining the error terms $\chr$ and $\tr$ by
\begin{eqnarray}
&&\chr(\cy,\cs)=\vc{\chr_1(\cy,\cs)}{\chr_2(\cy,\cs)}=\vc{W(\cy,\cs)}{\ps W(\cy,\cs)}-\sum_{i=1}^k (-1)^{i}\vc{\kappa(D_i(\cs),\cy)}{0},\label{defg**}\\
&&\tr(\ty,\ts)=\vc{\tr_1(\ty,\ts)}{\tr_2(\ty, \ts)}\nonumber\\
&=& \vc{w_{\tx}(\ty,\ts)}{\ps w_{\tx}(\ty,\ts)}- \sum_{i=1}^k (-1)^{i}\kappa^*\left(\cd_i(s),\tmu_i(s),\ty\right),\label{defh**}
\end{eqnarray}
we get the following estimate from \aref{defw**}:
\begin{cl}\label{clalgebre0}{\bf (Transformation of the error terms on $(y_1(s),1)$)} For all $\k\in[\hk,k]$, $L_\k>0$ 
and $|x|$ small enough with $x<0$, we have for all $s\in [s_{\k+1}, s_\k]$,
\[
\left\|\vc{\tr_1}{\tr_2}(\ts)\right\|_{\H(\ty > \ty_1(\ts))} \le C(L_\k) \left\|\vc{\chr_1}{\chr_2}(\cs)\right\|_{\H}
\]
where $S=S(s)$ is defined in \aref{defs2*}.
\end{cl}
{\it Proof}: see Appendix \ref{apptrans}.\Box

\medskip

Since $\|(\chr_1, \chr_2)(\cs)\|_{\H} \to 0$ as $\cs\to \infty$ from \aref{cprofile**} and \aref{defg**} and 
\[
\forall s\in [s_{\k+1}, s_\k],\;\;S(s) \ge S_{\k+1}\ge \lzero-1
\]
for $|x|$ small enough from \aref{defSm} and \aref{defsk+1}, this concludes the proof of Lemma \ref{lemleft}. \Box

\bigskip

{\bf Step 2: Decomposition of $(w_x(\ts), \partial_s w_x(s))$ on $(-1,1)$ for $s\in[s_{k+1}, \tsa]$}

In this step, we consider the equation \aref{eqw} satisfied by $w$ and we linearize it around $\sum_{i=1}^\k(-1)^i\kappa^*(\cd_i(s), \tmu_i(s))$ (which is an approximate solution of \aref{eqw} in the considered regime) in order to extend the estimate of Lemma \ref{lemleft} to the whole interval $(-1,1)$ and prove Proposition \ref{propall}. 

\medskip

We first consider in the following claim the soliton $\kappa^*(\cd_i(s), \tmu_i(s))$ where $i\in[2,k]$ appearing in the decomposition we aim at proving in Proposition \ref{propall}, and show that it ``vanishes'' at time $s=s_i$ defined in \aref{defsm} when $i\ge \hat k+1$, and at time $s=l+L_\hk$ when $i\in[2,\hat k]$. By ``vanishing'', we mean that the norm of the $i$-th soliton becomes smaller than $Ce^{-\frac L{p-1}}$ (where $L=L_i$ if $i\ge \hk+1$ and $L=L_\hk$ if $i\le \hk$), which may be made as small as we wish by taking $L$ large enough. More precisely, we have the following:
\begin{cl}\label{clvanish} {\bf (Vanishing of $\kappa^*(\cd_i(s), \nu_i(s))$ for large $s$)} The following holds for $x<0$ and $|x|$ small enough:\\
- if $i\in [\hat k+1,k]$ and $s\in[s_i,l+L_\hk]$, then $\|\kappa^*(\cd_i(s), \tmu_i(s))\|_{\H} \le C\lambda(\cd_i(s), \tmu_i(s))\le C e^{-\frac {L_i}{p-1}}$;\\
- if $i\in[2,\hat k]$, then $\|\kappa^*(\cd_i(l+L_\hk), \tmu_i(l+L_\hk))\|_{\H} \le C\lambda(\cd_i(l+L_\hk), \tmu_i(l+L_\hk))\le C e^{-\frac {L_\hk}{p-1}}$,\\
where 
\begin{equation}\label{deflambda}
\lambda= \lambda(d,\nu)=\frac{(1-d^2)^{\frac 1{p-1}}}{[(1+\nu)^2-d^2]^{\frac 1{p-1}}}.
\end{equation}
\end{cl}
{\it Proof}: See Appendix \ref{apptrans}.\Box\\ 
{\bf Remark}: As we announced in the strategy of the proof given in section \ref{strategy}, when $s=s_{k+1}$ defined in \aref{defsk+1}, $w_x(s)$ is close to a sum of $k$ decoupled solitons. As time increases, we start loosing the solitons beginning with the $k$-th at time $s_k$, up to the $(\hat k+1)$-th at time $s_{\hat k+1}$, and finishing with all the solitons for $i=2$ to $\hat k$, which vanish at the same time $s=l+L_\hk$, leaving only the first soliton in the decomposition of $w_x(s)$. As a matter of fact, the values of $s_\k$ for $\k\in[\hat k,k]$ we took in \aref{defsm} and \aref{deftsa} were chosen on purpose to guarantee this gradual vanishing. One may wander why we didn't take $s_{\hat k}=l+L_\hk$ in \aref{deftsa}. The reason is that the first soliton may blow-up at some time $s\in[s_{\hat k+1},l+L_\hk]$, as one can see from direct computations based on the expressions of $\cd_1$, $\tmu_1$ \aref{defs2} and $D_1$ \aref{equid**}. As a matter of fact, a whole step of our argument is devoted to the proof that this is not the case (see (ii) of Claim \ref{lemup} below where we will prove that indeed $s_{\hat k}=l+L_\hk$).

\medskip

We also claim the following:
\begin{cl}\label{clprep} 
There exists $C_0>0$ such that for all $\k\in[\hat k, k]$, $L_\k>0$, if $\lzero$ is large enough and $x<0$ close enough to $0$, then for all $s\in[s_{\k+1}, s_\k]$, we have
\begin{eqnarray}
\forall i=1,...,\k,&&-1+e^{-(p-1)L_\k}\le \frac{\tmu_i(s)}{1-|\cd_i(s)|}\le C_0 e^L_\k,\label{first}\\
\forall i=1,...,\k-1,&&\czeta_{i+1}^*(s) - \czeta_i^*(s)\ge \frac{(p-1)}7 \log \lzero,\label{second}
\end{eqnarray}
where 
\begin{equation}\label{defcz}
\czeta_i^*(s) =-\argth\left(\frac{\cd_i(s)}{1+\tmu_i(s)}\right).
\end{equation}
\end{cl}
{\it Proof}: See Appendix \ref{apptrans}. \Box

\bigskip

 {\it Proof of Proposition \ref{propall}}: We proceed by induction for $\k$ decreasing from $k$ to $\hat k$. Let us assume that
\begin{equation}\label{alternative}
\begin{array}{rcl}
\mbox{\bf either}&& \k=k,\\
\mbox{\bf or}&& \hat k\le \k\le k-1\mbox{ and Proposition \ref{propall} holds at the level }\k+1.
\end{array}
\end{equation}
We then consider some $L_\k>0$.
We first need to check that Proposition \ref{propall} holds for $s=s_{\k+1}$.
\begin{lem}\label{lemstep1}{\bf (Estimate of  $(w_{\tx}(s_{\k+1}), \partial_{\ts} w_{\tx}(s_{\k+1}))$ on $(-1,1)$)} Under \aref{alternative}, it holds that
\[
\lim_{\lzero\to \infty}
\left(
\lim_{x\to 0^-}
\left\|\vc{w_{\tx}(s_{\k+1})}{\ps w_{\tx}(s_{\k+1})}- \sum_{i=1}^{\k}(-1)^{i}\kappa^*\left(\cd_i(s_{\k+1}),\tmu_i(s_{\k+1})\right)\right\|_{\H}
\right)
=0.
\]
\end{lem}
{\it Proof}:\\
{\bf Case $\k=k$}: We will in fact prove that for all $\epsilon>0$, there is $L^*>0$ such that for all $L_{k+1}\ge L^*$, there exists $\delta(L_{k+1})>0$ such that for all $x\in[-\delta,0)$, 
\[
\left\|\vc{w_{\tx}(s_{k+1})}{\ps w_{\tx}(s_{k+1})}- \sum_{i=1}^{k}(-1)^{i}\kappa^*\left(\cd_i(s_{k+1}),\tmu_i(s_{k+1})\right)\right\|_{\H}\le \epsilon.
\]
Consider $\epsilon>0$ and from \aref{cprofile**} take $L_{k+1}$ large enough such that for all 
$\cs\ge L_{k+1}-1$, we have
\begin{equation}\label{nesrine}
\left\|\vc{W(\cs)}{\ps W(\cs)} - \vc{\d\sum_{i=1}^{k} (-1)^{i}\kappa(D_i(\cs))}0\right\|_{\H} \le \epsilon.
\end{equation}
Note that for all $i=1,...,k$,
\begin{equation}\label{contk*}
\left\|\kappa^*\left(\cd_i(s_{k+1}), \tmu_i(s_{k+1})\right)-\vc{\kappa(D_i(s_{k+1}))}{0}\right\|_{\H} \to 0\mbox{ as }x\to 0^-
\end{equation}
(where we recall from \aref{defsk+1} that $s_{k+1}=L_{k+1}$),
from the definition \aref{defs2} of $\cd_i$ and $\tmu_i$ together with the continuity estimate in (ii) of Lemma \ref{lemboundk*}.\\
Then, using the fact that $\tB>0$ (see \aref{chap**}) and the solution of the Cauchy problem for equation \aref{equ}, we see  that for some $K=K(L_{k+1})>0$ and $\delta=\delta(L_{k+1})>0$, we have for all $\ct\in [-e^{-L_{k+1}+1}, - e^{-L_{k+1}-1}]$, 
\[
\|(u(\ct),u_t(\ct))\|_{H^1\times L^2(|x|<(1+\delta)(-\ct))}\le K(L_{k+1}).
\]
Using the selfsimilar transformation, we see that for some $K'(L_{k+1})>0$ and for all $\cs\in [s_{k+1}-1, s_{k+1}+1]$, we have
\[
\|(W(\cs), \ps W(\cs))\|_{H^1\times L^2(|\cy|<1+\delta)}\le K'(L_{k+1})
\]
and $W$ uniformly continuous for $\cs\in [s_{k+1}-1, s_{k+1}+1]$.
Using the transformation \aref{defw**} (in particular \aref{rach1} and the relations in (i) of Claim \ref{clalgebre}), we see that when $\k=k$, Lemma \ref{lemstep1} follows from Lebesgue's Theorem, \aref{contk*} and \aref{nesrine}.

\medskip

\noindent {\bf Case $\hat k\le \k \le k-1$}: Consider $\epsilon>0$. Since $k\ge \k+1\ge \hk+1$, Claim \ref{clvanish} applies and we can fix $\lzero>0$ large enough so that 
\begin{equation}\label{train}
\left\|\kappa^*\left(\cd_{\k+1}(s_{\k+1}),\tmu_{\k+1}(s_{\k+1})\right)\right\|_{\H}\le \frac \epsilon 2
\end{equation}
and Proposition \ref{propall} applies at the level $\k+1$ (induction hypothesis). Hence, it follows that fixing $L_{m+2}$ large enough and taking $|x|$ small enough, we have
\[
\left\|\vc{w_{\tx}(s_{\k+1})}{\ps w_{\tx}(s_{\k+1})}- \sum_{i=1}^{\k+1}(-1)^{i}\kappa^*\left(\cd_i(s_{\k+1}),\tmu_i(s_{\k+1})\right)\right\|_{\H}\le \frac \epsilon 2.
\]
Thus, using \aref{train}, we see that
\[
\left\|\vc{w_{\tx}(s_{\k+1})}{\ps w_{\tx}(s_{\k+1})}- \sum_{i=1}^\k(-1)^{i}\kappa^*\left(\cd_i(s_{\k+1}),\tmu_i(s_{\k+1})\right)\right\|_{\H}\le\epsilon,
\]
which concludes the proof of Lemma \ref{lemstep1}.\Box

\bigskip

From Lemma \ref{lemstep1} and Claim \ref{clprep}, we see that for $\lzero$ large enough, $|x|$ small enough and $s=s_{\k+1}$, $w_x(s_{\k+1})$ is close to a sum of $\k$ decoupled solitons. We will in fact show that this decomposition propagates for all $s\in[s_{\k+1}, s_\k]$. To do so, we need a Lyapunov functional to control the error between our solution and such a sum. The modulation technique of Proposition \ref{lemode0} is crucial in obtaining such a Lyapunov function.

\medskip

Using Lemma \ref{lemstep1} and Claim \ref{clprep}, we see that for $\lzero$ large enough and $|x|$ small enough, the modulation technique of Proposition \ref{lemode0} applies at $s=s_{\k+1}$ (in fact, it applies on a small interval to the right of $s_{\k+1}$,  thanks to the continuity in $\H$ of $(w_x(s), \ps w_x(s))$ as a function of $s$). More precisely, there are $(d_i(s), \nu_i(s))$ such that for all $\epsilon>0$ small enough, there exists $\bar s=\bar s(\epsilon) $ maximal in $(s_{\k+1}, s_\k]$ such that for all $s\in[s_{\k+1}, \bar s]$,
\begin{eqnarray}
\mbox{for }l=0,\;1, 
&&\pi_l^{d_i^*(s)}(q(s))=0,\label{mode0}\\
&&\|q(s_{\k+1})\|_{\H} \le \beta\epsilon,\label{small0}\\
&&\|q(s)\|_{\H} \le \epsilon,\label{small}\\
\forall i=1,...,\k-1,&&\zeta_{i+1}^*(s) - \zeta_i^*(s)\ge \frac{(p-1)}9 \log \lzero,\label{gap}\\
\mbox{and}&&
-1+\frac{e^{-(p-1)L_\k}}2\le \frac{\nu_i(s)}{1-|d_i(s)|}\le C_0e^{L_\k}+1,\label{apriori} 
\end{eqnarray}
where $\pi_l^{{d_i^*}}$ is defined in \aref{defpdi},
\begin{equation}\label{defq}
q=\vc{q_1}{q_2}=\vc{w_{\tx}}{\partial_s w_{\tx}}-\sum_{i=1}^\k (-1)^i\kappa^*(d_i,\nu_i),
\end{equation}
$d_i^*=\frac{d_i}{1+\nu_i}=-\tanh \zeta_i^*$, $C_0$ is given in Claim \ref{clprep} and $\beta\in(0,1)$ will be fixed later small enough. One should keep in mind that $q$ and the above mentioned variables depend also on $x$ and $m$. For simplicity in the notation, we omit that dependence.

\bigskip

Since $\bar s$ is maximal in $(s_{\k+1}, s_\k]$, we have
\begin{eqnarray}
\mbox{\bf either}&&\bar s=s_\k,\label{dicho}\\
\mbox{\bf or }&&\bar s <s_\k\mbox{ and there is an equality case in \aref{small}-\aref{apriori}}.\nonumber
\end{eqnarray}
Our goal is to prove that $\bar s=s_\k$.
Using Lemmas \ref{lemboundk*} and \ref{lemleft}, we have the following:
\begin{cl}\label{clproche}
If $L_\k$ is fixed large enough, $\epsilon\le \epsilon_1(L_\k)$, $L_{\k+1}\ge L_{\k+1,1}(L_\k,\epsilon)$ and $|x|\le a_1(L_\k,\epsilon, L_{\k+1})$ for some $\epsilon_1>0$, $L_{\k+1,1}>0$ and $a_1$, 
then we have for all $s\in [s_{\k+1}, \bar s]$:\\
(i) For all $i=1,...,\k$, 
\[
\left|\frac{\nu_i(s)}{1-|d_i(s)|}- \frac{\tmu_i(s)}{1-|\cd_i(s)|}\right|+|\argth d_i(s) - \argth\cd_i(s)|+|\zeta_i^*(s) - \czeta_i^*(s)|\le C(L_\k) \epsilon.
\]
(ii) $\|\bar q(s)\|_{\H} \le C(L_\k) \epsilon$ where $\bar q=\vc{\bar q_1}{\bar q_2}=\vc{w_{x}}{\partial_s w_{x}}-\d\sum_{i=1}^\k (-1)^i\kappa^*(\cd_i,\tmu_i)$.
\end{cl}
{\it Proof}: See Appendix \ref{apptrans}.\Box

\bigskip

From Claim \ref{clproche}, we see that it is enough to prove that
\[
\bar s = s_{\k}
\]
in order to conclude the proof of Proposition \ref{propall}, 
provided that $L_\k$ is fixed large enough, $\epsilon$ is small enough in terms of $L_\k$, $\lzero$ is large enough in terms of $L_\k$ and $\epsilon$ and $|x|$ is small enough in terms of $L_\k$, $\epsilon$ and $\lzero$.
From \aref{dicho}, we need to prove that all the inequalities in \aref{small}-\aref{apriori} are strict. This is the aim of the remaining part of the proof.

\bigskip

We claim the following:
\begin{cl} \label{lemlyap} {\bf (Remarkable functionals for $q$)}
If $L_\k$ is large enough, then
there exist $C^1$ functions $h_1$ and $h_2$ defined for all $s\in[s_{\k+1}, \bar s]$ and constants $C_*(L_\k)>0$, $\delta^*(L_\k)>0$ and $\eta_0(L_\k)>0$ such that 
if $\epsilon\le \epsilon_2(L_\k)$, $L_{\k+1} \ge L_{\k+1,2}(L_\k,\epsilon)$ and $|x|\le a_2(L_\k, \epsilon, L_{\k+1}$) for some $\epsilon_2>0$, $L_{\k+1,2}>0$ and $a_2>0$, then
for all $s\in [s_{\k+1},\bar s]$,\\
(i) $\frac 1{C_*} \|q(s)\|_{\H}^2 \le h_1(s) \le C_* \|q(s)\|_{\H}^2$ and $h_1'(s) \le \frac 1{\delta^*} h_1(s)+\frac 1{\delta^*} J_\k(s)^2$ where $J_\k$ is defined in \aref{defjbar}.\\
(ii) $\left|(J_\k(s)^2)'\right|\le \d\frac 1{\delta^*}\left(\epsilon +\d\sum_{i=1}^\k\frac{|\nu_i(s)|}{1-{|d_i^*(s)|}}\right)J_\k(s)^2$ and $\left|(J_\k(s)^2)'\right|\le \frac 1{\delta^*} J_\k(s)^2$.\\
(iii) If
\begin{equation}\label{cond}
\sum_{i=1}^\k \frac{|\nu_i(s)|}{1-|d_i^*(s)|}\le \eta_0,
\end{equation}
then $\frac 1{C_*} \|q(s)\|_{\H}^2 \le h_2(s) \le C_* \|q(s)\|_{\H}^2$ and $h_2'(s) \le -\delta^*h_2(s)+\frac 1{\delta^*} J_\k(s)^2$.
\end{cl}
{\it Proof}: The proof is similar to our analogous estimate in Section 3.2 of \cite{MZajm10}. For that reason, we give it in Section \ref{sub2} in Appendix \ref{appproj9}, stressing only the differences with respect to \cite{MZajm10}.\Box

\bigskip

We now claim the following:
\begin{cl}\label{cl1} If $L_\k$ is fixed large enough, $M_0=M_0(L_\k)>0$ is large enough, $\epsilon \le \epsilon_3(L_\k)$, $\lzero\ge L_{\k+1,3}(L_\k,\epsilon)$ and $|x|\le a_3(L_\k,\epsilon, \lzero)$ for some $\epsilon_3>0$, $L_{\k+1,3}>0$ and $a_3>0$, then:\\
(i)
\[
\forall s\in[s_{\k+1}, \min(\bar s, s_\k-M_0)],\;\;
\d\sum_{i=1}^\k \frac{|\nu_i(s)|}{1-|d_i^*(s)|}\le \eta_0
\]
 where $\eta_0$ is given in Claim \ref{lemlyap}.\\
(ii) If we fix
\begin{equation}\label{fixb}
\beta=\beta(L_\k) =\frac 12\min\left(1, \frac 1{C_*}, \d\frac{e^{-\frac{M_0}{\delta^*}}}{4C_*^2}\right) 
\end{equation}
where the different constants appear in 
Claim \ref{lemlyap}, then for all $s\in[s_{\k+1}, \bar s]$ and $i=1,...,\k-1$, we have 
\begin{equation}\label{obj}
\begin{array}{rl}
&\zeta_{i+1}^*(s) - \zeta_i^*(s)\ge \frac{(p-1)}8 \log \lzero,\;\;
-1+\frac 34e^{-(p-1)L_\k}\le \frac{\nu_i(s)}{1-|d_i(s)|}\le C_0e^{L_\k}+\frac 12,\\
\mbox{ and }&
\|q(s)\|_{\H}\le \frac{\epsilon}{\sqrt 2}
\end{array}
\end{equation}
where $C_0$ appears in \aref{apriori}.
\end{cl}
{\it Proof}:\\
(i) See Appendix \ref{apptrans}.\\ 
(ii) Note that the two first estimates of \aref{obj} directly follow from Claims \ref{clprep} and \ref{clproche}, if
$L_\k$ is fixed large enough, $\epsilon \le \epsilon_4(L_\k)$, $\lzero\ge L_{\k+1,4}(L_\k,\epsilon)$ and $|x|\le a_4(L_\k, \epsilon, \lzero)$ for some $\epsilon_4>0$, $L_{\k+1,4}>0$ and $a_4>0$.
Thus, we only focus on the third estimate which is more delicate.\\
Consider $M_0>0$ to be fixed later in terms of $L_\k$.
If $s\in[s_{\k+1}, \min(\bar s, s_\k-M_0)]$, then we see from (i) of this claim, (ii) and (iii) of Claim \ref{lemlyap} that for $\epsilon$ and $\eta_0$ small enough, we have
\[
\left(h_2-\frac 2{{\delta^*}^2}J_\k^2\right)'\le -\delta^*\left(h_2-\frac 2{{\delta^*}^2}J_\k^2\right).
\]
Integrating this inequality and using (iii) of Claim \ref{lemlyap}, we write
\[
\|q(s)\|_{\H}^2\le C_*^2e^{-\delta^*(s-s_{\k+1})}\|q(s_{\k+1})\|_{\H}^2+2\frac{C_*}{{\delta^*}^2} J_\k(s)^2.
\]
Using the definition \aref{defjbar} of $J_\k$, \aref{small0} and \aref{gap},
we write
\[
\|q(s)\|_{\H}^2\le C_*^2e^{-\delta^*(s-s_{\k+1})}\beta^2\epsilon^2+2\frac{C_*}{{\delta^*}^2} (k-1)^2
\lzero^{-1/9}\le 2C_*^2\beta^2 \epsilon^2\le \frac{\epsilon^2}2, 
\]
provided that $\lzero=\lzero(L_\k,\epsilon)$ is large enough and $\beta=\beta(L_\k)$ is fixed by \aref{fixb}. Thus, we get the conclusion when $s\in[s_{\k+1}, \min(\bar s, s_\k-M_0)]$.\\
Now, if $\bar s>s_\k-M_0$ and $s\in[s_\k-M_0, \bar s]$, we have just proved that
\begin{equation}\label{dur}
\|q(s_\k-M_0)\|_{\H}^2 \le 2C_*^2\beta^2 \epsilon^2.
\end{equation}
Using (i) and (ii) of Claim \ref{lemlyap}, we write
\[
(h_1+J_\k^2)'\le \frac 2{\delta^*}(h_1+J_\k^2).
\]
Integrating this between $s_\k-M_0$ and $\bar s$, we write from (i) of Claim \ref{lemlyap}, 
\begin{equation}\label{duro}
\|q(\bar s)\|_{\H}^2 \le C_* e^{\frac 2{\delta^*}(\bar s-s_\k+M_0)}\left(C_* \|q(s_\k-M_0)\|_{\H}^2+J_\k(s_\k-M_0)^2\right).
\end{equation}
Since $\bar s\le s_m$ by definition of $\bar s$, we get from \aref{duro}, \aref{dur}, the definition \aref{defjbar} of $J_\k$, \aref{gap} and \aref{fixb},
\[
\|q(\bar s)\|_{\H}^2 \le C_* e^{\frac{2M_0}{\delta^*}}\left(2C_*^3\beta^2\epsilon^2+(k-1)^2 \lzero^{-1/9}\right)
\le 4C_*^4\beta^2e^{\frac{2M_0}{\delta^*}}\epsilon^2\le \frac{\epsilon^2}2,
\]
provided that $\lzero=\lzero(L_\k,\epsilon)$ is large enough. 
This concludes the proof of Claim \ref{cl1}.\Box

\bigskip

In conclusion, for all $\epsilon>0$, $s\in[s_{\k+1}, \bar s]$ and $i=1,...,\k-1$, we have just proved in Claim \ref{cl1} that
\[
\zeta_{i+1}^*(s) - \zeta_i^*(s)\ge \frac{(p-1)}8 \log \lzero,\;\;
-1+\frac 34e^{-(p-1)L_\k}\le \frac{\nu_i(s)}{1-|d_i(s)|}\le C_0e^{L_\k}+\frac 12,
\]
and 
\[
\|q(s)\|_{\H}\le \frac{\epsilon}{\sqrt 2}
\]
for well chosen parameters.
From the alternative \aref{dicho}, we see that
\[
\bar s = s_\k.
\]
Using (ii) of Claim \ref{clproche}, we see that for all $s\in [s_{\k+1}, s_\k]$, we have 
\[
\left\|\vc{w_{x}(s)}{\partial_s w_{x}(s)}-\d\sum_{i=1}^\k (-1)^i\kappa^*(\cd_i(s),\tmu_i(s))\right\|_{\H} \le C(L_\k)\epsilon.
\]
This concludes the proof of Proposition \ref{propall}.\Box

\subsection{Proof of Theorems \ref{thfini} and \ref{pdes}}\label{sub2*}
This section is devoted to the proof of the following result which directly implies Theorems \ref{thfini} and \ref{pdes}:
\begin{prop}[More properties of $\SS$]\label{propmore}
There exists $\delta_0>0$ such that if $-\delta_0 \le x<0$,
then for some $C_0>0$, we have:
\begin{eqnarray*}
(i)&& E(w_\tx(\ts))<2 E(\kappa_0)\mbox{ for some }\ts \ge - \log T(\tx),\\
(ii)&&\frac 1{C_0|\log |\tx||^{\gamma_1}} \le T'(\tx)-1 \le \frac{C_0}{|\log |\tx||^{\gamma_1}}.
\end{eqnarray*}
\end{prop}
Let us first derive Theorems \ref{thfini} and \ref{pdes} from this proposition and then prove it.

\bigskip

Consider $u(x,t)$ a blow-up solution of equation \aref{equ} and $x_0\in \SS$. Since equation \aref{equ} is invariant under space and time translation as well as symmetry in space, we may assume that $x_0=0$ and $T(x_0)=0$ and focus on the case $x<0$. 

\medskip

{\it Proof of Theorem \ref{thfini} assuming that Proposition \ref{propmore} holds}:
Applying Proposition \ref{propmore}, we see that for some $\delta_0>0$, if $-\delta_0\le x<0$, then $E(w_\tx(\ts))<2 E(\kappa_0)$ for some $\ts \ge -\log T(\tx)$. 
Using (i) of Proposition \ref{threg}, we see that $\tx \in \RR$. Since the same property holds for $x>0$ by symmetry, we see that $0$ is an isolated characteristic point and Theorem \ref{thfini} follows from Proposition \ref{propmore}.\Box

\bigskip

{\it Proof of Theorem \ref{pdes} assuming that Proposition \ref{propmore} holds}: 
It follows from (ii) of Proposition \ref{propmore} and the definition \aref{equid**} of $\gamma_1$.

\medskip

It remains to prove Proposition \ref{propmore}. 

\medskip

{\it Proof of Proposition \ref{propmore}}:

We work with the formulation given in the beginning of Section \ref{secfini} and in Section \ref{sub1*}.
 We proceed in 3 steps: In Step 1, we consider $x<0$ and give an upper bound on the energy of $w_{\tx}(\tsa)$ where $\tsa$ is defined in \aref{deftsa}. Using the lower bound of \aref{corcriterion}, we show that $\tsa = l+L_\hk$ in \aref{deftsa}, if $L_\hk$ is large enough and $|\tx|$ is small enough. Then, in Step 2, we conclude the proof of Proposition \ref{propmore}.

\bigskip

{\bf Step 1: An upper bound on $E(w_{\tx}(\tsa))$ for small $|\tx|$ and applications}

This step is dedicated to the proof of the following:
\begin{cl}\label{lemup} {\bf (An upper bound on $E(w_{\tx}(\tsa))$ for small $|\tx|$ and applications)} There exists $L_{\hk,1}>0$ 
such that for all $L_\hk\ge L_{\hk,1}$ and for all $\epsilon^*>0$, there exists $\delta(L_\hk, \epsilon^*)>0$ such that for all $\tx \in [-\delta,0)$, we have:\\
(i) 
$E(w_{\tx}(\tsa))\le E(\kappa_0)\left(\frac{p+1}{p-1}\bl1^2 - \frac 2{p-1} \bl1^{p+1}+\epsilon^*\bl1^{p+1}\right)+C(\epsilon^*+\sum_{i=2}^{\hat k}\bli)(1+\bl1^p)$,
where
\begin{equation}\label{defli*}
\bli = \lambda(\cd_i(\tsa), \tmu_i(\tsa))
\end{equation}
and $\lambda(d,\nu)$ is defined in \aref{deflambda}.\\
(ii) $\tsa=l+L_\hk$,\\
(iii) For all $i=2,...,\hk$, $\bli \le Ce^{-\frac {L_\hk}{p-1}}$ and $|\bl1-1|^2\le C(\epsilon^*+e^{-\frac {L_\hk}{p-1}})$.
\end{cl}
{\bf Remark}: If $k=2$, then $\hat k=\d E\left(\frac {k+1}2\right)=1$ and $\d\sum_{i=2}^{\hat k}\bli=0$.\\
{\it Proof}: Consider $L_\hk>0$ and $\epsilon^*\in (0, 1]$.

\medskip

(i) From Proposition \ref{propall} applied with $\k=\hk$,
we consider $\delta_1=\delta_1(L_\hk,\epsilon^*)>0$ such that for all $\tx \in [-\delta_1,0)$, we have
\begin{equation}\label{corinne}
\left\|\vc{w_{\tx}(\ty,\tsa)}{\ps w_{\tx}(\ty,\tsa)}- \sum_{i=1}^{\hat k}(-1)^{i}\kappa^*\left(\cd_i(\tsa), \tmu_i(\tsa), \ty\right)\right\|_{\H}\le \epsilon^*.
\end{equation}
Since we have for some $\delta_2=\delta_2(L,\epsilon^*)>0$ and for all $\tx \in [-\delta_2,0]$ (see below for a proof of this fact), 
\begin{equation}\label{fatigue1}
\frac{|\tmu_1(\tsa)|}{\sqrt{1-\cd_1(\tsa)^2}}\le \epsilon^*
\mbox{ and }\forall i=1,...,\hat k,\;\;
\|\kappa^*(\cd_i(\tsa), \tmu_i(\tsa), \cdot)\|_{\H} \le C\bli
\end{equation}
with
\begin{equation}\label{bli}
\bli\le 1\mbox{ when }2\le i \le \hat k,
\end{equation}
we write from \aref{corinne}
\begin{equation}\label{olivia}
\left\|\vc{w_{\tx}(\ty,\tsa)}{\ps w_{\tx}(\ty,\tsa)}+\kappa^*(\cd_1(\tsa),\tmu_1(\tsa), \ty)\right\|_{\H}\le \epsilon^*+C\sum_{i=2}^{\hat k}\bli.
\end{equation}
Applying (iii) of Claim \ref{lemsobolev-9}, we write from \aref{fatigue1}, \aref{bli} and \aref{olivia}
\[
E(w_{\tx}(\tsa))\le E(\kappa^*(\cd_1(\tsa),\tmu_1(\tsa), \cdot))+C(\epsilon^*+\sum_{i=2}^{\hat k}\bli)(1+\bl1^p).
\]
Using (i) of Lemma \ref{lemboundk*} and \aref{fatigue1}, we get 
(i) of Claim \ref{lemup} with $\delta=\min(\delta_1, \delta_2)$. It remains to prove \aref{fatigue1} and \aref{bli}.

\medskip

{\it Proof of \aref{fatigue1} and \aref{bli}}: Note first from \aref{deftsa}, \aref{defs2}, \aref{equid**} and \aref{defs2*} that 
\begin{equation}\label{ka1}
\htsa\to \infty\mbox{ and }\tsa \to \infty\mbox{ as }\tx\to 0^-,
\end{equation}
with (use \aref{chap**})
\begin{equation}\label{yol0}
\tsa \le l+L_\hk,\mbox{ and }\htsa\le-\log(-\tx(1-\tB))\le l+C
\end{equation}
for $|\tx|$ small enough. Therefore, using \aref{equid**}, we see that
\begin{equation}\label{yol}
\forall i=1,...,\hat k,\;\;\frac C{l^{\gamma_i}}\le 1-\cd_i(\tsa)\to 0\mbox{ as }\tx \to 0^-.
\end{equation}
Using the bound \aref{chap**} on $\tB$, we write from the definition \aref{defs2} of $\tmu_i(\tsa)$, \aref{yol0} and \aref{yol}
\begin{eqnarray}
\frac{|\tmu_1(\tsa)|}{\sqrt{1-\cd_1(\tsa)^2}}&\le& \left(\frac{\tB}{\sqrt{1-\cd_1(\tsa)}}+\sqrt{1-\cd_1(\tsa)}\right)e^{L_\hk}\nonumber\\
 &\le& \left(\frac C{l^{\gamma_1/2}}+\sqrt{1-\cd_1(\tsa)}\right)e^{L_\hk} \to 0\label{ande}
\end{eqnarray}
as $\tx \to 0$, and the first part of \aref{fatigue1} is proved.\\
Since we see from \aref{deftsa} and \aref{defs2} that
\begin{equation}\label{pam93}
1+\tmu_1(\tsa) -\cd_1(\tsa)\ge e^{-(p-1)L_\hk}(1-\cd_1(\tsa)),
\end{equation}
we write from the definitions \aref{deflambda} and \aref{defli*} of $\lambda$ and $\bli$, \aref{yol} and \aref{ande}, for $|\tx|$ small enough with $\tx<0$, 
\[
\bl1^{p-1} = \frac{1-\cd_1(\tsa)}{1+\tmu_1(\tsa) - \cd_1(\tsa)}\cdot
\frac{1+\cd_1(\tsa)}{1+\tmu_1(\tsa) + \cd_1(\tsa)}
\le e^{(p-1)L_\hk}\cdot 2^{p-1},
\]
hence
\[
\bl1 \le2 e^{L_\hk}.
\]
Taking $|\tx|$ small enough with $\tx <0$, we see from (i) of Lemma \ref{lemboundk*} and \aref{ande} that the second identity in \aref{fatigue1} holds for $i=1$.\\
When $2\le i \le \hat k$, using the bound \aref{chap**} on $\tB$ and \aref{yol0}, we see from \aref{yol} and by definition \aref{defs2} that $\tmu_i(\tsa)\ge 0$, hence by definition \aref{deflambda}, $\bli\le 1$ and \aref{bli} follows. Applying (i) of Lemma \ref{lemboundk*}, we see that the second identity in \aref{fatigue1} holds for $2\le i \le \hat k$. This concludes the proof of \aref{fatigue1}, \aref{bli} and the proof of (i) in Claim \ref{lemup} too.

\medskip

(ii) Taking $\epsilon^*\le \frac 1{p-1}$ and using the fact that $\bli\le 1$ when $2\le i \le \hat k$ (see \aref{bli}), we write from 
(i) of Claim \ref{lemup}
and 
\aref{corcriterion}
\[
E(\kappa_0)\le E(\kappa_0)\left(\frac{p+1}{p-1}\bl1^2 - \frac 1{p-1} \bl1^{p+1}\right)+C(1+\bl1^p)\equiv h_0(\bl1).
\]
Since $\bl1>0$ by the definition \aref{deflambda} and
\[
h_0(\lambda)\to -\infty\mbox{ as }\lambda \to \infty,
\]
it follows that for some $C_0>0$,
\begin{equation}\label{bound}
0<\bl1\le C_0
\end{equation}
on the one hand.\\
On the other hand, using \aref{fatigue1} and \aref{yol},
we see from \aref{deflambda} and \aref{defs2'} that if Case 2 occurs, then 
\[
\bl1^{p-1} = \frac{1-\cd_1(\tsa)}{1+\tmu_1(\tsa) - \cd_1(\tsa)}\cdot
\frac{1+\cd_1(\tsa)}{1+\tmu_1(\tsa) + \cd_1(\tsa)}
\ge e^{(p-1)L_\hk}\cdot \frac 1{2^{p-1}},
\]
hence
\begin{equation}\label{lina1}
\bl1 \ge \frac{e^{L_\hk}}2.
\end{equation}
for $|\tx|$ small enough. If $L_\hk\ge L_{\hk,1}$ where we fix
\[
e^{L_{\hk,1}}=4C_0
\]
and $C_0$ is given in \aref{bound}, then we see from \aref{bound} and \aref{lina1} that Case 2 given in \aref{defs2'} cannot occur. Therefore, Case 1 (given in \aref{deftsa}) holds and $\tsa=l+L_\hk$.
Thus, (ii) of Claim \ref{lemup} is proved.  

\medskip

(iii) Since $\tsa=l+L_\hk$ by (ii) of Claim \ref{lemup}, we see from Claim \ref{clvanish} and the definition \aref{defli*} of $\bli$ that
\begin{equation}\label{mpetit}
\bli \le Ce^{-\frac {L_\hk}{p-1}}\mbox{ when }2\le i \le \hat k.
\end{equation}
Using \aref{mpetit}, 
(i) of Claim \ref{lemup}, \aref{corcriterion} and \aref{bound}, we see that
\begin{equation}\label{battikh}
E(\kappa_0) - C(\epsilon^*+ e^{-\frac {L_\hk}{p-1}})\le E(\kappa_0)\left(\frac{p+1}{p-1}\bl1^2 - \frac 2{p-1} \bl1^{p+1}\right)\equiv E(\kappa_0)h(\bl1).
\end{equation}
Since the maximum of $h$ is $h(1)=E(\kappa_0)$ and since it is achieved only for $\bl1=1$, we get from the fact that $h'(1)=0$ and $h"(1)\neq 0$,
\begin{equation}\label{petit}
|\bl1-1|^2 \le C(\epsilon^*+ e^{-\frac {L_\hk}{p-1}}).
\end{equation}
This concludes the proof of Claim \ref{lemup}.\Box

\bigskip

{\bf Step 2: Conclusion of the proof of Proposition \ref{propmore}}:

Using (i) of Lemma \ref{lemboundk*}, \aref{fatigue1} and \aref{petit}, we see that \aref{olivia} and \aref{mpetit} yield
\begin{equation}\label{trap*}
\left\|\vc{w_{\tx}(\ty,\tsa)}{\ps w_{\tx}(\ty,\tsa)}+\vc{\kappa\left(\cd^*_1(\tsa),y\right)}{0}\right\|_{\H} \le C_1(\sqrt{\epsilon^*}+ e^{-\frac {L_\hk}{2(p-1)}})
\end{equation}
for some $C_1>0$, where $\cd^*_1(\tsa)=\frac{\cd_1(\tsa)}{1+\tmu_1(\tsa)}=\frac{\cd_1(\tsa)}{1+\tmu_1(\tsa)}$.
From (i), (ii) and (iii) of Claim \ref{lemsobolev-9}, we see that
\begin{equation}\label{elton}
E(w_{\tx}(\tsa)) \le E(\kappa_0)+C_2(\sqrt{\epsilon^*}+ e^{-\frac L{2(p-1)}})
\end{equation}
for some $C_2>0$. Fixing 
\begin{equation}\label{defL}
L_\hk=\max\left(L_{\hk,1},2(p-1)\left(|\log\left|\frac{\epsilon_0}{2C_1}\right|\right|, 2(p-1))\left|\log\left|\frac{E(\kappa_0)}{4C_1}\right|\right|\right)
\end{equation}
where $\epsilon_0$ appears in the trapping result of \cite{MZjfa07} stated in (ii) of Proposition \ref{threg}
and fix $\epsilon^*$ small enough so that
\begin{equation}\label{defe*}
\epsilon^*\le \min\left(\frac 1{p-1},\left(\frac{\epsilon_0}{2C_1}\right)^2, \left(\frac{E(\kappa_0)}{4C_1}\right)^2\right),
\end{equation}
we see from \aref{elton} that
\[
E(w_{\tx}(\tsa))\le E(\kappa_0) + \epsilon_0 \le \frac 32E(\kappa_0),
\]
(note that we also have
\begin{equation}\label{piege}
\left\|\vc{w_{\tx}(\ty,\tsa)}{\ps w_{\tx}(\ty,\tsa)}+\vc{\kappa\left(\cd^*_1(\tsa),y\right)}{0}\right\|_{\H} \le \epsilon_0),
\end{equation}
which concludes the proof of (i) of Proposition \ref{propmore}.

\bigskip

(ii) Using (ii) of Proposition \ref{threg}, we see 
that for $|\tx|$ small enough, we have for some $C_0>0$,
\begin{equation}\label{planchon}
\left|\argth(T'(\tx))-\argth\left(\cd^*_1(\tsa)\right)\right|\le C_0\epsilon_0
\end{equation}
Since we have from \aref{pam93}, (ii) of Claim \ref{lemup} and the positivity of $\tB$ (see \aref{chap**})
\[
 e^{-(p-1)L_\hk}(1-\cd_1(\tsa))\le 1+\tmu_1(\tsa) -\cd_1(\tsa)\le (1-\cd_1(\tsa))(e^{L_\hk}+1),
\]
 and from \aref{ka1}, \aref{equid**} and \aref{fatigue1}
\[
\frac 1{C l^{\gamma_1}}\le 1-\cd_1(\tsa) \le \frac C{l^{\gamma_1}},\;\;\tmu_1(\tsa) \to 0\mbox{ as }\tx \to 0^-,
\]
we write from the definition \aref{defs2} of $\tmu_1(\tsa)$,
\begin{eqnarray*}
\argth\left(\cd^*_1(\tsa)\right)&=&\argth\left(\frac{\cd_1(\tsa)}{1+\tmu_1(\tsa)}\right)= 
\frac 12 \log\left(\frac{1+\tmu_1(\tsa) +\cd_1(\tsa)}{1+\tmu_1(\tsa) -\cd_1(\tsa)}\right)\\
&=&-\frac 12 \log\left(1-\cd_1(\tsa)+\tmu_1(\tsa)\right)+O(1)\\
&=& - \frac 12 \log(1-\cd_1(\tsa))+O(1)= \frac{\gamma_1}2 \log l +O(1)
\end{eqnarray*}
as $\tx \to 0^-$. Using \aref{planchon}, we see that
\[
T'(\tx) =-\tanh(\frac{\gamma_1}2 \log l +O(1))=-1 +l^{-\gamma_1}e^{O(1)}\mbox{ as }\tx \to 0^-
\]
and (ii) follows. This concludes the proof of Proposition \ref{propmore} as well as Theorems \ref{thfini} and \ref{pdes}.\Box

\subsection{The blow-up speed in $L^\infty$ near a characteristic point}\label{secspeed}
We prove Corollary \ref{corspeed} here. Note that it is a consequence of our anaysis in this section. In particular, the lower bound on $T(x)$ in \aref{tx} is crucial for our argument.

\medskip

{\it Proof of Corollary \ref{corspeed}}: Consider $x_0\in\SS$. As in section \ref{secfini}, we assume from translation invariance of equation \aref{equ} that $x_0=T(x_0)=0$, hence,
\[
0\in\SS\mbox{ and }T(0)=0.
\]
We will also use the notation $W(Y,S)$ instead of $w_0(y,s)$ defined in the selfsimilar transformation \aref{defw}. 
Up to replacing $u$ by $-u$, we see from Proposition \ref{thsing} that \aref{cprofile**} and \aref{equid**} hold for some integer $k=k(0)\ge 2$, for some continuous functions $D_i(S)$, $C_0>0$ and $S_0\in\R$. 
Using the Hardy-Sobolev estimate of \aref{hs} and introducing
\begin{equation}\label{defbw}
\bar W(\xi,S) = (1-Y^2)^{\frac 1{p-1}}W(Y,S),\;\;Y=\tanh \xi\mbox{ and }\zeta_i(S)=-\argth D_i(S),
\end{equation}
we see that
\begin{equation}\label{barprofile}
\left\|\bar W(\xi,S)-\sum_{i=1}^k(-1)^i
\kappa_0 \cosh^{-\frac 2{p-1}}(\xi-\zeta_i(S))
\right\|_{L^\infty} \to 0\mbox{ as }S\to \infty
\end{equation}
From \aref{defw}, our goal that for $S$ large enough,
\begin{equation}\label{but}
\frac{S^{\frac{k-1}2}}C \le \|W(S)\|_{L^\infty(-1,1)}\le CS^{\frac{k-1}2}.
\end{equation}
We first prove the lower bound, then the upper bound.

\medskip

{\it - The lower bound}: Since
\[
\|W(S)\|_{L^\infty(-1,1)}\ge |W(-D_1(S),S)|
\]
and 
from \aref{equid**}
\begin{equation}\label{estD1}
\frac{S^{-\frac{(p-1)(k-1)}2}}C \le 1-D_1(S)\le CS^{-\frac{(p-1)(k-1)}2}
\end{equation}
for $S$ large enough, we write 
from \aref{barprofile}, \aref{defbw} and the definition \aref{defkd} of $\kappa_d$,
\[
|\bar W(\zeta_1(S),S)| \ge \frac{\kappa_0}2,\mbox{ hence }|W(-D_1(S),S)|\ge \frac 12 \kappa_{D_1(S)}(-D_1(S),S)\ge \frac{S^{\frac{k-1}2}}C,
\]
and the lower bound in \aref{but} follows.

\medskip

{\it - The upper bound}: The following uniform estimate is crucial for your argument. It follows from our analysis in the proof of Theorems \ref{thfini} and \ref{pdes}.
\begin{lem}\label{lemunif}{\bf (Uniform bound on $w_x(s)$)}If $|x|$ is small enough and $s$ is large enough, then  
\[
\|w_x(s)\|_{\H} +\|(1-y^2)w_x(s)\|_{L^\infty(-1,1)}\le C.
\]
\end{lem}
Let us first use this lemma to prove the upper bound in \aref{but}, then we will prove it.\\
Up to making the symmetry $x\mapsto -x$ in equation \aref{equ}, we may assume that 
\begin{equation}\label{neg}
-1<Y\le 0.
\end{equation}
Introducing 
\begin{equation}\label{defXXs}
X= -\frac{e^{-S}}4,
\end{equation}
we see from \aref{defw**} and \aref{rach1} that
\begin{equation}\label{defWW}
W(Y,S) = (1+(1-b)Xe^S)^{-\frac 2{p-1}}w_X(y,s)
\end{equation}
where 
\begin{equation}\label{change}
s=S-\log(1+(1-b)Xe^S),\;\;
y=\frac{Y-Xe^S}{1+(1-b)Xe^S},
\mbox{ and }b=1-\frac{T(X)}X.
\end{equation}
Since we have from \aref{change} and \aref{defXXs} for $S$ large enough,
\[
\frac{|\log|X||^{-\frac{(p-1)(k-1)}2}}C\le b\le C|\log|X||^{-\frac{(p-1)(k-1)}2},\mbox{ hence }
\frac{S^{-\frac{(p-1)(k-1)}2}}C\le b\le CS^{-\frac{(p-1)(k-1)}2}
\]
by Theorem \ref{pdes} and 
\begin{eqnarray}
1+(1-b)Xe^S &=& \frac 34 +\frac b4 = \frac 34+O(S^{-\frac{(p-1)(k-1)}2}),\nonumber\\
y &=&\frac{4Y+1}{3+b}\le \frac 13,\nonumber\\
1+y &=& \frac{4+4Y+b}{3+b} \ge \frac b4 \ge \frac{S^{-\frac{(p-1)(k-1)}2}}C,\nonumber
\end{eqnarray}
by \aref{defXXs}, \aref{change}, \aref{neg} and \aref{estD1}, we see from 
\aref{defWW} and Lemma \ref{lemunif} that
\[
|W(Y,S)|\le C|w_X(y,s)|\le C(1-y^2)^{-\frac 1{p-1}}\le C(1+y)^{-\frac 1{p-1}} \le CS^{\frac{k-1}2}
\]
which is the desired conclusion. It remains to prove Lemma \ref{lemunif} in order to conclude the proof of Corollary \ref{corspeed}.

\medskip

{\it Proof of Lemma \aref{lemunif}}: The bound on the weighted $L^\infty$ space follow from the bound in $\H$ thanks to the Hardy-Sobolev estimate \aref{hs}. Thus, we only prove the bound in $\H$.\\
- When $x=0$, the result follows from 
Proposition 3.5 page 66 in \cite{MZjfa07}. Of course, the reader may see this bound as a consequence of the decomposition into solitons \aref{cprofile**} (which is much stronger) and the boundedness of the solitons in $\H$ stated in (i) of Claim \ref{lemsobolev-9}.\\
- When $x\neq 0$, 
up to replacing $x$ by $-x$ in equation \aref{equ}, we may assume that $x<0$. Therefore, the previous part of Section \ref{secfini} applies.\\ 
If $s\in[L_{k+1},\tsa]$, 
then the result follows from the application of Proposition \ref{propall}, for all $\k\in[k,\hk]$ and (i) of Claim \ref{lemsobolev-9}, provided that the different constants $L_{k+1}$,...,$L_{\hk+1}$ are large enough and $|x|$ is small enough.\\
If $s\ge \tsa$, then we recall from \aref{piege} that for $|x|$ small enough, we have
\begin{equation*}
\left\|\vc{w_{\tx}(\ty,\tsa)}{\ps w_{\tx}(\ty,\tsa)}+\vc{\kappa\left(\cd^*_1(\tsa),y\right)}{0}\right\|_{\H} \le \epsilon_0.
\end{equation*}
Thus, the trapping result of (ii) of Proposition \ref{threg} applies and we get for some positive $\mu_0$ and $C_0$ and for all $s\ge \tsa$:
\begin{equation*}
\left\|\vc{w_x(s)}{\partial_s w_x(s)}-\theta(x_0)\vc{\kappa(T'(x))}{0}\right\|_{\H}\le C_0 e^{-\mu_0(s-s_0)}.
\end{equation*}
Using 
again (i) of 
Claim \ref{lemsobolev-9}, we get the conclusion. This concludes the proof of Lemma \ref{lemunif} and Corollary \ref{corspeed} too.\Box

\appendix

\section{Properties of $\kappa(d,y)$ and $\kappa^*(d,\nu,y)$}\label{appprop}
In this section, we give some properties of $\kappa(d,y)$ and $\kappa^*(d,\nu,y)$ defined in \aref{defkd} and \aref{defk*}.
We first recall the following:
\begin{cl}\label{lemsobolev-9} 
(i) {\bf (Boundedness of $\kappa(d,y)$ in several norms)}
For all $d\in (-1,1)$, it holds that
\[
\|\kappa(d,y)\|_{L^{p+1}_\rho(-1,1)}+ \|\kappa(d,y)(1-y^2)^{\frac 1{p-1}}\|_{L^\infty(-1,1)}
\le C\|\kappa(d,y)\|_{\H_0}\le CE(\kappa_0).
\]
(ii) {\bf (Same energy level for $\kappa(d,y)$)} For all $d\in (-1,1)$, it holds that
\[
E(\pm\kappa(d,y))=E(\kappa_0).
\]
(iii) {\bf (Continuity of the Lyapunov functional)} If $(w_i(y,s), \ps w_i(y,s))\in \H$ for $i=1$ and $2$ and for some $s\in \R$, then 
\begin{eqnarray*}
&&|E(w_1(s))-E(w_2(s))|\\
&\le& C\left\|\vc{w_1(s)}{\ps w_1(s)}-\vc{w_2(s)}{\ps w_2(s)}\right\|_{\H}\left(1+\left\|\vc{w_1(s)}{\ps w_1(s)}\right\|_{\H}^p+\left\|\vc{w_2(s)}{\ps w_2(s)}\right\|_{\H}^p\right).
\end{eqnarray*}
\end{cl}
{\it Proof}: 
For (i), use identity (49) page 59 in \cite{MZjfa07} and the following Hardy Sobolev identity: 
\begin{equation}\label{hs}
\forall h\in \H_0,\;\;\|h\|_{L^2_{\frac \rho{1-y^2}(-1,1)}}+\|h\|_{L^{p+1}_\rho(-1,1)}+ \|h(1-y^2)^{\frac 1{p-1}}\|_{L^\infty(-1,1)}
\le C\|h\|_{\H_0}.
\end{equation}
For (ii), see (ii) of Proposition 1 page 47 in \cite{MZjfa07}. 
The proof of (iii) is straightforward from the definition \aref{defenergy} of $E(w)$.\Box

\bigskip

\noindent Now,
we give in the following lemma some  properties of $\kappa^*(d,\nu,y)$ defined in \aref{defk*} which are useful for the proof.
\begin{lem}[Properties of $\kappa^*(d,\nu,y)$]\label{lemboundk*} For all $d\in (-1,1)$ and $\nu >-1+|d|$, we have:\\
(i)
\begin{eqnarray*}
\forall y\in(-1,1),\;\;0\le\d\kappa_1^*(d,\nu,y)=\lambda \kappa\left(\frac d{1+\nu},y\right)&\le & \frac{\kappa_0 \lambda}{(1-y^2)^{\frac 1{p-1}}},\\
\d\left\|\kappa^*\left(d,\nu\right)\right\|_{\H} &\le& C\lambda +C1_{\{\nu<0\}}\frac{|\nu|}{\sqrt{1-d^2}}\lambda^{\frac{p+1}2},\\
\d\left\|\kappa^*\left(d,\nu\right)-\vc{\kappa\left(\frac d{1+\nu}\right)}{0}\right\|_{\H} &\le& C|\lambda-1|+C\frac{|\nu|}{\sqrt{1-d^2}}\lambda^{\frac{p+1}2}\\
E(\kappa_0)\ge E(\kappa^*(d,\nu))=E(\kappa_0)\left(\frac{p+1}{p-1}\lambda^2-\frac 2{p-1} \lambda^{p+1}\right.&+&\left.\frac 2{(p-1)}\frac{\nu^2}{(1-d^2)}\lambda^{p+1}\right)
\end{eqnarray*}
where $\lambda$ is defined in \aref{deflambda}.\\
(ii) For all $A\ge 2$, there exists $C(A)>0$ such that if $(d_1,\nu_1)$ and $(d_2, \nu_2)$ satisfy 
\begin{equation}\label{condA}
\frac {\nu_1}{1-|d_1|},\frac {\nu_2}{1-{|d_2|}}\in [-1+\frac 1A, A],\;\;
\end{equation}
then
\begin{equation}\label{defze}
\|\kappa^*(d_1,\nu_1)-\kappa^*(d_2, \nu_2)\|_{\H} \le C(A)\left(\left| \frac {\nu_1}{1-|d_1|} -\frac {\nu_2}{1-{|d_2|}}\right| +\left|\argth d_1 - \argth d_2\right|\right).
\end{equation}
Moreover, there exists $\epsilon_0(A)>0$ such that if $\|\kappa^*(d_1,\nu_1)-\kappa^*(d_2, \nu_2)\|_{\H_{1,2}}
\le \epsilon_0(A)$ where 
\begin{equation*}
\|V\|_{\H_{1,2}}^2 =\int_{\tanh(\argth(\max_{i=1,2}(-\frac {d_i}{1+\nu_i}))+A)}^1 \left(V_1^2+{V_1'}^2(1-y^2)+V_2^2\right)\rho dy,
\end{equation*}
then 
\begin{equation}\label{inf2}
\left| \frac {\nu_1}{1-|d_1|} -\frac {\nu_2}{1-|d_2|}\right| +\left|\argth d_1 - \argth d_2\right|\le  C(A)\|\kappa^*(d_1,\nu_1)-\kappa^*(d_2, \nu_2)\|_{\H_{1,2}}.
\end{equation}
\end{lem}
{\it Proof}: Consider $d\in (-1,1)$ and $\nu>-1+|d|$.

\medskip

(i) We first introduce the transformations
\begin{equation}\label{deftr}
r(y) \mapsto \bar r(\xi)= r(y)(1-y^2)^{\frac 1{p-1}}\mbox{ and }
r(y) \mapsto \hat r(\xi)=r(y)(1-y^2)^{\frac 1{p-1}+\frac 12}\mbox{ with }y=\tanh \xi
\end{equation}
and for $r=(r_1,r_2)$, the notation $\tilde r = (\bar r_1,\hat r_2)$.
We then recall the following estimate from \cite{MZajm10}:
\begin{cl}\label{cltrans}{\bf (Continuity of the transformation 
defined in \aref{deftr} and its inverse)} 
There exists $C_0>0$ such that for all $r\in \H$, we have
$\frac 1{C_0} \|r\|_{\H} \le \|\tilde r\|_{H^1\times L^2(\R)}\le C_0 \|r\|_{\H}$.
\end{cl}
{\it Proof}: See Claim B.2 in \cite{MZajm10}.\Box

\medskip

Using the transformation \aref{deftr} and definitions \aref{defkd} and \aref{defk*} of $\kappa(d,\cdot)$ and $\kappa^*$, we see that
\begin{equation}\label{67B}
\bar \kappa^*_1(d,\nu,\xi)=\lambda \bar \kappa_0(\xi-\zeta^*)\mbox{ and }
\bar \kappa\left(\frac d{1+\nu},\xi\right) = \bar \kappa_0(\xi - \zeta^*)
\end{equation}
where
\begin{equation}\label{defbk0} 
\bar \kappa_0(\xi)=\kappa_0 \cosh^{-\frac 2{p-1}}\xi,\;\;
\tanh \zeta^* = -\frac d{1+\nu}
\end{equation}
and $\lambda$ is defined in \aref{deflambda}.
Since $\|\bar \kappa_0\|_{L^\infty(\R)}\le \kappa_0$, the first identity follows from \aref{deftr}.\\
Since we have from the definition \aref{defk*} of $\kappa^*$, 
\[
\kappa^*_2(d,\nu,y)=-\sqrt{\frac 2{p+1}}\frac{\nu}{\sqrt{1-d^2}} \kappa^*_1(d,\nu, y)^{\frac{p+1}2},
\]
it follows from \aref{deftr} and \aref{67B} that
\begin{equation}\label{sami}
\widehat{\kappa^*_2}(d,\nu,\xi)=
-\sqrt{\frac 2{p+1}}\frac{\nu}{\sqrt{1-d^2}}\lambda^{\frac{p+1}2}\bar \kappa_0(\xi-\zeta^*)^{\frac{p+1}2}.
\end{equation}
Therefore, we have from \aref{67B} and \aref{sami},
\begin{eqnarray}
&&\|\bar \kappa^*_1(d,\nu,\cdot)\|_{H^1(\R)}\le C \lambda,\;\;\|\widehat{\kappa^*_2}(d,\nu,\cdot)\|_{L^2(\R)}\le C\frac{|\nu|}{\sqrt{1-d^2}}\lambda^{\frac{p+1}2}\label{penn}\\
\mbox{and}&&\left\|\bar \kappa^*_1(d,\nu,\cdot)-\bar \kappa\left(\frac d{1+\nu},\cdot\right)\right\|_{H^1(\R)}\le C |\lambda-1|.\label{penn0}
\end{eqnarray}
Since we have from the definition \aref{deflambda} of $\lambda$ and \aref{penn}, if $\nu\ge 0$ then
\[
\frac{\nu^2 \lambda^{p-1}}{1-d^2}=\frac{\nu^2}{(1+\nu)^2-d^2}\le 1,\mbox{ hence }\|\widehat{\kappa^*_2}(d,\nu,\cdot)\|_{L^2(\R)}\le C\lambda,
\]
the second and the third identities of (i) directly follow from  \aref{penn} and \aref{penn0} together with Claim \ref{cltrans}.

\medskip

For the left-hand side of the forth identity, using the fact that for $\mu=\pm1$, $\kappa^*(d,\mu e^s,\cdot)$ is a particular solution of equation \aref{eqw1} which goes to  $(\kappa(d,y),0)$ in $\H$ as $s\to -\infty$, we get the result from the monotonicity of $E$ and (ii) of Claim \ref{lemsobolev-9}. 

\medskip

For the right-hand side of the forth identity, using the definition \aref{defenergy} of the Lyapunov functional $E(w)$ and integration by parts, we see that
\[
E(\kappa^*) = \frac 12 \iint (\kappa^*_2)^2 \rho + \frac 12 \iint \kappa^*_1\left(-\L \kappa^*_1 + \frac{2(p+1)}{(p-1)^2} \kappa^*_1-\frac 2{p+1} {\kappa^*_1}^p\right)\rho.
\]
Performing the change of variables \aref{deftr} and using the fact that 
\[
\partial^2_\xi \bar \kappa^*_1 - \frac 4{(p-1)^2} \bar \kappa^*_1 +\frac 2{p+1}{(\bar \kappa^*_1)}^p = (1-y^2)^{\frac p{p-1}}\left[\L \kappa^*_1 - \frac{2(p+1)}{(p-1)^2} \kappa^*_1+\frac 2{p+1} {\kappa^*_1}^p\right]
\]
(see pages 59 and 60 in \cite{MZjfa07} for a proof of this fact), we see from the transformation \aref{deftr} that
\begin{eqnarray*}
E(\kappa^*)&=&\frac 12 \int_{\R} (\widehat{\kappa^*_2})^2 + \frac 12 \int_{\R} \bar \kappa^*_1\left(-\partial^2_\xi \bar \kappa^*_1 + \frac 4{(p-1)^2} \bar \kappa^*_1 -\frac 2{p+1}(\bar \kappa^*_1)^p\right)\\
&=& \int_{\R} \left( \frac 12 (\widehat{\kappa^*_2})^2 +  \frac 12 (\partial_\xi \bar\kappa^*_1)^2 +\frac 2{(p-1)^2} (\bar \kappa^*_1)^2 -\frac 1{p+1}(\bar \kappa^*_1)^p\right)
\end{eqnarray*}
where $\widehat{\kappa^*_2}$ and $\bar \kappa^*_1$ are given in \aref{sami} and \aref{67B}. Therefore,
\begin{equation}\label{zouhour}
E(\kappa^*(d,\nu,\cdot))=\beta_1 \lambda^2-\beta_2\lambda^{p+1}+ \beta_2 \frac{\nu^2}{1-d^2}\lambda^{p+1},
\end{equation}
for some $\beta_1>0$ and 
\begin{equation}\label{bb1}
\beta_2 = \frac 1{p+1}\int_{\R} \bar \kappa_0(\xi)^{p+1}d\xi = \frac{1}{p+1}\kappa_0^{p+1}\iint \rho(y) dy
\end{equation}
by the change of variables $y=\tanh \xi$. 
Since $\kappa^*(0,0,y)=(\kappa_0,0)$ and $\lambda =1$ when $d=\nu=0$, we see from \aref{zouhour} that 
\begin{equation}\label{b3}
E(\kappa_0) = \beta_1-\beta_2.
\end{equation}
Since 
\[
E(\kappa_0) = \left(\frac{(p+1)}{(p-1)^2}\kappa_0^2 - \frac{1}{p+1}\kappa_0^{p+1}\right)\iint \rho(y) dy=\frac{p-1}{2(p+1)}\kappa_0^{p+1}\iint \rho(y) dy
\]
from the definitions \aref{defenergy} and \aref{defkd} of $E(\kappa_0)$ and $\kappa_0$, we see from \aref{bb1} and \aref{b3} that $\beta_2 = \frac 2{p-1}E(\kappa_0)$ and $\beta_1 = \frac{p+1}{p-1} E(\kappa_0)$ and the result follows from \aref{zouhour}. This concludes the proof of (i) of Lemma \ref{lemboundk*}.

\bigskip

(ii) Take $A\ge 2$ and assume that \aref{condA} holds. Introducing the following new variables
\begin{equation}\label{defze0}
\zeta=-\argth d,\;\;\mu=\frac \nu{1-|d|},\;\;
\eta =\frac \nu{1-d^2},
\end{equation}
and for $i=1$ and $2$, $\zeta_i=\zeta(d_i,\nu_i)$, $\mu_i=\mu(d_i,\nu_i)$, $\eta_i=\eta(d_i,\nu_i)$ as well as $\Bli=\lambda(d_i,\nu_i)$ defined in \aref{deflambda},
we see from \aref{condA} and \aref{deflambda} that for $i=1$ and $2$,
\begin{equation}\label{condA2}
\mu_i,\eta_i\in [-1+\frac 1{A}, A],\;\;\Bli \in [(1+A)^{-\frac 2{p-1}}, A^{\frac 2{p-1}}]
\end{equation}
and 
\begin{equation}\label{lip}
(\zeta,\mu)\mapsto (\zeta,\eta)\mbox{ is Lipschitz and so is its inverse.}
\end{equation}
Due to this fact and to the fact that $\eta$ is $C^1$ as a function of $(d,\nu)$, unlike $\mu$, we will prove \aref{defze} and \aref{inf2} with the variable $\eta$ instead of $\mu$.\\
We first prove the upper bound \aref{defze} then the lower bound \aref{inf2}.

\medskip

{\it The upper bound \aref{defze} on $\|\kappa^*(d_1,\nu_1) - \kappa^*(d_2,\nu_2)\|_{\H}$}: 
From Claim \ref{cltrans}, \aref{67B}, \aref{sami} and \aref{condA2}, we see that
\begin{equation}\label{chu}
\|\kappa^*(d_1,\nu_1) - \kappa^*(d_2,\nu_2)\|_{\H} \le C(A)\left(|\zeta^*_1-\zeta^*_2|+|\Bl1 - \lambda_2|+\left|\frac {\nu_1}{\sqrt{1-d_1^2}}- \frac {\nu_2}{\sqrt{1-d_2^2}}\right|\right)
\end{equation}
where for $i=1$ and $2$,
\begin{equation}\label{defz0}
\zeta_i^*=-\argth\left(\frac {d_i}{1+\nu_i}\right).
\end{equation}
Since $\frac \nu{\sqrt{1-d^2}}=\eta \sqrt{1-d^2} = \eta \cosh^{-1} \zeta$ from \aref{defze0}, we write from \aref{condA2}
\begin{equation}\label{te}
\left|\frac {\nu_1}{\sqrt{1-d_1^2}}- \frac {\nu_2}{\sqrt{1-d_2^2}}\right|\le C(A)|\eta_1-\eta_2|+C(A)|\zeta_1-\zeta_2|.
\end{equation}
Therefore, clearly the upper bound follows from \aref{chu}, \aref{te} and \aref{defze0}, if we prove that under \aref{condA2}
\begin{equation}\label{jaco}
\|\Jac_{\zeta,\eta}(\zeta^*, \lambda)\|\le C(A).
\end{equation}
It remains then to prove \aref{jaco}.\\
 From the definitions \aref{defz0} and \aref{deflambda} of $\zeta^*$ and $\lambda$, we write
\begin{equation}\label{jac1}
\Jac_{d,\nu}(\zeta^*, \lambda)=
\frac 1{(1+\nu)^2-d^2}
\left(
\begin{array}{ll}
-(1+\nu)&d\\
\frac{2d\lambda^{2-p}(1-(1+\nu)^2)}{(p-1)[(1+\nu)^2-d^2]}&-\frac{2(1+\nu)\lambda^{2-p}(1-d^2)}{(p-1)[(1+\nu)^2-d^2]}
\end{array}
\right).
\end{equation}
Using \aref{jac} and \aref{condA2}, we see that
%
\begin{equation}\label{jac2}
\|\Jac_{d,\nu}(\zeta^*, \lambda)\|\le \frac {C(A)}{1-d^2}\mbox{ and }\|\Jac_{\zeta, \eta}(d,\nu)\|\le C(A)(1-d^2).
\end{equation}
By multiplication, we get \aref{jaco}. This concludes the proof of the upper bound in (ii).

\medskip

{\it The lower bound \aref{inf2} on $\|\kappa^*(d_1,\nu_1) - \kappa^*(d_2,\nu_2)\|_{\H}$}: Using the transformation \aref{deftr} and the embedding of $H^1(\xi>A)$ into $L^\infty(\xi>A)$, the notation \aref{defze0}, \aref{defz0} and symmetry, we see that it is enough to prove that for some $\epsilon_0(A)>0$, if for some $\zeta^*_1<0$, $\Bl1$ and $\lambda_2$ satisfying \aref{condA2}, we have
\begin{equation}\label{petit9}
E\equiv\|\Bl1 \bk(\cdot - \zeta^*_1)-\lambda_2 \bk\|_{L^\infty(\xi>A)}\le \epsilon_0(A),
\end{equation}
then
\begin{equation}\label{coer9}
|\zeta^*_1|+|\Bl1-\lambda_2| \le C(A)E\mbox{ and }\|\Jac_{\zeta, \eta}(\zeta^*, \lambda)^{-1}\|\le C(A).
\end{equation}
Taking $\xi=A-\zeta^*_1\ge A$ then $\xi=A$ in \aref{petit9}, we write
\begin{equation}\label{2points}
|\Bl1\cosh^{-\frac 2{p-1}}(A-2\zeta^*_1)-\lambda_2\cosh^{-\frac 2{p-1}}(A-\zeta^*_1)|+|\Bl1 \cosh^{-\frac 2{p-1}}(A-\zeta^*_1) - \lambda_2\cosh^{-\frac 2{p-1}}A|\le 
\frac E{\kappa_0}.
\end{equation}
Using \aref{condA2}, we see that 
\[
|\zeta^*_1|\le C(A)\mbox{ and }|\cosh^{-\frac 2{p-1}}(A-2\zeta^*_1)\cosh^{-\frac 2{p-1}}A- \cosh^{-\frac 4{p-1}}(A-\zeta^*_1)|\le C(A)E,
\]
hence
\[
 C(A)E\ge |\cosh(A-2\zeta^*_1)\cosh A- \cosh^2(A-\zeta^*_1)|=\frac 12|1-\cosh (2\zeta^*_1)|,
\]
and the first estimate of \aref{coer9} follows from \aref{2points}, 
provided that $\epsilon_0$ is small enough.\\
Now, we prove the second estimate of \aref{coer9}.
Since we know that $\|\Jac_{\zeta,\eta}(d,\nu)^{-1}\|\le \frac {C(A)}{1-d^2}$ from \aref{jac-0}, it is enough to prove that under \aref{condA2}, we have
\begin{equation}\label{jac-2}
\|\Jac_{d,\nu}(\zeta^*,\lambda)^{-1}\|\le C(A)(1-d^2)
\end{equation}
From \aref{jac1}, we see that
\[
\Jac_{d,\nu}(\zeta^*,\lambda)^{-1}=-
\left(
\begin{array}{cc}
(1+\nu)(1-d^2)&\frac{d\lambda^{p-2}(p-1)}2\left[(1+\nu)^2-d^2\right]\\
d(1-(1+\nu)^2)&\frac{(1+\nu)\lambda^{p-2}(p-1)}2\left[(1+\nu)^2-d^2\right]
\end{array}
\right)
\]
Using \aref{condA2} and the definitions \aref{defze0} and \aref{deflambda} of $\mu_i$, $\eta_i$ and $\lambda_i$, we see that 
\aref{jac-2} follows, and so does the second estimate of \aref{coer9}.
This concludes the proof of 
Lemma \ref{lemboundk*}.\Box

\section{Technical computations on the solitons}\label{apptrans}
This Appendix is devoted to the proofs of Claims \ref{clalgebre0}, \ref{clvanish}, \ref{clprep}, \ref{clproche} and (i) of Claim \ref{cl1}.

\medskip

{\it Proof of Claim \ref{clalgebre0}}: From the transformation \aref{defw**}, Claim \ref{clalgebre0} is a direct consequence of the following:
\begin{cl}[Transformation of the error terms]\label{clalgebre}$ $\\
(i) It holds that
\begin{eqnarray}
\tr_1(\ty, \ts)&=&(1-(1-\tB) \tx e^{\ts})^{-\frac 2{p-1}} \chr_1(\cy,\cs),\label{faou1}\\
\partial_\ty\tr_1(\ty, \ts)&=&(1-(1-\tB) \tx e^{\ts})^{-\frac {p+1}{p-1}} \pY \chr_1(\cy,\cs),\label{faou2}\\
\tr_2(\ty, \ts)&=& (1-(1-\tB)\tx e^{\ts})^{-\frac {p+1}{p-1}}\left[\chr_2(\cy,\cs)+\frac {2(1-\tB)\tx e^{\ts}}{p-1}\chr_1(\cy,\cs)\right.\nonumber\\
&&\left.+\frac{\tx e^{\ts}(1+\ty (1-\tB))}{1-(1-\tB)\tx e^{\ts}}\pY \chr_1(\cy,\cs)\right]\label{faou3}
\end{eqnarray}
where $(\cy,\cs)$ and $(\ty, \ts)$ are linked by \aref{defw**}.\\
(ii) For all $\k\in[\hk,k]$, 
$L_\k>0$ large enough and $|\tx|$ small enough with $\tx <0$, for all 
$\ts\in[s_{\k+1},s_\k]$
and $\ty \in [\ty_1(\ts), 1]$, we have
\[
\begin{array}{rclrcl}
-e^{L_\k}&\le&(1-\tB)\tx e^{\ts}\le 0,&\d\frac 1{1+e^{L_\k}}&\le& (1-(1-\tB) \tx e^{\ts})^{-1}\le 1,\\
\\
1-\ty&\le & (1+e^{L_\k})(1-\cy),&1+\ty&\le & 2(1+e^{L_\k})(1+\cy),\\
\\
\frac{|\tx|e^{\ts}(1+\ty (1-\tB))}{1-(1-\tB)\tx e^{\ts}} &\le & (1+2e^{L_\k})(1-{\cy}^2).
\end{array}
\] 
\end{cl}
Thus, the proof reduces to the proof of Claim \ref{clalgebre}.

\medskip

{\it Proof of Claim \ref{clalgebre}}:\\ 
(i)- {\it Proof of \aref{faou1}}: It follows directly from \aref{rach1} and \aref{rach2} thanks to the linear character of $\tT$ and the definitions \aref{defg**}, \aref{defh**} and \aref{defs2} of $\chr_1$, $\tr_1$, $\tmu_i$ and $\cd_i$.\\
- {\it Proof of \aref{faou2}}: It follows from the differentiation of \aref{faou1}.\\
{\it Proof of \aref{faou3}}: Differentiating \aref{rach1} and \aref{rach2} with respect to $\ts$, we get from the definition \aref{defk*} of $\kappa^*$
\[
\begin{array}{rcl}
\partial_\ts w_\tx(\ty, \ts)&=& (1-(1-\tB)\tx e^{\ts})^{-\frac {p+1}{p-1}}\left[\pS W(\cy,\cs)+\frac {2(1-\tB)\tx e^{\ts}}{p-1}W(\cy,\cs)\right.\nonumber\\
&&\left.+\frac{\tx e^{\ts}(1+\ty (1-\tB))}{1-(1-\tB)\tx e^{\ts}}\pY W(\cy,\cs)\right]\\
\kappa^*_2(d,\nu,\ty)&=& (1-(1-\tB)\tx e^{\ts})^{-\frac {p+1}{p-1}}\left[\frac {2(1-\tB)\tx e^{\ts}}{p-1}\kappa(d,\cy)+\frac{\tx e^{\ts}(1+\ty (1-\tB))}{1-(1-\tB)\tx e^{\ts}}\py \kappa(d,\cy)\right]
\end{array}
\]
where $d$ is arbitrary in $(-1,1)$ and $\nu =[\tB-(1-d)]\tx e^{\ts}$.
Using the definitions \aref{defg**} and \aref{defh**} of $\chr$ and $\tr$, we get the conclusion of (i).

\medskip

\noindent (ii) Consider 
$\k\in[\hk,k]$, $L_\k>0$ and $\ts\in[s_{\k+1}, s_\k]$.  
Since we have from the definition \aref{defsm}-\aref{deftsa} of $s_\k$, 
\[
s_\k \le l+L_\k,\mbox{ hence }|\tx| e^{s_\k}\le e^{L_\k},
\]
we write from \aref{chap**} for $|\tx|$ small enough,
\[
0<-(1-\tB)\tx e^{\ts}\le |\tx| e^{\ts}\le |\tx| e^{s_\k}\le e^{L_\k},
\]
which yields the first and the second inequality.\\
Using \aref{defw**}, \aref{chap**} and the second inequality, we write for $|\tx|$ small enough,
\[
1-\cy = \frac{1-\ty-(2-\tB)\tx e^\ts}{1-(1-\tB) \tx e^{\ts}}\ge \frac{1-\ty}{1+e^{L_\k}},
\]
and the third inequality follows.\\
Using again \aref{defw**}, we write
\begin{equation}\label{arab}
1+\cy = \frac{1+\ty+\tB\tx e^\ts}{1-(1-\tB) \tx e^{\ts}}.
\end{equation}
Since we have by definition \aref{defy1} of $\ty_1(\ts)$, 
\[
\mbox{if } \ty\ge \ty_1(\ts),\mbox{ then }|\tB \tx e^{\ts}|=\frac 12 (1+\ty_1(\ts))\le \frac 12 (1+\ty),
\] 
it follows from \aref{arab} and the second inequality that
\[
1+\cy \ge \frac{1+\ty}{2(1+e^{L_\k})},
\]
and the forth inequality follows.\\
Since $\ty \in (\ty_1(\ts), 1)$, it follows from \aref{defw**}, \aref{chap**}, the first and the second inequality that $\cy\in (\cy_1,\cy_2)$ with
\[
\cy_1= \frac{-1+\tx e^{\ts}(1-2\tB)}{1-(1-\tB) \tx e^{\ts}}\mbox{ and }
\cy_2 = \frac{1+\tx e^{\ts}}{1-(1-\tB) \tx e^{\ts}},
\]
hence
\begin{eqnarray}
1-\cy&>& 1-\cy_2 = \frac{(2-\tB)|\tx| e^{\ts}}{1-(1-\tB) \tx e^{\ts}}\ge \frac{|\tx|e^{\ts}}{1+e^{L_\k}},\label{dan1}\\
\frac{1+\cy}{1+(1-\tB)\cy}&>& \frac{1+\cy_1}{1+(1-\tB)\cy_1}=\frac{|\tx| e^{\ts}}{1-2(1-\tB) \tx e^{\ts}}\ge \frac{|\tx|e^{\ts}}{1+2e^{L_\k}}\label{dan2}
\end{eqnarray}
Since we have from \aref{defw**},
\begin{equation*}
\frac{|\tx|e^{\ts}(1+\ty (1-\tB))}{1-(1-\tB)\tx e^{\ts}} = |\tx|e^{\ts} (1+\cy (1-\tB)),
\end{equation*}
using \aref{dan1} and \aref{dan2} we get the fifth inequality.
%
%
%
This concludes the proof of  Claim \ref{clalgebre} and Claim \ref{clalgebre0} too.\Box

\bigskip

Now, we prove Claim \ref{clvanish}.

\medskip

{\it Proof of Claim \ref{clvanish}}: Note first from \aref{equid**} that when $1\le i<\frac{k+1}2 <j \le k$ and $S\ge S_0$ is large enough, we have
\begin{equation}\label{expDi}
\frac {S^{-\gamma_i}}{C_0}\le 1-D_i(S) \le C_0S^{-\gamma_i},\;\;|D_{\frac{k+1}2}(S)|\le 1-\frac 1{C_0}\mbox{ and }\frac {S^{\gamma_j}}{C_0}\le 1+D_j(S) \le C_0S^{\gamma_j}
\end{equation}
(of course, the middle estimate holds only when $\frac{k+1}2\in\N$).\\
Now, consider 
$L_{k+1}>0$,...,$L_\hk>0$, 
$x<0$ and $\k\in[\hat k,k]$.
Using the definition \aref{defSm} of $S_\k$,
we see from the definitions \aref{defsk+1}, \aref{defsm} and \aref{deftsa} that when $|x|\to 0$, we have
\begin{eqnarray}
\cs_{k+1}&=&L_{k+1}+O(x),\nonumber\\
\cs_\k&=&l+\gamma_\k\log l +L_\k +O(l^{\gamma_\k})\mbox{ if }\hat k+1\le \k \le k,\label{sm}\\
\cs_{\hat k}\le S(l+L_\hk)&=&l-\log(1+e^{-L_\hk})+O(l^{-\gamma_1}).\nonumber
\end{eqnarray}
Take first $i\ge \hat k+1$ and $s\in[s_i,l+L_\hk]$. Using \aref{sm}, \aref{expDi}, the bound \aref{chap**} on $\tB$ and the definition \aref{defs2} of $\tmu_i$, we see that as $x\to 0$,
\begin{equation}\label{phi}
\frac l2\le S(s)\le 2l,\;\;\frac{l^{\gamma_i}}C \le 1+\cd_i(s) \le C l^{\gamma_i}\mbox{ and }\tmu_i(s) \ge \frac{e^{L_i} l^{\gamma_i}}C>0.
\end{equation}
Recalling from the definition \aref{deflambda} of $\lambda$ that
\begin{equation}\label{a1}
\lambda(d,\nu)^{1-p} = \left(1+\frac{\nu}{1-d}\right)\left(1+\frac{\nu}{1+d}\right),
\end{equation}
we see from \aref{phi} that $\lambda(\cd_i(s_i), \tmu_i(s_i))^{1-p}\ge \frac{e^{L_i}}C$ and the conclusion follows from (i) of Lemma \ref{lemboundk*}.\\
Now, if $i=2,...,\hat k$ and $s=l+L_\hk$, then
using \aref{sm}, \aref{chap**} and \aref{expDi}, we see that
\[
2\tB \le \frac C{l^{\gamma_1}} \le 1- \cd_i(l+L_\hk),\mbox{ hence }\tmu_i(l+L_\hk) \ge \frac {(1-\cd_i(l+L_\hk))e^{L_\hk}}2>0
\]
from the definition \aref{defs2} of $\tmu_i$. Using \aref{a1}, we see that $\lambda(\cd_i(l+L), \tmu_i(l+L))^{1-p}\ge \frac{e^{L_\hk}}C$, and the conclusion follows again from (i) of Lemma \ref{lemboundk*}. This concludes the proof of Claim \ref{clvanish}.\Box

\bigskip

Now, we prove Claim \ref{clprep}.

\medskip 

{\it Proof of Claim \ref{clprep}}:\label{ppclprep} Consider 
$\k\in[\hk,k]$,
$L_\k>0$, $\lzero>0$, $x<0$, 
and $s\in[s_{\k+1}, s_\k]$. Using \aref{sm}, we see that for 
$|x|$ small enough, we have
\begin{equation}\label{fayssal}
\begin{array}{rcccll}
\lzero&\le& \ts& \le& l+\min(\gamma_\k,0) \log l +L_\k,&\\
\lzero-1&\le& \cs(s)& \le& l+1&\mbox{and for }\k\le k-1,\;\;|\cs -l|\le 1.
\end{array}
\end{equation}
We first prove \aref{first} then \aref{second}.

\medskip
{\it Proof of \aref{first}}: Consider $i=1,...,\k$. Using the definition \aref{defs2} of $\cd_i$ and $\tmu_i$, \aref{expDi}, \aref{fayssal} and the bound \aref{chap**} on $\tB$, we write for 
$\lzero$
large enough and $|x|$ small enough:\\
- If $\frac{k+1}2<i\le \k$, then $\cd_i(s)<0$ and $\gamma_\k<0$ by definition \aref{equid**}. Therefore,
\begin{equation}\label{nadia}
0\le \frac{\tmu_i(s)}{1-|\cd_i(s)|}= \frac{|x|e^s}{1+\cd_i(s)}[1-\cd_i(s)-\tB]\le \frac{2|x|e^\ts}{1+D_\k(\cs)}\le 2C_0 \cs^{-\gamma_\k}|x|e^\ts\le 3C_0 e^{L_\k}
\end{equation}
and \aref{first} follows.\\
- If $i=\frac{k+1}2$,
then
\begin{equation}\label{nadia2}
0\le \frac{\tmu_i(s)}{1-|\cd_i(s)|}\le C|x|e^\ts \le Ce^{L_\k}
\end{equation}
and \aref{first} follows in this case.\\
- If $i<\frac{k+1}2$, then 
\begin{equation}\label{parrot}
\cd_i(s)>0\mbox{ and }\frac{\tmu_i(\ts)}{1-|\cd_i(\ts)|}=|x|e^\ts\left[1-\frac \tB{1-D_i(\cs)}\right]
\le |x|e^\ts \le e^{L_\k}
\end{equation}
and 
\begin{equation}\label{adrien}
\frac \tB{1-D_i(\cs)}\le C_0^2 l^{\gamma_i-\gamma_1}.
\end{equation}
If $2\le i<\frac{k+1}2$, then it follows by definition \aref{equid**} of $\gamma_i$ that $\frac \tB{1-D_i(\cs)}\le \frac 12$
for $|x|$ small enough. Using \aref{parrot}, we see that
\begin{equation}\label{ghoul-1}
\frac{\tmu_i(\ts)}{1-|\cd_i(\ts)|} \ge 0\mbox{ when }2\le i <\frac{k+1}2
\end{equation}
and \aref{first} follows from \aref{parrot} and \aref{ghoul-1}.\\
If $i=1$
and $\k=\hk$,
recalling 
that $\ts\le \tsa$ for $|x|$ small enough,
 we see by definition \aref{deftsa} of $\tsa$ and \aref{parrot} that
\begin{equation}\label{ghoul}
\frac{\tmu_1(\ts)}{1-|\cd_1(\ts)|}=\frac{\tmu_1(\ts)}{1-\cd_1(\ts)}\ge -1+e^{-(p-1)L_{\hk}}
\end{equation}
and \aref{first} follows from \aref{parrot} and \aref{ghoul}.\\
If $i=1$ and $\k\in[\hk+1,k]$, then we see from \aref{parrot} and the definitions \aref{defsm} and \aref{equid**} of $s_\k$ and $\gamma_\k$ that
\begin{equation}\label{luba}
\frac{\tmu_1(\ts)}{1-|\cd_1(\ts)|} \le |x|e^s \le |x| e^{s_\k} = l^{\gamma_\k}e^{L_\k} \le l^{-\frac{(p-1)}2}e^{L_\k}\to 0\mbox{ as }x\to 0.
\end{equation}
Therefore, \aref{first} follows in this case too.\\
This concludes the proof of \aref{first}.

\medskip 

{\it Proof of \aref{second}}: 
 Using the definition \aref{defcz} of $\czeta_i^*(s)$, we write for $i=1,...,\k$,
\begin{eqnarray}
&&\czeta_i^*(s)=-\argth\left(\frac{\cd_i(s)}{1+\tmu_i(s)}\right)=\frac 12 \log\left[\frac{1-\cd_i(s)+\tmu_i(s)}{1+\cd_i(s)+\tmu_i(s)}\right]\label{ghoul1}\\
&=& \frac 12 \log\left(\frac{1-\cd_i(s)}{1+\cd_i(s)}\right)+\frac 12 \log\left(1+\frac{\tmu_i(s)}{1-\cd_i(s)}\right)-\frac 12 \log \left(1+\frac{\tmu_i(s)}{1+\cd_i(s)}\right).\label{ghoul2}
\end{eqnarray}
Since the two last terms are bounded by $CL_\k$ from \aref{first}, and $\cd_i(s)=D_i(\cs(s))$ by \aref{defs2}, we write from \aref{equid**}
\[
\left|\czeta_i^*(s)+\frac{\gamma_i}2 \log \cs(s)\right|\le C^*L_\k\mbox{ where }\gamma_i = (p-1) \left(\frac{k+1}2-i\right)
\]
for some $C^*>0$.\\
If $\frac{(p-1)}2\log \cs(s) \ge 3 C^*L_\k$, then for all $i=1,...,\k-1$, $\czeta^*_{i+1}(s) - \czeta^*_i(s)\ge \frac{(p-1)}6 \log S(s) \ge \frac{(p-1)}6 \log(\lzero-1)$ by \aref{fayssal} and \aref{second} follows by taking $L_{\k+1}$ large enough.\\
If $\frac{(p-1)}2\log \cs(s) < 3 C^*L_\k$, then for $|x|$ small enough, we have
\begin{equation}\label{fouillade}
\lzero-1 \le \cs(s) \le e^{\frac{6C^*}{p-1}L_\k}\mbox{ and }\lzero \le \ts \le C(L_\k)
\end{equation}
for some $C(L_\k)>0$, where we used \aref{fayssal} and the expression \aref{defs2*} of $\cs(s)$. Using \aref{nadia}, \aref{nadia2} and \aref{parrot}, we see that
\[
\frac{\tmu_i(\ts)}{1-|\cd_i(\ts)|}\le C(L_\k) |x| \le \frac 1{10}
\]
for $|x|$ small enough. Recalling from \aref{nadia}, \aref{nadia2} and \aref{ghoul-1} that
\begin{equation}\label{i2}
\tmu_i(s) \ge 0\mbox{ when }i\ge 2,
\end{equation}
we derive from \aref{ghoul2}, \aref{defs2} and \aref{equid**}:
\begin{equation}\label{mouna}
\mbox{for }2\le i\le \k,\;\;\left|\czeta_i^*(s)+\frac{\gamma_i}2 \log \cs(s)\right|\le C\mbox{ where }\gamma_i = (p-1) \left(\frac{k+1}2-i\right).
\end{equation}
If $i=1$ and $\tmu_1(s)\ge 0$, then \aref{mouna} holds for $i=1$ by the same argument as for $i\ge 2$.\\
If $i=1$ and $\tmu_1(s)< 0$, then, recalling from \aref{parrot} that $\cd_1(s) >0$, we see by 
the expression
\aref{ghoul1} of $\czeta^*_1(s)$, \aref{defs2} and \aref{equid**} that
\begin{equation}\label{thomas}
\czeta^*_1(s) \le -\argth (\cd_1(s)) \le -\frac{\gamma_1}2 \log \cs(s) +C.
\end{equation}
Using \aref{mouna}, the line after \aref{mouna}, \aref{thomas} and \aref{fouillade}, we see that for 
$|x|$ small enough, for all $i=1,...,\k-1$, 
\[
\czeta_{i+1}^*(s) - \czeta_i^*(s)\ge \frac{(p-1)}2 \log S(s) -C \ge \frac{(p-1)}2 \log (\lzero-1)-C,
\]  
and estimate \aref{second} follows when $\frac{(p-1)}2 \log \cs(s) < 3 C^*L_\k$, provided that $\lzero$ is large enough. This concludes the proof of \aref{second} and Claim \ref{clprep} as well. \Box

\bigskip

Now, we prove Claim \ref{clproche}.

\medskip

{\it Proof of Claim \ref{clproche}}: From the upper bound in (ii) of Lemma \ref{lemboundk*}, Claim \ref{clprep} and \aref{apriori}, we clearly see that (ii) follows from (i) and estimate \aref{small}. Thus, we only prove (i). In fact, we will bound only the two first terms in (i) since the third is bounded by their sum thanks to \aref{lip} and \aref{jaco}. Take $s\in[s_{\k+1}, \bar s]$. Using \aref{ordre}, Lemma \ref{lemleft}, Claims \ref{clvanish} and \ref{clprep}, and \aref{gap}, we see that
\begin{equation}\label{partiel}
\left\|\sum_{i=1}^\k (-1)^i\left(\kappa^*(d_i(s), \nu_i(s))- \kappa^*(\cd_i(s), \tmu_i(s))\right)\right\|_{\H(y>y_1(s))}\le 2\epsilon
\end{equation}
where $y_1(s)$ is given in \aref{defy1}
and for all $i=1,...,m-1$,
\begin{equation}\label{gap2}
\czeta_{i+1}^*(s) -\czeta_i^*(s)\ge 
\frac{(p-1)}7 \log \lzero
\mbox{ and }\zeta_{i+1}^*(s) -\zeta_i^*(s)\ge 
\frac{(p-1)}9 \log \lzero
\end{equation}
if 
$\lzero$ is large enough and $|x|$ is small enough. Note first that for $L_\k$ large enough and $|x|$ small enough, we have
\begin{equation}\label{y1left}
\forall i=1,...,\k,\;\;\argth y_1(s) \le -\argth\left(\frac{\cd_i(s)}{1+\tmu_i(s)}\right)+pL_\k
\end{equation}
(see below for a proof of this fact).

\medskip

Using \aref{y1left}, we see that \aref{partiel} yields by truncation
\begin{equation}\label{kcutoff}
\|\kappa^*(d_\k(s), \nu_\k(s))-\kappa^*(\cd_\k(s), \tmu_\k(s))\|_{\H(y>\tanh(\argth(\max(-\frac{\cd_\k(s)}{1+\tmu_\k(s)},-\frac{d_\k(s)}{1+\nu_\k(s)}))+pL_\k))} \le 2\epsilon.
\end{equation}
From Claim \ref{clprep} and \aref{apriori}, we see that the lower bound in (ii) of Lemma \ref{lemboundk*} applies provided that $\epsilon\le \epsilon_0(L_\k)$ for some $\epsilon_0(L_\k)$, and we get the conclusion for $i=\k$.\\
Using now the upper bound in (ii) of Lemma \ref{lemboundk*}, we see that \aref{kcutoff} actually holds in $\H$, in the sense that
\[
\|\kappa^*(d_\k(s), \nu_\k(s))-\kappa^*(\cd_\k(s), \tmu_\k(s))\|_{\H} \le C(L_\k)\epsilon.
\]
Using this, we see that \aref{partiel} now holds with a sum running from $i=1$ up to $i=\k-1$. Iterating the above argument by induction for $i$ decreasing from $\k$ to $1$, one can easily conclude the proof of of (i) in Claim \ref{clproche}. It remains to prove \aref{y1left} in order to conclude.

\medskip

{\it Proof of \aref{y1left}}: From \aref{gap2} and the definition \aref{defcz} of $\czeta_i^*(s)$, it is enough to prove this for $i=1$. 
Using \aref{fayssal} and \aref{expDi}, we see that
\begin{equation}\label{pool}
\cs\le l+ 1\mbox{ and }\frac{l^{-\gamma_1}}C \le 1-\cd_1(s)
\end{equation}
for 
$\lzero$ large enough and $|x|$ small enough.\\
- {\bf Case $\hk+1\le \k\le k$}: Using the definition \aref{defs2} of $\tmu_i(s)$ and \aref{pool}, we see that
\begin{equation}\label{hamm}
-\frac{\cd_1(s)}{1+\tmu_1(s)}=-\frac{\cd_1(s)}{1+[b-(1-\cd_1(s))xe^s]}\ge -\frac{1-\frac{l^{-\gamma_1}}C}{1+[b-\frac{l^{-\gamma_1}}C]xe^s}.
\end{equation}
Since $s\le s_\k$, we see from the bound \aref{chap**} on $b$ and the definition \aref{defsm} of $s_\k$ that $\left|[b-\frac{l^{-\gamma_1}}C]xe^s\right|\le Cl^{-\gamma_1}l^{-\frac{(p-1)}2}e^{L_\k}$. Using \aref{hamm}, we write for $|x|$ small enough,
\begin{equation}\label{bb27}
-\frac{\cd_1(s)}{1+\tmu_1(s)}\ge -1 +\frac{l^{-\gamma_1}}C,\mbox{ hence }-\argth \left(\frac{\cd_1(s)}{1+\tmu_1(s)}\right)\ge - \frac{\gamma_1}2 \log l-C.
\end{equation}
Using the definition \aref{defy1} and \aref{defsm} of $y_1(s)$ and $s_\k$ together with \aref{equid**}, we write
\[
y_1(s) = -1 +2b|x|e^s \le -1 +C l^{-\gamma_1}l^{-\frac{(p-1)}2}e^{L_\k},
\]
hence
\begin{equation}\label{bb28}
\argth y_1(s) \le -\left(\frac{\gamma_1}2+\frac{p-1}4\right)\log l + \frac{L_\k}2+C.
\end{equation}
Thus, \aref{y1left} follows from \aref{bb27} and \aref{bb28} for $|x|$ small enough, in the case where $\hk+1\le \k\le k$.\\
- {\bf Case $\k= k$}: Since $s\le \tsa$, we see from the definition \aref{deftsa} of $\tsa$ that
\begin{equation*}
-\frac{\cd_1(s)}{1+\tmu_1(s)}\ge -\frac{\cd_1(s)}{\cd_1(s)+e^{-(p-1)L_\k}(1-\cd_1(s))}\ge -1 +\frac{e^{-(p-1)L_\k}C} l^{-\gamma_1},
\end{equation*}
hence
\begin{equation}\label{down}
-\argth\left(\frac{\cd_1(s)}{1+\tmu_1(s)}\right) \ge -\frac{\gamma_1}2 \log l - \frac{(p-1)}2 L_\hk-C.
\end{equation}
Using \aref{chap**}, the definitions \aref{defsm}, \aref{equid**} and \aref{defy1} of $s_{\k}$, $\gamma_i$ and $y_1(s)$, we see that 
\begin{equation}\label{deux}
y_1(s) = -1 - 2bxe^s \le -1 +Cl^{-\gamma_1} e^{L_\hk}\mbox{ hence }\argth y_1(s) \le -\frac {\gamma_1}2 \log l +\frac {L_\hk}2+C.
\end{equation}
Thus, from \aref{down} and \aref{deux}, we see that \aref{y1left} holds for $L_\k$ large enough. This concludes the proof of \aref{y1left} and (i) of Claim \ref{clproche} too. Since (ii) of Claim \ref{clproche} follows from (i) as we said in the beginning of the proof, this concludes the proof of Claim \ref{clproche}.\Box

\bigskip

Now, we prove (i) of Claim \ref{cl1}.

\medskip

{\it Proof of (i) of Claim \ref{cl1}}: Consider 
$M_0>0$, $x<0$ and $s\in[s_{\k+1}, s_\k-M_0]$. From \aref{apriori}, 
\aref{id9A} (applied with $A=CL_\k$)
and Claim \ref{clproche}, we see that
\[
\d\sum_{i=1}^\k \frac{|\nu_i(s)|}{1-|d_i^*(s)|}\le 
C(L_\k)\d\sum_{i=1}^\k \frac{|\nu_i(s)|}{1-|d_i(s)|}\le 
C(L_\k) \d\sum_{i=1}^\k \frac{|\tmu_i(s)|}{1-\cd_i(s)^2}+C(L_\k) \epsilon. 
\]
Taking $\epsilon$ small enough so that $C(L_\k) \epsilon \le \frac{\eta_0}2$, it is enough to show that for all $i=1,...,\k$,
\begin{equation}\label{rere}
\frac{|\tmu_i(s)|}{1-\cd_i(s)^2}\le \frac{\eta_0}{2kC(L_\k)}.
\end{equation}
 Note first from the definitions \aref{defsk+1}-\aref{deftsa} of $s_{\k+1}$ and $s_\k$ that for $|x|$ small enough,
\begin{equation}\label{bornes2}
\lzero \le s \le l+L_\k-M_0.
\end{equation}
Using the definition \aref{defs2*} of $S=S(s)$, we see that if $M_0 \ge L_\k$ and $|x|$ is small enough, then 
\begin{equation}\label{bornecs2}
\lzero-1 \le \cs \le l+2.
\end{equation}
From the definition \aref{defs2} of $\tmu_i$ and $\cd_i$, we write
\begin{equation}\label{habib}
\frac{\tmu_i(s)}{1-\cd_i(s)^2} =\frac{|x|e^s}{1+D_i(\cs)}\left(1-\frac \tB{1-D_i(\cs)}\right).
\end{equation}
Using \aref{bornecs2}, \aref{expDi}, the bound \aref{chap**} on $\tB$ and the definitions \aref{defsk+1}-\aref{deftsa} of $s_\k$, we write
\[
\frac{|\tmu_i(s)|}{1-\cd_i(s)^2} \le \frac{C|x|e^s}{1+D_i(\cs)} \le \frac{C|x|e^s}{1+D_\k(\cs)}\le C_0e^{L_\k-M_0}
\]
and \aref{rere} follows if $M_0=M_0(L_\k)$ is large enough.\\
This concludes the proof of (i) of Claim \ref{cl1}.\Box

\section{Analysis of the linearization of equation \aref{eqw1} around a sum of $\kappa^*(d_i,\nu_i)$}\label{appproj9} 
This Appendix is devoted to the study of the dynamics of equation \aref{eqw1}  near a decoupled sum of $\pm\kappa^*(d_i, \nu_i)$ and to the proof of Claim \ref{lemlyap}. This is in fact a generalization of the case where all $\nu_i=0$ already treated in Appendix C in \cite{MZajm10}, and of the case of one soliton with $\nu_i=0$ treated in Section 5 page 99 in \cite{MZjfa07}.

\medskip

We have 3 steps in the following:\\
- In Section \ref{sub1}, we give the equation satisfied by the error term $q$ defined in \aref{defq} when we make the linearization. We also give a decomposition of the error term, well-adapted to the spectrum of the linearization.\\
- In Section \ref{sub2}, we project the equation on the various modes and prove Claim \ref{lemlyap}.\\
- Finally, we give in Section \ref{sub3} a table for some integrals appearing in the proof.

\subsection{An equation satisfied by $q$}\label{sub1}
Using the fact that when $\nu = \mu e^s$, the function $(y,s)\mapsto \kappa^*(d,\mu e^s,y)$ is a particular solution of equation \aref{eqw1}, we see from equation \aref{eqw1} and the change of variables \aref{defq} that $q$ defined in \aref{defq} satisfies the following equation for all $\ts \in [s_{\k+1},\bar s]$:
\begin{eqnarray}
\d\frac \partial {\partial s}
\left(
\begin{array}{l}
q_1\\
q_2
\end{array}
\right)
&=&\ll
\left(
\begin{array}{l}
q_1\\
q_2
\end{array}
\right)
-\sum_{j=1}^\k(-1)^j\left[(\nu_j'(s)-\nu_j(s))\pnu\kappa^*+d_j'(s)\partial_d \kappa^*\right](d_j(s),\nu_j(s),y)\nonumber\\
&+&\vc{0}{R}
+\left(
\begin{array}{l}
0\\
f(q_1)\label{eqq*}
\end{array}
\right)
\end{eqnarray}
where
\begin{eqnarray}
\ll\vc{q_1}{q_2}&=&\vc{q_2}{\L q_1+\psi q_1-\frac{p+3}{p-1}q_2-2y\py q_2},\label{defL9}\\
\psi(y,s)&=&p|K^*_1(y,s)|^{p-1} -\frac{2(p+1)}{(p-1)^2},\;\;
K^*_1(y,s) = \sum_{j=1}^\k  (-1)^i\kappa^*_1(d_j(s),\nu_j(s),y),\label{defpsi}\\
f(q_1)&=&|K^*_1+q_1|^{p-1}(K^*_1+q_1)- |K^*_1|^{p-1}K^*_1- p|K^*_1|^{p-1} q_1,\nonumber\\
R&=& |K^*_1|^{p-1}K^*_1- \sum_{j=1}^\k (-1)^j\kappa^*_1(d_j(s),\nu_j(s),y)^p.\nonumber
\end{eqnarray}
Let us remark that equation \aref{eqq*} can be localized near the center 
$d^*_i(s)$
of $\kappa^*(d_i(s), \mu_i(s))$ (which is the same as the center of $\kappa\left(\frac {d_i(s)}{1+\nu_i(s)}\right)$ by (i) of Lemma \ref{lemboundk*}) for each $i=1,...,\k$, which allows us to view it locally as a perturbation of the case of $\kappa(d,y)$ already treated in \cite{MZjfa07}. For this, given $i=1,...,\k$, we need to expand the linear operator of equation \aref{eqq*} as 
\begin{equation}\label{expl}
\ll(q) = \LL_{d_i^*(s)}(q) + (0, \bar V_i(y,s) q_1)+ (0, V_i^*(y,s) q_1)
\end{equation}
with
\begin{eqnarray}
\LL_{d}\vc{q_1}{q_2}&=&\vc{q_2}{\L q_1+\psi^*(d)q_1-\frac{p+3}{p-1}q_2-2yq_2'},\label{defld}\\
\psi^*(d,y)&=&p\kappa(d,y)^{p-1}-\frac{2(p+1)}{(p-1)^2},\label{defpsi*}\\
\bar V_i(y,s)&=&p \kappa^*(d_i(s), \nu_i(s),y)^{p-1}-p\kappa(d_i^*(s),y)^{p-1},\label{defbvi}\\
V_i^*(y,s)&=& p|K^*_1(y,s)|^{p-1} - p \kappa^*(d_i(s), \nu_i(s),y)^{p-1}.\label{defvi-09}
\end{eqnarray}
Since we will work in the regime where the solitons' sum is decoupled (see the estimate \aref{gap}),
the properties of $\LL_{d_i(s)}$ will be essential in our analysis.\\
From section 4 in \cite{MZjfa07}, we know that for any $d\in (-1,1)$, the operator $\LL_d$ has $1$ and $0$ as eigenvalues, the rest of the eigenvalues are negative. More precisely, 
we have for $l=0$ and $1$,
\begin{equation}\label{79bis}
\LL_d (F_l(d)) = l F_l(d)\mbox{ and }
\|F_1(d)\|_{\H}+\|F_0(d)\|_{\H}\le C
\end{equation}
where $F_l(d)$ is defined in \aref{deffld}.
The projection on $F_l(d)$ is defined in \aref{defpdi} by $\pi_l^d(r) =\phi\left(W_l(d), r\right)$ where $W_l(d)$ is defined in \aref{defWl2-0} and \aref{eqWl1-0}. Of course,
\begin{equation}\label{79bis*}
\LL_d^* W_l(d) = l W_l(d)
\end{equation}
where $\LL_d^*$ is the conjugate of $\LL_d$ with respect to the inner product $\phi$ defined in \aref{defPhi}. We also recall from \aref{normw} the following orthogonality condition
\begin{equation}\label{orth}
\pi_l^d(F_{l'}(d)) = \phi(W_l(d), F_{l'}(d)) = \delta_{l,l'}.
\end{equation}
In the following, we give a decomposition of the solution which is well adapted to the proof:
\begin{lem}[Decomposition of $q$]\label{lemdecomp9}
If we introduce for all  $r$ and $\r$ in $\H$ the operator
 $\pi_-(r)\equiv r_-(y,s)$ defined by
\begin{equation}\label{decomp9}
r(y,s)= \sum_{i=1}^{k}\sum_{l=0}^1 \pi^{d_i^*(s)}_l(r) F(d_i^*(s),y) + \pi_-(r),
\end{equation}
and 
\begin{equation}\label{defphi9}
\varphi(r,\r) = \iint \left(r_1'\r_1' (1-y^2)-\psi r_1\r_1+r_2\r_2\right)\rho dy,
\end{equation}
then the following holds for all $\lzero \ge L_{\k+1,5}(L_\k)$ and $|x|\le a_5(L_\k,L_{\k+1})$ for some $L_{\k+1,5}$ and $a_5>0$:\\
(i) For all $s\in[s_{\k+1},\bar s]$ 
and for all $r$ and $\r$ in $\H$, we have
\begin{equation*}
\left| \varphi(r,\r) \right| \le C(L_\k)\|r\|_{\H} \|\r\|_{\H}\mbox{ and }q(y,s)=q_-(y,s).
\end{equation*}
(ii) There exists $\eta_0>0$ such that for all $s\in[s_{\k+1},\bar s]$, if \aref{cond} holds,
then 
\begin{equation}\label{defa-9}
\frac 1{C_0}\|q(s)\|_{\H}^2\le \m(s)\le C(L_\k) \|q(s)\|_{\H}^2\mbox{ where }\m(s)=\varphi(q(s),q(s)).
\end{equation}
\end{lem}
{\bf Remark}: The operator $\pi_-$ depends on the time variable $s$. The constants $C_0>0$ and $\eta_0>0$ are universal.\\
{\it Proof}: Using \aref{apriori}, we easily see that 
\begin{equation}\label{defli}
\frac {e^{-\frac{2L_\k}{p-1}}}C \le \lambda_i(s) \le Ce^{2L_\k}
\mbox{ and }-1+\frac{e^{-(p-1)L_\k}}2 \le \nu_i \le 1+C_0e^{L_\k},
\end{equation}
where $\lambda_i(s) = \lambda(d_i(s), \nu_i(s))$ is defined in \aref{deflambda}, hence
\begin{equation}\label{cle}
\forall y\in(-1,1),\;\;
\frac {e^{-\frac{2L_\k}{p-1}}}C \kappa(d_i^*(s),y)\le \kappa^*_1(d_i(s), \nu_i(s),y)\le Ce^{2L_\k}\kappa(d_i^*(s),y)
\end{equation}
from (i) of Lemma \ref{lemboundk*}. 

\medskip

(i) Using \aref{cle}, the proof of the continuity of $\varphi$ reduces to the case where all the $\nu_i=0$ already treated in (i) of Lemma 3.10 of \cite{MZajm10}. As for the identity $q=q_-$, it is obvious by definition \aref{decomp9} of $q_-$ and the orthogonality relation \aref{mode0}. 

\medskip

(ii) The right-hand inequality follows from (i) of this lemma.
It remains then to prove the left-hand side of to conclude.\\
Clearly, from \aref{cle}, \aref{gap} and \aref{fayssal}, the interaction between the solitons in the definition \aref{defphi9} of $\varphi$ can be controlled as in Lemma 3.10 of \cite{MZajm10}, provided that $\lzero>0$ is large enough. Therefore, we only give the proof in the case of one soliton and refer to Lemma 3.10 of \cite{MZajm10} for the extension to the case of more solitons.\\
In the case of one soliton, we recall from the definition \aref{defphi9} that
\[
\varphi(r,\r) = \iint \left(r_1'\r_1' (1-y^2)-\left(p\kappa^*_1(d_i, \nu_i,y)^{p-1}-\frac{2(p+1)}{p-1}\right)r_1\r_1+r_2\r_2\right)\rho dy.
\] 
Introducing
\[
\varphi_{d^*_1}(r,\r) = \iint \left(r_1'\r_1' (1-y^2)-\left(p\kappa(d^*_1,y)^{p-1}-\frac{2(p+1)}{p-1}\right)r_1\r_1+r_2\r_2\right)\rho dy,
\]
we recall from Proposition 4.7 page 90 in \cite{MZjfa07} and the orthogonality relation \aref{mode0} that
\begin{equation}\label{coer}
\varphi_{d^*_1}(q,q)\ge \frac 1{C_1} \|q\|_{\H}^2
\end{equation}
where $C_1>0$ is independent of $d^*_1$. Using (i) of Lemma \ref{lemboundk*}, we see that
\begin{eqnarray*}
&&\varphi(q,q)- \varphi_{d^*_1}(q,q) = p(1-\lambda_1^{p-1})\iint \kappa(d^*_1,y)^{p-1} q_1^2 \rho dy.
\end{eqnarray*}
Now, if $\frac{|\nu_1(s)|}{1-|d^*_1(s)|} \le \eta_0$, then we see from the definition \aref{deflambda} of $\lambda_1$ that
\[
|\lambda_1-1|\le C \eta_0.
\]
Therefore, if $\eta_0$ is small enough, then we write from (i) of Lemma \ref{lemboundk*} and the Hardy-Sobolev inequality \aref{hs}, 
\begin{equation}\label{lee}
\varphi(q,q) - \varphi_{d^*_1}(q,q) \ge - C\eta_0 \iint q_1^2 \frac\rho{1-y^2} dy
\ge -\frac 1{2C_1}\|q\|_{\H}^2.
\end{equation}
Using \aref{coer}, we get the left-hand side. This concludes the proof of Lemma \ref{lemdecomp9}.\Box

\subsection{A dynamical study of the equation satisfied by $q$} \label{sub2}
We prove Claim \ref{lemlyap} here.
We claim that it follows from the following lemma derived from equation \aref{eqq*}:
\begin{lem}\label{lemproj*}If we fix $L_\k$ large enough and take $\epsilon \le \epsilon_6(L_\k)$, $\lzero \ge L_{\k+1,6}(L_\k,\epsilon)$ and $|x| \le a_6(L_\k,\epsilon, \lzero)$ for some $\epsilon_6>0$, $L_{\k+1,6}>0$ and $a_6>0$, then, for all $s\in[s_{\k+1},\bar s]$, 
we have (below, $C=C(L_\k)$):\\ 
(i) {\bf (Control of the modulation parameters)} 
\begin{equation*}
\forall i=1,...,\k,\;\;\frac{|d_i'(s)|+\left|\nu_i'(s)- \nu_i(s)\right|}{1-d_i^*(s)^2}\le C\left(\ep+\frac{|\nu_i(s)|}{1-d_i^*(s)^2}\right)\|q(s)\|_{\H} + C J_\k(s)
\end{equation*}
where $J_\k(s)$ is defined in \aref{defjbar}.\\
(ii) {\bf (Control of $\m(s) = \varphi(q,q)$)}
\begin{eqnarray*}
\left(\frac 12 \m+R_-\right)'&\le& - \frac 3{p-1}\iint q_{2}^2 \frac \rho{1-y^2} dy
+C\|q(s)\|_{\H}^2\left(\ep+\sum_{i=1}^\k \frac{|\nu_i(s)|}{1-d_i^*(s)^2}\right)\\
&& + C\|q(s)\|_{\H}J_\k(s) + C J_\k(s)^2
\end{eqnarray*}
where $R_-(s)= -\iint \F(q_1) \rho dy$ with 
\begin{equation}\label{defF}
\F(q_1) = \int_0^{q_1}f(\xi) d\xi = \frac{|K^*_1+q_1|^{p+1}}{p+1}-\frac{|K^*_1|^{p+1}}{p+1}-|K^*_1|^{p-1}K^* q_1 - \frac p2 |K^*_1|^{p-1}q_1^2
\end{equation}
and $R_-$ satisfies
\begin{equation*}
|R_-(s)|\le C \|q(s)\|_{\H}^{\bar p+1}\mbox{ where }\bar p = \min(p,2).
\end{equation*}
(iii) {\bf (Control of $\iint q_1q_2 \rho$)} Assume that \aref{cond}. Then, we have
\begin{equation*}
\frac d{ds}\iint q_1q_2 \rho \le 
-\frac 7{10}\m+CJ_\k(s)^2
+ C \iint q_{2}^2 \frac \rho{1-y^2}.
\end{equation*}
(iv) {\bf (A rough estimate)}
\[
\left(\frac 12 \|q(s)\|_{\H}^2+R_-\right)'\le C\|q(s)\|_{\H}^2+ C J_\k(s)^2.
\]
\end{lem}
{\bf Remark}: If we compare this statement with the analogous one in Lemma 3.11 in \cite{MZajm10}, we see that we have two new features:\\
- Here, since we allow $\nu_i$ to be different from $0$, the modulation allows us to kill for each $i$ one additional projection of $q$, namely $\pi^{d^*_i}_1(q)$ (see \aref{mode0}). As a consequence, we obtain a differential inequality satisfied by all the new modulation parameters $\nu_i$ in (i) of Lemma \ref{lemproj*} above;\\
- If $\nu_i$ is small enough so that $\kappa^*(d_i,\nu_i)$ is close to $\kappa(d^*_i)$, then, we expect to get the same bounds as in Lemma 3.11 in \cite{MZajm10}. But if $\kappa^*(d_i,\nu_i)$ is far from $\kappa(d^*_i)$, then we expect new terms to appear in the bounds. As a matter of fact, the distance between $\kappa^*(d_i,\nu_i)$ and $\kappa(d^*_i)$ can be bounded by $\d\frac{|\nu_i|}{1-{d^*_i}^2}$, which is the new term appearing in the bounds above.

\medskip 

Indeed, let us first derive Claim \ref{lemlyap} from Lemma \ref{lemproj*}, then, we will prove the latter.

\bigskip

{\it Proof of Claim \ref{lemlyap} assuming Lemma \ref{lemproj*}}: This part is similar to the proof of Proposition 3.8 in \cite{MZajm10}, when all the $\nu_i=0$.\\
Let us introduce for all $s\in [\ts_{\k+1},\bar s]$,
\begin{equation}\label{defhh}
\begin{array}{rcl}
h_1(s) &=& \frac 12\|q\|_{\H}^2-\iint \F(q_1)\rho dy,\\
h_2(s) &=& \frac 12 \varphi(q,q)-\iint \F(q_1)\rho dy+\eta \iint q_1 q_2 \rho,
\end{array}
\end{equation}
where $\varphi$ and $ \F(q_1)$ are given in \aref{defphi9} and \aref{defF}, 
and $\eta$ will be fixed small enough. 
In the following, we will show that Claim \ref{lemlyap} holds with the above defined functions $h_1(s)$ and $h_2(s)$.

(i) From the definition \aref{defhh} of $h_1$, (ii) of Lemma \ref{lemproj*} and \aref{small}, we write
\[
|h_1(s) - \frac 12\|q(s)\|_{\H}^2|\le |R_-(s)|\le C\|q(s)\|_{\H}^{\bar p+1}\le \frac{\|q(s)\|_{\H}}{100}
\]
 if $\epsilon$ is small enough, and the first identity in (i) follows.\\ 
The second identity of follows then directly from (iv) of Lemma \ref{lemproj*}.

\medskip

(ii) By definition \aref{defjbar} of $J_\k$, we write
\begin{equation*}
|J_\k'| = \frac 1{p-1} \left|\sum_{i=1}^{\k-1}({\zeta^*_{j+1}}'-{\zeta^*_j}')e^{- \frac 1{p-1}(\zeta^*_{j+1}-\zeta^*_j)}\right|\le \frac{2J_\k}{p-1}\sum_{j=1}^\k |{\zeta^*_j}'|.
\end{equation*}
Since $\zeta^*_j = -\argth {d^*_j}$ and ${d^*_j} = \frac{d_j}{1+\nu_j}$, we write from (i) of Lemma \ref{lemproj*} and \aref{defli}
\begin{equation*}
|{\zeta^*_j}'| = \frac{|{d^*_j}'|}{1-{d^*_j}^2}= \frac{|d_j'(1+\nu_j) - d_j \nu_j'|}{(1+\nu_j)^2(1-{d^*_j}^2)}\le C\left(\frac{|\nu_j|}{1-{|d^*_j|}}+\epsilon\right).
\end{equation*}
Thus, (ii) follows.

\medskip

(iii) Assume that \aref{cond} holds. From the definition \aref{defhh} of $h_2$, (ii) of Lemma \ref{lemproj*}, the coercive estimate in (ii) of Lemma \ref{lemdecomp9} and \aref{small}, we write
\[
|h_2(s) - \frac 12\m(s)|\le |R_-(s)| +\eta \|q(s)\|_{\H}^2 \le C\m(s)^{\frac{\bar p+1}2}+C \eta \m(s)\le \frac{\m(s)}{100}
\]
 if $\eta$ and $\epsilon$ are small enough, and the first identity follows from (ii) of Lemma \ref{lemdecomp9}.\\
As for the second identity, we write from Lemma \ref{lemproj*} and 
the first identity
\begin{eqnarray}
\left(\frac \m2+R_-\right)'& \le& -\frac 3{p-1} \iint q_2^2 \frac \rho{1-y^2}dy +(\frac \eta 5 +C(\epsilon+\eta_0))h_2+\frac C{\eta} J_\k^2,\nonumber\\
\frac d{ds} \iint q_1q_2 \rho &\le & -\frac 7{10}h_2+CJ_\k^2 +C_1 \iint  q_2^2 \frac \rho{1-y^2}dy\nonumber
\end{eqnarray}
for some $C_1>0$. Therefore, by definition \aref{defhh} of $h_2$, we write
\begin{eqnarray*}
h_2'&\le & -\left(\frac 3{p-1}- C_1 \eta\right)\iint q_2^2 \frac \rho{1-y^2}dy
-h_2\left[-\frac \eta 5+\frac 7{10}\eta -C(\ep+\eta_0)\right]\\
&& +\left(\frac C{\eta} +C\eta\right) J_\k^2.
\end{eqnarray*}
Fixing 
\begin{equation}\label{defeta}
\eta = \frac 3{(p-1) C_1},
\end{equation}
we see that fixing $\eta_0(L_\k)$ small enough and taking $\epsilon \in (0, \bar \epsilon]$ for some $\bar\epsilon(L_\k)>0$, we get
\begin{equation*}
h_2' \le -\frac {2\eta}5 h_2 +C_2 J_\k^2
\end{equation*}
for some $C_2>0$. This concludes the proof of Claim \ref{lemlyap} assuming Lemma \ref{lemproj*}.\Box

\bigskip

Now, we give the proof of Lemma \ref{lemproj*}.

\medskip

{\it Proof of Lemma \ref{lemproj*}}: If we set aside (iv), then we see from  \aref{defli}, \aref{cle} and \aref{gap} that our case is no more than a straightforward adaptation of the case where all $\nu_i=0$ (hence $\kappa^*_1(d_i(s), \nu_i(s),y)=\kappa(d_i^*(s),y)$ and $d_i^*(s) = d_i(s)$) treated in \cite{MZajm10}, except for 3 new terms of equation \aref{eqq*} : $(0,\bar V_i)$ given in \aref{defbvi} and the modulation terms involving $\partial_d\kappa^*$ and $\partial_\nu \kappa^*$. Therefore, we proceed in 4 parts to prove (i), (ii), (iii) and (iv) giving details only for (iv) and the three new terms. The different constants $C$ and $C^*$ below depend on $L_\k$.

\bigskip

{\bf Proof of (i): Projection of equation \aref{eqq*} on $F_1^{d_i^*}$ and $F_0^{d_i^*}$}   

We prove (i) of Lemma \ref{lemproj*} here. Fixing some $i=1,...,\k$ and $l=0$ or $1$, we project in the following each term of equation \aref{eqq*} with the projector $\pi^{d_i^*}_l$ defined in \aref{defpdi}, expanding $\ll(q)$ as in \aref{expl}.

- Note first that thanks to \aref{defli} and \aref{cle}, the analysis of Appendix C in \cite{MZajm10} holds here and gives the following estimates:
\[
\begin{array}{l}
\left|\pi^{d_i^*}_l(\partial_s q)\right|
\le C\d\frac{|{d_i^*}'|\cdot\|q\|_{\H}}{1-{d_i^*}^2}\le C\d\frac{\|q\|_{\H}}{1-{d_i^*}^2}(|d_i'|+|\nu_i'-\nu_i|+|\nu_i|),\\
\\
\pi^{d_i^*}_l\left(\LL_{d_i^*}(q)\right)=l\pi^{d_i^*}_l\left(q\right)=0,\\
\\
\left|\pi^{d_i^*}_l\vc{0}{V_i^*q_1+R+f(q_1)}\right|\le  C^*\|q\|_{\H}^2 +C^*\d\sum_{j=1}^{\k-1} e^{-\frac 2{p-1}(\zeta^*_{j+1}-\zeta^*_j)} \le C^*(\|q\|_{\H}^2 +J_\k)
\end{array}
\]
where $J_\k$ is defined in \aref{defjbar}.

\medskip

- From (i) of Lemma \ref{lemboundk*} and the definition \aref{defbvi}, we see that $\bar V_i(y,s) = p(\lambda_i^{p-1}-1)\kappa({d^*_i},y)^{p-1}$ where $\lambda_i$ is defined right after \aref{defli}. Since we see from the definition \aref{defkd} and \aref{defWl2-0} of $\kappa(d_i^*,y)$ and $W_{l,2}(d_i^*,y)$ that $W_{l,2}(d_i^*,y)\le C\kappa(d_i^*,y)$, we write from the definition \aref{defpdi} of the projector $\pi_l^{d_i^*}$, \aref{hs} and (i) of Claim \ref{lemsobolev-9}, 
\begin{eqnarray*}
&&\left|\pi_l^{d_i^*}\vc{0}{\bar V_i q_1}\right|= \left|\iint W_{l,2}(d_i^*)\bar V_i q_1 \rho dy \right|\le C|\lambda_i^{p-1}-1|\cdot \left|\iint \kappa(d_i^*)^pq_1 \rho\right|\\
&\le & \frac{C|\nu_i|}{1-{d_i^*}^2}\left(\iint \kappa(d_i^*)^{p+1} \rho\right)^{\frac p{p+1}}\left(\iint q_1^{p+1} \rho\right)^{\frac 1{p+1}}\le \frac{C|\nu_i|}{1-{d_i^*}^2}\|q\|_{\H}.
\end{eqnarray*}

- Using Claim \ref{clnond'} for the projections of $\pnu \kappa^*(d_j, \nu_j)$ and $\pd \kappa^*(d_j, \nu_j)$, we write from equation \aref{eqq*} and the above estimate the following (we first project with $\pi_0^{d^*_i}$ then with $\pi_1^{d^*_i}$): 
\begin{eqnarray*}
\frac{|\nu_i'-\nu_i|+|d_i'|}{1-{d_i^*}^2}&\le & \frac {C\|q\|_{\H}}{1-{d_i^*}^2}
(|d_i'|+|\nu_i'-\nu_i|+|\nu_i|)+C^*\|q\|_{\H}^2+C^* J_\k\\
&&+C\left(\sum_{j\neq i}\frac{|\nu_j'-\nu_j|+|d_j'|}{1-{d^*_j}^2}\right)J_\k.
\end{eqnarray*}
Summing-up this estimate in $i$ and using the smallness of $\|q\|_{\H}$ and $J_\k$ (see \aref{small} and \aref{gap}),
we get the conclusion of (i).

\bigskip

{\bf Proof of (ii): A differential inequality on $\m(s)$}:

 Using equation \aref{eqq*}, we write
\begin{eqnarray}
\varphi(\ps q, q)&=& \varphi(\ll q,q)-\sum_{j=1}^\k (-1)^j[ (\nu_j-\nu_j') \varphi(\pnu \kappa^*(d_j, \nu_j),q)+d_j' \varphi(\pd \kappa^*(d_j, \nu_j),q)]\nonumber\\
&&+ \varphi((0,R),q)+\varphi((0, f(q_1)),q).\label{mohsen}
\end{eqnarray}
Using \aref{cle}, we see that the analysis performed in Appendix C of \cite{MZajm10} for the case where all the $\nu_i=0$ holds and gives
\begin{equation}\label{estR}
\begin{array}{ll}
\varphi(\ll q,q) = -\d\frac 4{p-1}\iint q_2^2\frac \rho{1-y^2} dy,&
\iint R^2 \rho(1-y^2) dy\le C J_\k,\\
\iint \kappa(d_i^*,y)|f(q_1)|\rho dy \le C\|q\|_{\H}^2,&
|R_-|\le C\|q\|_{\H}^{\bar p+1}
\end{array}
\end{equation}
where $\bar p = \min(p,2)$, $\F(q_1)$ and $J_\k$ are given in \aref{defF} and \aref{defjbar}.

- By definition \aref{defa-9} and \aref{defphi9} of $\m(s)$ and $\varphi$, we write
\begin{equation}\label{anas}
\frac 12 \m'(s) = \varphi(\ps q, q)- \frac{p(p-1)}2 \sum_{i=1}^\k (-1)^i (d_i' I_i +\nu_i' \bar I_i)
\end{equation}
where
\[
I_i= \iint \pd \kappa^*_1(d_i, \nu_i) |K^*_1|^{p-3}K^*_1 q_1^2 \rho dy\mbox{ and }\bar I_i = \iint \pnu \kappa^*_1(d_i, \nu_i) |K^*_1|^{p-3}K^*_1 q_1^2 \rho dy.
\]
Using 
\aref{bordev},
\aref{hs} and (ii) of Lemma \ref{lemtech} below, we see that
\begin{equation}\label{ghrab}
|I_i|+|\bar I_i|\le \frac C{1-{d_i^*}^2}\iint \kappa(d_i^*,y)|K^*_1|^{p-2}dy\|q_1(1-y^2)^{\frac 1{p-1}}\|_{L^\infty}^2\le \frac C{1-{d_i^*}^2} \|q\|_{\H}^2.
\end{equation}
Therefore, from \aref{anas}, (i) of Lemma \ref{lemproj*}, \aref{small} and \aref{gap}, we write
\begin{eqnarray}
\left|\frac 12 \m'(s) - \varphi(\ps q, q)\right|&\le& C\|q\|_{\H}^2 \left(C\sum_{i=1}^\k\frac{|\nu_i|}{1-{d_i^*}^2}+C\ep\|q\|_{\H} + CJ_\k\right)\nonumber\\
&\le &C\|q\|_{\H}^2 \left(\ep+C\sum_{i=1}^\k\frac{|\nu_i|}{1-{d_i^*}^2}\right).\label{amel}
\end{eqnarray}

- Since $|F_{1,l}({d^*},y)|\le C\kappa({d^*},y)$ by definitions \aref{deffld} and \aref{defkd}, using \aref{dnu}, \aref{defL1}, \aref{expdd}, \aref{defD}, \aref{boundD}, \aref{boundlnu} and \aref{condnu1}, we write
\begin{eqnarray}
&&(1-{d^*}^2)\left(\|\pnu \kappa^*(d,\nu)\|_{\H} +\|\pd \kappa^*(d,\nu)\|_{\H}\right)\nonumber\\
&\le & C +|\nu|\left\|\frac{(1-{d^*}^2)^{\frac 1{p-1}}}{(1+{d^*} y)^{\frac 2{p-1}+1}}\right\|_{L^2_\rho}\le C +\frac{C|\nu|}{\sqrt{1-{d^*}^2}}\le C.\label{h16A}
\end{eqnarray}
Therefore, 
using 
(i) of Lemma \ref{lemdecomp9}, (i) of Claim \ref{lemsobolev-9}, (i) of Lemma \ref{lemproj*} and \aref{apriori}, we write
\begin{eqnarray}
&&\left| (\nu_j-\nu_j') \varphi(\pnu \kappa^*(d_j, \nu_j),q)+d_j' \varphi(\pd \kappa^*(d_j, \nu_j),q)\right| \le \frac{C\|q\|_{\H}}{1-{d_j^*}^2}\left(|\nu_j'-\nu_j|+|d_j'|\right)\nonumber\\
&\le & C\|q\|_{\H} \left((\ep+\frac{|\nu_j|}{1-{d_j^*}^2})\|q\|_{\H} +J_\k\right)\le C (\ep+\frac{|\nu_j|}{1-{d_j^*}^2})\|q\|_{\H}^2 +C\|q\|_{\H} J_\k.\label{souissi}
\end{eqnarray}

- Using the definition \aref{defphi9} of $\varphi$ and \aref{estR}, we write
\begin{eqnarray}
\varphi((0,R),q) = \iint R q_2 \rho dy &\le & \frac 1{p-1} \iint q_2^2 \frac \rho{1-y^2} dy + C \iint R^2 \rho (1-y^2)\nonumber\\
&\le & \frac 1{p-1}  \iint q_2^2 \frac \rho{1-y^2} dy + C J_\k^2.\label{estRR}
\end{eqnarray}
Using the definition of $R_-$ given in \aref{estR}, we write from the definition \aref{defphi9} of $\varphi$ and the equation \aref{eqq*} satisfied by $q$,
\begin{eqnarray*}
&&R_-'+\varphi((0,f(q_1)),q) = R_-'+\iint q_2 f(q_1) \rho dy\\
&=& R_-'+\iint \ps q_1 f(q_1) \rho dy+\sum_{i=1}^\k (-1)^i d_i'\iint \pd \kappa^*_1(d_i, \nu_i) f(q_1) \rho dy\\
&& +(-1)^i (\nu_i'- \nu_i) \iint \pnu \kappa^*_1(d_i, \nu_i) f(q_1) \rho dy,\\
&=& \sum_{i=1}^\k d_i'\iint ((-1)^i \pd \kappa^*_1(d_i, \nu_i)f(q_1) - \partial_{d_i} \F(q_1)) \rho dy\\
&& +\nu_i'\iint ((-1)^i \pnu \kappa^*_1(d_i, \nu_i)f(q_1) - \partial_{\nu_i} \F(q_1)) \rho dy-(-1)^i \nu_i \iint\pnu \kappa^*_1(d_i, \nu_i) f(q_1) \rho dy\\
&=&\frac{p(p-1)}2 \sum_{i=1}^\k(-1)^i  d_i'\iint \pd \kappa^*_1(d_i, \nu_i)|K^*_1|^{p-3}K^*_1q_1^2\rho dy\\
&&+(-1)^i \nu_i'\iint \pnu \kappa^*_1(d_i, \nu_i) |K^*_1|^{p-3}K^*_1 q_1^2 \rho dy-(-1)^i \nu_i \iint \pnu \kappa^*_1(d_i, \nu_i) f(q_1) \rho dy.
\end{eqnarray*} 
Using 
\aref{bordev},
(ii) of Lemma \ref{lemtech},
\aref{estR}, (i) of Lemma \ref{lemproj*}, \aref{small} and \aref{gap}, we write
\begin{equation}\label{najet}
|R_-'+\varphi((0, f(q_1)), q)|\le C\|q\|_{\H}^2 \sum_{i=1}^\k \frac{|d_i'|+|\nu_i'|+|\nu_i|}{1-{d_i^*}^2}\le  C\|q\|_{\H}^2 \left(\ep+ \sum_{i=1}^\k \frac{|\nu_i|}{1-{d_i^*}^2}\right)
\end{equation}
Gathering the estimates \aref{mohsen}, \aref{estR}, \aref{amel}, \aref{souissi}, \aref{estRR} and \aref{najet},
we get to the conclusion of (ii) of Lemma \ref{lemproj*}.

\bigskip

{\bf Proof of (iii): An additional estimate}

Using equation \aref{eqq*}, we write
\begin{eqnarray*}
&&\frac d{ds} \iint q_1 q_2 \rho dy = \iint \ps q_1 q_2 \rho dy + \iint \ps q_2 q_1 \rho dy\\
&=& - \sum_{i=1}^\k(-1)^i\left(d_i'\iint \partial_{d_i} \kappa^*(d_i, \nu_i) \cdot (q_2,q_1)\rho dy + (\nu_i'- \nu_i)\iint \partial_{\nu_i} \kappa^*\cdot (q_2, q_1) \rho dy\right)\\
&&+\iint q_2^2 \rho dy + \iint q_1 (\L q_1 + \psi q_1 - \frac{p+3}{p-1} q_2 - 2y \py q_2 +f(q_1)+R)\rho dy
\end{eqnarray*}
where the dot ``$\cdot$'' stands for the usual inner product coordinate by coordinate. 

- Arguing as in pages 112 of \cite{MZjfa07},
we write
\[
\begin{array}{l}
\iint q_2^2 \rho dy \le \iint q_2^2 \frac \rho{1-y^2} dy,\;\;\iint q_1(\L q_1 + \psi q_1)\rho dy \le - \m + \iint q_2^2 \frac \rho{1-y^2} dy,\\
\\
\left|- \frac{p+3}{p-1} \iint q_1 q_2 \rho dy - 2 \iint q_1 \py q_2 \rho dy\right|\le \frac{\m}{100} +C^2 \iint q_2^2 \frac \rho{1-y^2} dy,\\
\\
\left|\iint q_1 f(q_1)\rho dy\right|\le C\|q\|_{\H}^{\bar p+1} \le \frac{\m}{100}\mbox{ where }\bar p = \min(p,2).
\end{array}
\]

- Using 
\aref{h16A},
the Cauchy-Schwartz inequality and (i) of Lemma \ref{lemproj*}, we write
\begin{eqnarray*}
&&\left|d_i'\iint \pd \kappa^*(d_i, \nu_i) \cdot (q_2, q_1) \rho dy 
+ (\nu_i'-\nu_i)  \iint \pnu \kappa^*(d_i, \nu_i) \cdot (q_2, q_1) \rho dy\right|\\ 
&\le & \left(|d_i'|\|\pd \kappa^*(d_i, \nu_i)\|_{\H} +(\nu_i'-\nu_i) \|\pnu \kappa^*(d_i, \nu_i)\|_{\H}\right)\|q\|_{\H}\\
&\le & \frac{C\|q\|_{\H}}{1-{d_i^*}^2}(|d_i'|+|\nu_i'- \nu_i|)
\le C\left(\left(\ep+\frac{|\nu_i|}{1-{d_i^*}^2}\right)\|q\|_{\H}^2 +J_\k \|q\|_{\H}\right)\\
&\le &\m\left(\frac 1{100} +\frac{C|\nu_i|}{1-{d_i^*}^2}\right)+CJ_\k^2\le \frac{\m}{50}+C J_m^2.
\end{eqnarray*}

 - Using \aref{hs} and \aref{estR}, we write
\begin{eqnarray*}
\iint q_1R\rho dy &\le& \left(\iint q_1^2 \frac \rho{1-y^2} dy \right)^{\frac 12}\left(\iint R^2 \rho (1-y^2) dy\right)^{\frac 12}\\
&\le& C\|q\|_{\H} \left(\iint R^2 \rho (1-y^2) dy\right)^{\frac 12}
\le \frac{\m}{10} + CJ_\k^2.
\end{eqnarray*}
Gathering the above estimates, we conclude the proof of (iii) in Lemma \ref{lemproj*}.

\bigskip

{\bf Proof of (iv): A rough estimate}

Using the definitions \aref{defnh0} and \aref{defPhi} of $\|q(s)\|_{\H}$ and the inner product $\phi$, we write from the equation \aref{eqq*} satisfied by $q$:
\begin{eqnarray}
\frac 12 \frac d{ds}\|q\|_{\H}^2&=&\phi(\ps q, q)\\
&=& \phi(\ll q,q)-\sum_{j=1}^\k (-1)^j[ (\nu_j-\nu_j') \phi(\pnu \kappa^*(d_j, \nu_j),q)+d_j' \phi(\pd \kappa^*(d_j, \nu_j),q)]\nonumber\\
&&+ \phi((0,R),q)+\phi((0, f(q_1)),q).\label{mohsen2}
\end{eqnarray}

- By definition \aref{defL9} and \aref{defPhi} of $\ll (q)$ and $\phi$, we write
\begin{eqnarray*}
\phi(\ll q,q) &=&\int(-\L q_1 +q_1) q_2 \rho dy +\int\left(\L q_1 +\psi q_1 -\frac{p+3}{p-1} q_2 -2 y \py q_2\right)\rho dy\\
&=& \int(1+\psi) q_1 q_2 \rho dy -\frac{p+3}{p-1}\int q_2^2\rho dy -2 \int y \py q_2\rho dy.
\end{eqnarray*}
Since
\[
|\psi(y,s)|\le C \sum_{j=1}^\k \lambda_j^{p-1} \le C e^{2(p-1)L_\k}
\]
by definition \aref{defpsi} and the bound \aref{defli}, and
\[
-2 \int y \py q_2\rho dy=\int q_2^2 \left(1-\frac{4y^2}{(p-1)(1-y^2)}\right) \rho dy \le \frac{p+3}{p-1}\int q_2\rho dy -\frac 4{p-1}\int q_2^2 \frac \rho{1-y^2} dy
\]
from integration by parts, we conclude that
\begin{equation}\label{Lqq}
\phi(\ll q,q) \le C\|q\|_{\H}^2 -\frac 4{p-1}\int q_2^2 \frac \rho{1-y^2} dy.
\end{equation}
Since $\phi((0, f(q_1)),q)=\varphi((0, f(q_1)),q)$ by definitions \aref{defPhi} and \aref{defphi9}, we argue as for \aref{souissi}, \aref{estRR} and \aref{najet} and write 
\begin{eqnarray}
&&\sum_{j=1}^\k (-1)^j[ (\nu_j-\nu_j') \phi(\pnu \kappa^*(d_j, \nu_j),q)+d_j' \phi(\pd \kappa^*(d_j, \nu_j),q)]\nonumber\\
&&+ \phi((0,R),q)+\phi((0, f(q_1)),q)\nonumber\\
&\le & C\|q\|_{\H}^2 +C J_m^2+\frac 1{p-1}\int q_2^2 \frac \rho{1-y^2} dy-R-'\label{reste}
\end{eqnarray}
where $R_-$ is defined in Lemma \ref{lemproj*}. Gathering the information from \aref{mohsen2}, \aref{Lqq} and \aref{reste}, we conclude the proof of Lemma \ref{lemproj*}.\Box
\subsection{A table for some integrals involving the solitons}\label{sub3}

Let us now recall the following result from \cite{MZajm10}:
\begin{lem}\label{lemtech}$ $\\
(i) Consider $\alpha>0$, $\beta>0$, $d_i=-\tanh \zeta_i$, $d_j=-\tanh \zeta_j$ and\\
 $I_1= \d\iint \kappa(d_j)^\alpha\kappa(d_i)^\beta(1-y^2)^{\frac{\alpha+\beta}{p-1}-1} dy$. As $|\zeta_i-\zeta_j|\to \infty$, the following holds:\\
if $\alpha=\beta$, then $I_1\sim C_0|\zeta_i-\zeta_j|e^{-\frac {2\beta}{p-1}|\zeta_i-\zeta_j|}$;\\
if $\alpha\neq \beta$, then $I_1\sim C_0e^{-\frac 2{p-1}\min(\alpha, \beta)|\zeta_i-\zeta_j|}$  for some $C_0=C_0(\alpha, \beta)>0$.\\
(ii)
For any $i=1,...,\k$, $\iint \kappa(d_i^*,y)|K^*_1|^{p-2} dy\le C$, where $K^*_1$ is defined in \aref{defld}. 
\end{lem}
{\it Proof}: See Lemma E.1 in \cite{MZajm10}.\Box

\def\cprime{$'$}

\noindent{\bf Address}:\\
Universit\'e de Cergy Pontoise, D\'epartement de math\'ematiques, 
2 avenue Adolphe Chauvin, BP 222, 95302 Cergy Pontoise cedex, France.\\
\vspace{-7mm}
\begin{verbatim}
e-mail: merle@math.u-cergy.fr
\end{verbatim}
Universit\'e Paris 13, Institut Galil\'ee, 
Laboratoire Analyse, G\'eom\'etrie et Applications, CNRS UMR 7539,
99 avenue J.B. Cl\'ement, 93430 Villetaneuse, France.\\
\vspace{-7mm}
\begin{verbatim}
e-mail: Hatem.Zaag@univ-paris13.fr
\end{verbatim}

\end{document}